\documentclass[11pt,a4paper,openany,titlepage]{memoir}
\usepackage{bm}

\usepackage{pgfplots}
\pgfplotsset{compat=1.15}
\usepackage{mathrsfs}
\usetikzlibrary{arrows}
\usepackage[OT1]{fontenc}

\usepackage[english,mt]{titlepagecustom} % Special IDSC styles and ommands  

\usepackage[latin1]{inputenc}

\usepackage{amsmath,amssymb,amsfonts,mathrsfs}

\usepackage{enumitem}

\usepackage{amsthm}%{ntheorem}

\usepackage{graphicx}

\usepackage{soul}

\usepackage{pdfpages}

\usepackage{graphicx}
\usepackage{float}

\usepackage{yfonts}
\usepackage{tikz-cd}

\usepackage{fontawesome}

\usepackage{utfsym}
\usepackage{tikz}
\usepackage{varioref}

\usepackage{mathtools}

\usepackage[h]{esvect}

\usepackage{array}

\usepackage{listings}
\usepackage{booktabs}

\nonzeroparskip
\defaultlists

\makeatletter

\if@twoside
  \copypagestyle{chapter}{Ruled}
\else
  \copypagestyle{chapter}{ruled}
\fi
\makeoddhead{chapter}{}{}{}
\makeevenhead{chapter}{}{}{}
\makeheadrule{chapter}{\textwidth}{0pt}
\copypagestyle{abstract}{empty}

\makechapterstyle{bianchimod}{%
  \chapterstyle{default}
  \renewcommand*{\chapnamefont}{\normalfont\Large\rmfamily}
  
  \renewcommand*{\printchaptername}{%
    \chapnamefont\centering\@chapapp}

  }

\setlrmarginsandblock{3.5cm}{3.5cm}{1}
\setulmarginsandblock{3.5cm}{3.5cm}{*}
\checkandfixthelayout 

\chapterstyle{bianchimod}

\setsecheadstyle{\Large\bfseries\rmfamily}
\setsubsecheadstyle{\large\bfseries\rmfamily}
\setsubsubsecheadstyle{\bfseries\rmfamily}
\setparaheadstyle{\normalsize\bfseries\rmfamily}
\setsubparaheadstyle{\normalsize\itshape\rmfamily}
\setsubparaindent{0pt}

\captionnamefont{\rmfamily\bfseries\footnotesize}
\captiontitlefont{\rmfamily\footnotesize}
\makeatother

\numberwithin{equation}{chapter}

\newtheorem{theorem}{Theorem}[chapter]
\newtheorem{example}[theorem]{Example}
\newtheorem{remark}[theorem]{Remark}
\newtheorem{corollary}[theorem]{Corollary}
\newtheorem{definition}[theorem]{Definition}
\newtheorem{propdef}[theorem]{Proposition-Definition}
\newtheorem{lemma}[theorem]{Lemma}
\newtheorem{proposition}[theorem]{Proposition} %proposition uses the same counter as theorem
\newtheorem{fact}{Fact}

\theoremstyle{plain} %nonumberbreak?

\newcommand{\C}{\mathbb{C}}

\newcommand{\Q}{\mathbb{Q}}

\newcommand{\Z}{\mathbb{Z}}

\newcommand{\frx}{\mathfrak{X}}

\newcommand{\codim}{\text{codim}(}

\newcommand{\pone}{\mathbb{P}^1}
\newcommand{\ptwo}{\mathbb{P}^2}
\newcommand{\pr}{\mathbb{P}^r}
\newcommand{\ponex}{\pone\times\pone}

\newcommand{\m}{M_{0,n}(\pr,d)}
\newcommand{\mbar}{\overline{M}_{0,n}(\pr,d)}

\newcommand{\monebar}{\overline{M}_{0,n}(\ponex,(d,e))}

\DeclarePairedDelimiter\abs{\lvert}{\rvert}

\renewcommand{\epsilon}{\ensuremath\varepsilon}

\renewcommand{\phi}{\ensuremath{\varphi}}

\let\emptyset\varnothing

\usepackage{amssymb}

\usepackage{diagbox}
\usepackage{hhline}
\usepackage{boldline}
\usepackage{makecell}

\studentA{Greg Weiler}
\title{Kontsevich's Formula for Rational Curves from Classical and Quantum Perspectives} 
\supervision{Prof. Dr. Rahul Pandharipande\\ Prof. Dr. Junliang Shen}
\date{March 2023}

\begin{document}

\maketitle

\frontmatter

%% The abstract of your thesis.  Edit the file as needed.
\begin{abstract}
    Kontsevich's formula for rational plane curves is a recursive relation for the number $N_d$ of degree $d$ rational curves in $\ptwo$ passing through $3d-1$ general points. 
    We provide two proofs of this recursion: the first more direct and combinatoric, the second more abstract. 
    In order to achieve this, we introduce several moduli spaces, such as the Deligne-Mumford-Knudsen spaces and the Kontsevich spaces, and exploit their properties. In particular, the boundary structure of these spaces gives rise to certain fundamental relations crucial to both proofs.
    For the second proof, we reconsider the objects in question from the cohomological viewpoint and generalize the numbers $N_d$ to Gromov-Witten invariants. We introduce quantum cohomology and deduce Kontsevich's formula from the associativity of the quantum product.
    We also adapt these steps to the case of curves in $\ponex$, whose bidegrees lead to slightly more complicated but analogous results.
\end{abstract}

%% TOC with the proper setup, do not change.
\cleartorecto
\tableofcontents
\mainmatter

%% Your real content!

% Some commands used in this file
\newcommand{\package}{\emph}

%for minzi
%reread intros of book and of FP notes

\chapter{Introduction}

%p1xp1 simplest rational surface beyond p2

%explain main idea we exploit, mabe with picture like junliang.
%explain two possible "compactifications" of p1
%https://math.stackexchange.com/questions/983568/how-to-show-p1-times-p1-as-projective-variety-by-segre-embedding-is-not-is

%mention combinatorics

%what is a moduli problem
%projective spaces over base field \C
%reference Euclid, Whitehead's axioms 

%pose the moduli problem with historical perspectives on euclid and on the axioms of projective geometry

%say that we base pr stuff on invitiation

%Quantum cohomology first came up in the context of theoretical physics, see \cite{Witten:1990hr}.
%\textcolor{red}{explain some more}

The faculty of counting is one of the most essential components of cognition. From enumeration, the notions of units and measurements can flourish. It may then come as little surprise that enumerative geometry is an ancient science. Euclid's first axiom for planar geometry concerned the existence of a line between two arbitrary points. About two millenia later, Whitehead published his axioms for projective geometry, the second of which postulates a unique line through any two given points. Throughout the vast history of enumerative geometry, there has been no shortage of innovation. A myriad of minds developed entirely new mathematical disciplines whose techniques found application in the old art.

A more recent such development was the advent of the moduli space. Algebraic geometers, seeing that the objects of their interest often appeared in algebraically describable families, decided to apply their methods not only to the objects in question, but also to the spaces formed by these families. It became possible to rigorously formulate the quest for the nature of these families in so called moduli problems. In lucky cases, various types of moduli spaces corresponding to the problem at hand could be constructed, their properties studied and exploited. We will encounter one such case in this thesis -- our moduli problem will be that of the classification of rational curves of a given degree passing through finitely many given points in a complex projective space. We have already discussed the easiest relevant example of such a curve: the unique line passing through two general points. But much more can be done. Given any degree, we will determine the necessary number of point constraints and compute the number of rational curves of that degree passing through those points. This result is known as Kontsevich's formula. The main spaces in which we work will be $\ptwo$, the projective plane, and $\ponex$, the product of two projective lines. Although $\ponex$ is the simplest rational surface (other than $\ptwo$ itself) and the approach is similar, we shall encounter some significant differences between the two cases. 

Our approach will be twofold: at first, we will directly use the structure, in particular the boundary structure, of our moduli spaces. We will count by hand the number of points in certain intersections. 
Later, we shall introduce more sophisticated machinery. Combining the knowledge of our moduli spaces with notions of cohomological algebra will lead, via the introduction of Gromov-Witten invariants, to the formulation of quantum cohomology. These concepts have their origin in quantum field theory, hence the terminology. Ultimately, this slightly more abstract formalism will permit us to proceed with far greater elegance than before.
As one might expect, combinatorial arguments will be indispensable at every stage.

Our main source, especially for the theory sections concerning $\pr$, will be \cite{invitation}. This in turn serves as a guide to the rigorous constructions in \cite{FPnotes}.

\chapter{Primer on Moduli Spaces}\label{th0} 
%inspired heavily by \cite{invitation}

%should include def family, fine and coarse moduli space and an example. maybe early exercise with yoneda i did at steep. also intuition on when to expect fine or coarse. provide reference for rest. 

%For an introduction to moduli spaces, see chapters 0 and 1 of \cite{invitation}.

%green stuff from chapter 0
As mentioned in the Introduction, moduli spaces are essential tools in enumerative geometry. They parametrize families of the objects of interest, and the study of their structure can have powerful implications for the moduli problem in question -- this will be a guiding principle throughout this thesis.
We briefly introduce the most important theoretical aspects. From now on, we always work over the field $\C$.

    \begin{propdef} %family, pullback (0.2.2), sections (0.1.2)
       A family $\mathfrak{X}$ over a base variety $B$ is a morphism $\mathfrak{X}\rightarrow B$ (with some extra structure depending on the context, e.g. flatness). We concisely write $\mathfrak{X}/B$. Equivalence of families will be defined case by case. Usually, selecting the point space $\text{Spec}(\C)=\bullet$ as base variety $B$ will yield information on the objects within the families.

       Given another variety $B'$ and a morphism $\phi:B'\rightarrow B$, we define the pullback family $\phi^*\mathfrak{X}/B'$ as the following fiber product:
       \begin{center}
           \begin{tikzcd}
                B'\times_B \mathfrak{X} \arrow[rr] \arrow[dd,swap,"\pi'"] && \mathfrak{X} \arrow[dd,"\pi"]\\\\
                B' \arrow[rr,"\phi"] && B 
           \end{tikzcd}
       \end{center}
       This pullback operation is well behaved in regards to the identity $B\rightarrow B$, to compositions $ \psi\circ\phi$ and to equivalence of families.
    \end{propdef}
        
    \begin{definition}%universal family and fine moduli space: intuitive then categorical
        A universal family for a moduli problem is a family $U/M$ such that any other family $\mathfrak{X}/B$ can be obtained from it via pullback. More precisely, this means that given any family $\mathfrak{X}/B$, there exists a unique morphism $\kappa:B\rightarrow M$, called a classifying map, such that $\kappa^* U/B\cong\frx/B$. The base family $M$ of the universal family is called a fine moduli space.
    \end{definition}

    \begin{remark}[Parameter Space]%correspondence btw objs and pts of M
        The fine moduli space $M$ parametrizes the objects of the moduli problem: by the above definition, given a base variety $B$, there is a one to one correspondence between families over $B$ and morphisms $\kappa:B\rightarrow M$. In the special case $B=\bullet$, families over $\bullet$ correspond to objects and morphisms $\kappa:\bullet\rightarrow M$ correspond to points of $M$.
    \end{remark}

    \begin{remark}\label{nonexmod} %non-existence of moduli
        We illustrate a kind of moduli problem for which no fine moduli space can exist. Assume that there exists a family $\frx/B$ over an irreducible variety $B$ such that all members except one belong to a single equivalence class, the other member making up its own equivalence class. 
        Assume that there exists a coarse moduli space $(M,v)$. Then there is a classifying morphism $B\rightarrow M$. %\textcolor{red}{why in this direction? clear for fine}. 
        But this morphism is not continuous due to the presence of the special point, which is a contradiction since $B$ is an irreducible variety. Thus there can be no such moduli space.

        A concrete example of such a moduli problem is that of the classification of two points in $\pone$ up to linear equivalence. %x_0 x P1 -> P1
    \end{remark}

    We now rephrase the definition of fine moduli spaces in categorical terms.
    First, we encode the moduli problem in the so called moduli functor, a contravariant functor from the category of schemes $\mathbf{Sch}$ to the category of sets $\mathbf{Set}$. It send scheme $B$ to the set of equivalence classes of families over $B$, and a morphism $\phi:B'\rightarrow B$ to the pullback $\phi^*$.

    Recall that given a scheme $Y$, there exists a contravariant functor $h_Y:\mathbf{Sch}\rightarrow\mathbf{Set}$, the so called functor of points of $Y$. It sends a scheme $B$ to the set $\text{Hom}(B,Y)$ and $\phi:B'\rightarrow B$ to the precomposition with $\phi$.

    A functor $F$ is called representable if there is an isomorphism of functors $u:h_Y\rightarrow F$ for some scheme $Y$. In this case we say that the functor $F$ is represented by the pair $(Y,u)$. By the Yoneda Lemma, a natural transformation $u:h_Y\rightarrow F$ corresponds to an element $U\in F(Y)$. We therefore also say that $F$ is represented by $(Y,U)$.
   
    \begin{definition}\label{finemodcat}%fine mod space categ def
        Let $F$ be a moduli functor. Then we call a family $U/M$ universal and its base $M$ a fine moduli space for $F$ if $(M,U)$ represents $F$.
    \end{definition}
    
    If we weaken this definition by allowing non-representable functors, we arrive at the notion of a coarse moduli space.
    \begin{definition}%coarse moduli space: categorical
        A coarse moduli space for a moduli functor $F$ is a pair $(M,v)$ with $M\in\mathbf{Sch}$ and $v:F\rightarrow h_M$ is a natural transformation such that the following conditions are met:
        \begin{enumerate}[label=(\roman*)]
            \item $(M,v)$ is initial among such pairs
            \item The map of sets $v_\bullet:F(\bullet)\rightarrow\text{Hom}(\bullet,M)$ is bijective,
        \end{enumerate}
        where $\bullet=\text{Spec}(\C)$.
    \end{definition}
    By the Yoneda Lemma, there is a one to one correspondence between natural transformations $h_M\rightarrow h_{M'}$ and morphisms $M\rightarrow M'$. Using the same notation for both, the first condition then means that given another pair $(M',v')$, there exists a unique morphism $\psi:M\rightarrow M'$ such that $v'=\psi\circ v$. 
    Intuitively, every natural transformation $v':F\rightarrow h_{M'}$ factors through $v:F\rightarrow h_M$ and thus $h_M$ is the representable functor closest to $F$.

    The second condition means that as before, geometric points of $M$ parametrize equivalence classes of objects.

    Note that fine moduli spaces are examples of coarse moduli spaces: just take $v:=u^{-1}$, where $u:h_M\rightarrow F$ is the equivalence corresponding to $U$ in Definition \ref{finemodcat}.
   
    \begin{example}[Universal Family] %example: exc 14 or 15?
        Consider the moduli problem of classifying the subsets of a scheme: a family over a scheme $B$ is defined to be a subset of $B$, i.e. $A/B$ is given by the inclusion morphism $A\xhookrightarrow{} B$. Two subsets are defined to be equivalent if they are equal, i.e. if they contain the same elements.

        The universal family for this problem is given by the inclusion $\{1\}\xhookrightarrow{} \{1,0\}$, where the values $1$ and $0$ can be interpreted as 'true' and 'false', respectively. In particular, the set $\{1,0\}$ is a fine moduli space for this moduli problem.

        To show this, one must prove that the pair $(\{1,0\},\{1\})$ represents the moduli functor $\mathcal{P}:\mathbf{Sch}^\text{op}\rightarrow\mathbf{Set}$ sending a scheme $B$ to its power set $\mathcal{P}(B)$ and a morphism to the pullback via this morphism. More precisely, it sends $\alpha:A\rightarrow B$ to $\alpha^*:\mathcal{P}(B)\rightarrow \mathcal{P}(A)$, with $\alpha^{*}(V)=\alpha^{-1}(V)\subset A$ for a subset $V\subset B$.

        By the Yoneda Lemma, this is equivalent to the existence of an isomorphism of functors $p:h_{\{1,0\}}\rightarrow \mathcal{P}$. 

        Let us define such a natural transformation $p$ objectwise and then show that it is an isomorphism.
        Given a scheme $B$, define $p_B:\text{Hom}(B,\{1,0\})\rightarrow \mathcal{P}(B)$ as $p_B(f):=f^{-1}(\{1\})$.
        In the opposite direction, given $V\in\mathcal{P}(B)$, define $\chi_B(V)$ on $x\in B$ as follows:
        $$\chi_B(V)(x)=\begin{cases}1, \text{ if }x\in V,\\
        0, \text{ else.}
        \end{cases}$$
        Then $\chi_B\circ p_B(f)=\chi_B(f^{-1}(\{1\}))=f$ and $p_B\circ\chi_B(V)=\chi_B^{-1}(\{1\})=V$.
        Since $B\in\mathbf{Sch}$ was chosen arbitrarily, this concludes the argument.
        
    \end{example}

    \begin{remark}[Intuition on the Existence of Moduli Spaces]\label{intuition}
        %intuition on when to expect fine or coarse. provide reference for rest. p18 inv notes
        %every object is auto-free -> expect fine mod spcae
        %every obj has finite aut grp -> coarse mod sp
        %there is an obj with infinite aut grp -> not even coarse mod sp
        The existence of moduli spaces is correlated with the number of automorphisms of the objects in question (however, the following are not formal obstructions, and exceptions do exist). 
        
        If every object is automorphism-free (i.e. admits only the trivial automorphism), then one can expect the existence of a fine moduli space.
        If every object has a finite group of automorphisms, one can expect a coarse moduli space.
        However, if there exists an object with infinitely many automorphisms, not even a coarse moduli space is to be expected. Note that this principle applies to the concrete example in Remark \ref{nonexmod}.
        %give pts in pone exp?
    \end{remark}
\chapter{Relevant Spaces and their Properties}\label{th1} %Classical Theory for Kontsevich Formula
%construction (sketched as in book) of relevant moduli spaces, their dimensions, boundary properties, ...
%trees, stable maps, marks, transversality...

%MUST explain bonudary (divisors). what is being counted? how do irreducible components, divisors and cardinalities relate? DEFINITION OF DEGREE VIA INTERSECTION WITH GENERAL HYPERPLANE, EQUIVALENT DIVISORS HAVE SAME DEGREE?
%explain link between genus 0 and rationality. should be clear.

%green stuff from chapters 1,2

%say in tex that we omit proofs
%mention that this is all abbreviated
%mention the names of the spaces: deligne mumford knudsen, Kontsevich
    In this section, we describe the concrete moduli spaces needed to tackle our enumerative problems. We begin by considering only curves, and later introduce the additional structure of maps. We modify the corresponding spaces until we find one whose properties allow us to deduce the results we desire: the so called Kontsevich space.
    
    For proofs and further details see again \cite{invitation}. For rigorous constructions see \cite{FPnotes}.

     \section{Moduli Spaces of Curves and Maps for $\pr$}\label{th1pr}
        %define W(d),Wdeg,W* and mention 3 problems

        \begin{definition}\label{nptsmooth}%n-pted smooth rational curve, family thereof
            An $n$-pointed smooth rational curve $(C,p_1,\dots,p_n)$ is a projective smooth rational curve $C$ together with a selection of $n$ marked points (or simply, marks) $p_1,\dots,p_n\in C$. 
            Two such curves $(C,p_1,\dots,p_n)$ and $(C',p_1',\dots,p_n')$ are isomorphic if there exists an isomorphism $\phi:C\rightarrow C'$ with $\phi(p_i)=p_i'$ for $1\le i\le n$.

            More generally, a family $\pi:\mathfrak{X}\rightarrow B$ of $n$-pointed smooth rational curves is given by:
            
            \begin{center}
                \begin{tikzcd}%[row sep=huge]
                    \mathfrak{X} \arrow[dd,swap,shift right=2ex,"\pi"] \\\\
                     B \arrow[uu,swap,shift right=-0.5ex,"\dots",shorten >=1ex, shorten <=1ex] 
                    \arrow[uu,swap,shift right=2.25ex,"\sigma_i",shorten >=1ex, shorten <=1ex]
                \end{tikzcd}
            \end{center}
            where the map $\pi:B\rightarrow\mathfrak{X}$ is flat and proper and the maps $\sigma_i$ for $1\le i\le n$ are disjoint sections (i.e. $\pi\circ\sigma_i=\text{id}_B$) such that $\mathfrak{X}_b:=\pi^{-1}(b)$ is a projective smooth rational curve for all $b\in B$.
            Two such families $\pi:B\rightarrow\mathfrak{X}$ and $\pi':B\rightarrow\mathfrak{X}'$ with the same base are isomorphic if there is an isomorphism $\phi:\mathfrak{X}\rightarrow\mathfrak{X}'$ making the following diagram commute for every $i$:
            \begin{center}
                \begin{tikzcd}[column sep=1.5em]
                    \mathfrak{X} \arrow[dr,swap,shift right=.75ex,"\pi"] \arrow{rr}{\phi} && \mathfrak{X}'\arrow[dl,swap,shift right=.75ex,"\pi'"] \\
                    & B \arrow[ul,swap,shift right=.75ex,"\sigma_i",shorten >=1ex, shorten <=1ex] \arrow[ur,swap,shift right=.75ex,"\sigma_i'",shorten >=1ex, shorten <=1ex]
                \end{tikzcd}%add point here?
            \end{center}
        \end{definition}

        \begin{proposition}%existense of M0n fine mod space for iso classes of curves above (construct universal family)
            For $n\ge3$, there exists a fine moduli space $M_{0,n}$ for the classification problem of $n$-pointed smooth rational curves up to isomorphism.
        \end{proposition}
        The subscript $0$ means that we work with rational curves, i.e. curves of genus $0$.
        The restriction on $n$ stems from the fact that classifying $n$-pointed smooth rational curves up to isomorphism is equivalent to the classification of $n$-tuples of points in a fixed $\pone$ up to projective equivalence, and any triple of points in $\pone$ can be mapped to $(0,1,\infty)$ by a unique automorphism (see Chapter 0 of \cite{invitation}). In particular, the space $M_{0,3}$ contains only a single point and $M_{0,4}\cong\pone$.
        The universal family looks as follows:
        \begin{center}
            \begin{tikzcd}[row sep=huge]
                M_{0,n}\times\pone \arrow[d,swap,shift right=2ex,"\pi"] \\
                  M_{0,n} \arrow[u,swap,shift right=-0.5ex,"\dots",shorten >=1ex, shorten <=1ex] 
                \arrow[u,swap,shift right=2.25ex,"\sigma_i",shorten >=1ex, shorten <=1ex]
            \end{tikzcd}
        \end{center}
        with $n$ sections, three of which are the constant sections $0$, $1$, and $\infty$. For more details see section 1.1 in \cite{invitation}.

        \begin{definition}%tree of projective lines, twig, node
            A tree of projective lines is a connected curve satisfying the following conditions:
            \begin{itemize}
                \item each irreducible component (called a twig) is isomorphic to $\pone$,
                \item the points of intersection (called nodes) of the components are ordinary double points, %each twig diff tngt direction
                \item there are no loops: if $\delta$ is the number of nodes, then there are $\delta+1$ irreducible components.
            \end{itemize}
        \end{definition}
        Note that this is not a tree in the graph theoretical sense, but that the dual graph is! The three conditions amount to the fact that the curve has arithmetic genus $0$.

        \begin{definition}%stable n-pointed rational curve
            For $n\ge3$, we call $n$-pointed rational curve a tree $C$ with $n$ distinct marked points that are smooth points of $C$. Isomorphisms between $n$-pointed rational curves are defined as in the smooth case, see Definition \ref{nptsmooth}.

            An $n$-pointed rational curve is called stable if each of its twigs contains at least three special points, i.e. points that are nodes or marks. 

            A family of stable $n$-pointed rational curves is a flat and proper morphism $\pi:\mathfrak{X}\rightarrow B$ with $n$ distinct sections such that for all $b\in B$, the geometric fiber $\mathfrak{X}_b:=\pi^{-1}(b)$ is a stable $n$-pointed rational curve. Isomorphisms are defined as in the smooth case.
        \end{definition}

        \begin{remark}
            Note that the stability condition for $n$-pointed rational curves is equivalent to being automorphism-free. %1.2.4
        \end{remark}       
        We shall from now on often simply refer to curves, but this will always mean rational curves.

        \begin{theorem}[Knudsen]%existence of M0nbar fine mod space for curves above, del-mum-knud space
            For each $n\ge3$, there is a smooth projective variety $\overline{M}_{0,n}$ that is a fine moduli space for stable $n$-pointed curves. Furthermore, it contains the subvariety ${M}_{0,n}$ as a dense open set.
        \end{theorem}
        
        This moduli space is called the Deligne-Mumford-Knudsen space. The universal family $\overline{U}_{0,n}\rightarrow\overline{M}_{0,n}$ can be described recursively, in the sense that it is given by the morphism $\overline{M}_{0,n+1}\rightarrow\overline{M}_{0,n}$ forgetting the last mark. We will now discuss this morphism in more detail. 
        
        \begin{remark}%forgetful map and stabilization
            Given a stable $n$-pointed curve $(C,p_1,\dots,p_n)$ and a point $q\in C$, there is a canonical way of producing a stable $n+1$-pointed curve (called stabilization). The crucial part is to deal with the cases where $q$ coincides with an existing mark or node by adding certain twigs. 

            Similarly, given a stable $n+1$-pointed curve $(C,p_1,\dots,p_{n+1})$ for $n\ge3$, there is a canonical way of associating to it a stable $n$-pointed curve $(C,p_1,\dots,p_n)$. It consists in forgetting the last mark $p_{n+1}$. If $C$ is non-reduced, one must take care to contract non-stable twigs that may arise as a result of the removal of $p_{n+1}$. 

            Both of these processes can be generalized to families, see Section 1.3 of \cite{invitation} for further details.
        \end{remark}
        
        %bdry, bdry cycle, bdry divisor
        Consider the boundary $\overline{M}_{0,n}\setminus{M}_{0,n}$ of $\overline{M}_{0,n}$. Note that each of its points corresponds to a non-reduced curve. The boundary admits a stratification into the subsets $\Sigma_\delta$ of curves with $\delta\le n-3$ nodes. %why this constraint? stability?
        The dimension of $\Sigma_\delta$ is equal to $n+2\delta-3(\delta+1)=n-3-\delta$. To see this, sum up the possible ways of moving the special points on each twig: the number $n+2\delta$ then represents the number of special points (each node counting towards two twigs), and the number $3(\delta+1)$ is subtracted since each of the $\delta+1$ twigs is isomorphic to $\pone$, so we can map three of its special points to $(0,1,\infty)$ via a unique automorphism. This way of computing dimensions is justified in Remark \ref{recstruc}.%1.5.9 book

        \begin{example}%m06bar: draw pictures (all of them at once in geogebra!)
        In Figure \ref{expstrat}, we illustrate this stratification for $\overline{M}_{0,6}$. The rows in the figure correspond, in descending order, to the strata $\Sigma_\delta$ with $\delta=0,1,2,3$ and thus $\dim(\Sigma_\delta)=3,2,1,0$.
        Note that we do not label the marks; in fact each depicted configuration corresponds to several cycles (indeed as many as there are ways of distributing the marks). The first unlabeled configuration in the second row, for instance, corresponds to $\frac{6!}{2!4!}=15$ labeled configurations. 
        \newpage
    
        \definecolor{xdxdff}{rgb}{0.49019607843137253,0.49019607843137253,1}
        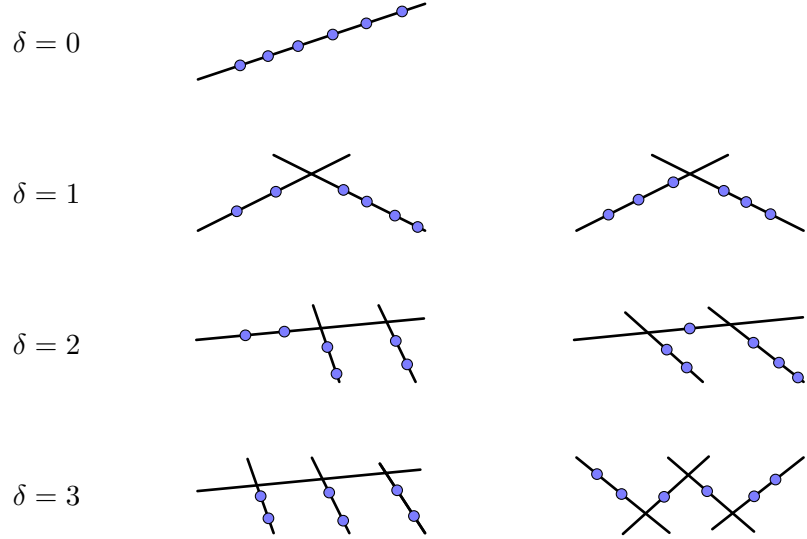
\begin{figure}
            \centering
            \caption{Stratification of the boundary of $\overline{M}_{0,6}$}
            \label{expstrat}
        
        \begin{tikzpicture}[line cap=round,line join=round,>=triangle 45,x=1cm,y=1cm]
                \clip(-8,-4) rectangle (9.3197922077922,6);
                \draw [line width=1pt] (-1,5)-- (-4,4);
                \draw[color=black] (-6,4.5) node {$\delta=0$};
                \draw [line width=1pt] (-4,2)-- (-2,3);
                \draw[color=black] (-6,2.5) node {$\delta=1$};
                \draw [line width=1pt] (-3,3)-- (-1,2);
                \draw [line width=1pt] (1,2)-- (3,3);
                \draw [line width=1pt] (2,3)-- (4,2);
                \draw[color=black] (-6,0.5) node {$\delta=2$};
                \draw [line width=1pt] (-4.021585202676114,0.5550609996064534)-- (-1.014894923258561,0.838414010232191);
                \draw [line width=1pt] (-1.6130846123573412,1.011574183392364)-- (-1.1250877607241256,0);
                \draw [line width=1pt] (-2.4788854781582073,1.011574183392364)-- (-2.1325651318378607,0);
                \draw [line width=1pt] (2.747403384494293,0.9800905155450598)-- (4,0);%plus 1.252596615505707 minus 0.9800905155450598
                \draw [line width=1pt] (1.6454750098386457,0.9171231798504513)-- (2.668694214876033,0);
                \draw [line width=1pt] (3.9910082644628098,0.854155844155843)-- (0.9685761511216049,0.5550609996064534);
                \draw[color=black] (-6,-1.5) node {$\delta=3$};
                \draw [line width=1pt] (-4.005843368752464,-1.444151908697362)-- (-1.062120425029519,-1.1607988980716244);
                \draw [line width=1pt] (-3.344686343959075,-1.0191223927587556)-- (-3,-2);
                \draw [line width=1pt] (-2.494627312081861,-1.0033805588351035)-- (-2,-2);
                \draw [line width=1pt] (-1.5973427784336909,-1.082089728453364)-- (-1,-2);
                \draw [line width=1pt] (-1,-2)-- (-1.5973427784336909,-1.082089728453364);
                \draw [line width=1pt] (2.747403384494288,-0.9876387249114514)-- (1.613991341991336,-2.010857929948837);
                \draw [line width=1pt] (2.212181031090116,-1.0033805588351035)-- (3.345593073593068,-1.979374262101533);
                \draw [line width=1pt] (1,-1)-- (2.227922865013768,-1.995116096025185);
                \draw [line width=1pt] (2.794628886265244,-1.9478905942542288)-- (4,-1);
                \begin{scriptsize}
                    \draw [fill=xdxdff] (-3.442518772136965,4.185827075954345) circle (2pt);
                    \draw [fill=xdxdff] (-3.074159858323506,4.308613380558832) circle (2pt);
                    \draw [fill=xdxdff] (-2.677465643447473,4.440844785517509) circle (2pt);
                    \draw [fill=xdxdff] (-2.2193782762691963,4.593540574576934) circle (2pt);
                    \draw [fill=xdxdff] (-1.7754585596222072,4.741513813459264) circle (2pt);
                    \draw [fill=xdxdff] (-1.3032035419126435,4.898932152695785) circle (2pt);
                    \draw [fill=xdxdff] (-3.4870186540732093,2.2564906729633956) circle (2pt);
                    \draw [fill=xdxdff] (-2.97068650137742,2.51465674931129) circle (2pt);
                    \draw [fill=xdxdff] (-2.0773878000787165,2.538693900039358) circle (2pt);
                    \draw [fill=xdxdff] (-1.0973878000787165,2.048693900039358) circle (2pt);
                    \draw [fill=xdxdff] (-1.7688478551751352,2.3844239275875676) circle (2pt);
                    \draw [fill=xdxdff] (-1.398450373868563,2.1992251869342814) circle (2pt);
                    \draw [fill=xdxdff] (1.4232528925619747,2.211626446280987) circle (2pt);
                    \draw [fill=xdxdff] (1.8199471074380078,2.409973553719004) circle (2pt);
                    \draw [fill=xdxdff] (2.279608658008649,2.6398043290043245) circle (2pt);
                    \draw [fill=xdxdff] (2.946224950806763,2.5268875245966185) circle (2pt);
                    \draw [fill=xdxdff] (3.2421714285714223,2.3789142857142886) circle (2pt);
                    \draw [fill=xdxdff] (3.5633048406139256,2.218347579693037) circle (2pt);
                    \draw [fill=xdxdff] (-3.373028916076586,0.6161814873488172) circle (2pt);
                    \draw [fill=xdxdff] (-2.8581214624379814,0.664706797115806) circle (2pt);
                    \draw [fill=xdxdff] (-2.2902297108973366,0.4605256939321214) circle (2pt);
                    \draw [fill=xdxdff] (-2.16937019897,0.10750467341718328) circle (2pt);
                    \draw [fill=xdxdff] (-1.3864497284734534,0.5417801738085065) circle (2pt);
                    \draw [fill=xdxdff] (-1.2367115633767463,0.23138624078735637) circle (2pt);
                    \draw [fill=xdxdff] (2.4961527201048446,0.7062274309120863) circle (2pt);
                    \draw [fill=xdxdff] (2.192632141422442,0.4266999294603773) circle (2pt);
                    \draw [fill=xdxdff] (2.456576501075185,0.19012355439178197) circle (2pt);
                    \draw [fill=xdxdff] (3.337868956986818,0.5180824754529452) circle (2pt);
                    \draw [fill=xdxdff] (3.6749157817040152,0.2543611847271696) circle (2pt);
                    \draw [fill=xdxdff] (3.9254351048051564,0.05834308161815784) circle (2pt);%plus 1.252596615505707 minus 0.9800905155450598
                    \draw [fill=xdxdff] (-3.1715538892088126,-1.511807266782669) circle (2pt);
                    \draw [fill=xdxdff] (-3.0688173490282913,-1.8041657355025054) circle (2pt);
                    \draw [fill=xdxdff] (-2.266168537404312,-1.4636997745088909) circle (2pt);
                    \draw [fill=xdxdff] (-2.0815440807892687,-1.8356976773392242) circle (2pt);
                    \draw [fill=xdxdff] (-1.368770459180322,-1.4333267186854213) circle (2pt);
                    \draw [fill=xdxdff] (-1.1470513348296008,-1.7740327069180133) circle (2pt);
                    \draw [fill=xdxdff] (1.2710187719241695,-1.2196352474173495) circle (2pt);
                    \draw [fill=xdxdff] (1.5960502941107266,-1.4830427534252968) circle (2pt);
                    \draw [fill=xdxdff] (2.725907446455874,-1.4457560831778389) circle (2pt);
                    \draw [fill=xdxdff] (2.1577496542305012,-1.5199650091773693) circle (2pt);
                    \draw [fill=xdxdff] (3.629326601546048,-1.291493486056833) circle (2pt);
                    \draw [fill=xdxdff] (3.347159202770343,-1.5133868268354709) circle (2pt);
                \end{scriptsize}
            \end{tikzpicture}
            \end{figure}
        \end{example}

        \begin{definition}%bdry cycle, bdry divisor
            For $\delta\ge1$, we call the closure of an element of $\Sigma_\delta$ a boundary cycle. It consists of the element in question along with boundary cycles corresponding to further degenerated elements from strata with more nodes. A codimension $1$ boundary cycle is called a boundary divisor. 
        \end{definition}
        %describe dab
        If $S=\{p_1,\dots,p_n\}$ is the set of marks, then for each partition $S=A\sqcup B$ with $\abs{A},\abs{B}\ge2$ there is a boundary divisor $D(A|B)$, with each of its points representing a curve $C$ with two twigs, $C_A$ and $C_B$, respectively carrying the marks from $A$ and $B$.

        \begin{proposition}%normal crossings
            The boundary of $\overline{M}_{0,n}$ is a divisor with normal crossings, i.e. irreducible components of $k$-fold intersections of divisors have codimension $k$.
        \end{proposition}

        \begin{example}%m05
            We illustrate the intersection of the boundary divisors $D(p_1p_2|p_3p_4p_5p_6)$ and $D(p_1p_2p_3|p_4p_5p_6)$ in $\overline{M}_{0,6}$. The only common degenerations of these divisors lie in the closure of the following codimension $2$ configuration:
            
            \definecolor{xdxdff}{rgb}{0.49019607843137253,0.49019607843137253,1}
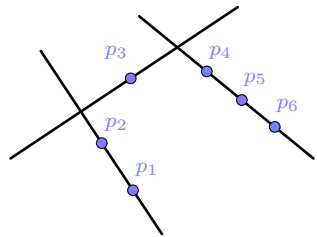
            \begin{tikzpicture}[line cap=round,line join=round,>=triangle 45,x=1cm,y=1cm]
                \clip(-5,-1.25) rectangle (15.405350303920823,2.5);
                \draw [line width=1pt] (0,0)-- (3,2);
                \draw [line width=1pt] (1.7,1.88)-- (4,0);
                \draw [line width=1pt] (0.4,1.42)-- (2,-1);
                \begin{scriptsize}
                    \draw [fill=xdxdff] (1.5876923076923077,1.0584615384615386) circle (2pt);
                    \draw[color=xdxdff] (1.4,1.4) node {$p_3$};
                    \draw [fill=xdxdff] (1.6175847155553442,-0.4215968822774583) circle (2pt);
                    \draw[color=xdxdff] (1.805591645368366,-0.15) node {$p_1$};
                    \draw [fill=xdxdff] (1.204752625825769,0.20281165343852448) circle (2pt);
                    \draw[color=xdxdff] (1.4,0.45) node {$p_2$};
                    \draw [fill=xdxdff] (2.592330356738135,1.1506169257966548) circle (2pt);
                    \draw[color=xdxdff] (2.775828357405361,1.4) node {$p_4$};
                    \draw [fill=xdxdff] (3.0551244277231313,0.772333076469788) circle (2pt);
                    \draw[color=xdxdff] (3.2362796783720706,1.05) node {$p_5$};
                    \draw [fill=xdxdff] (3.490082045238203,0.41680250215312076) circle (2pt);
                    \draw[color=xdxdff] (3.6638416192697294,0.7) node {$p_6$};
                \end{scriptsize}
            \end{tikzpicture}
        \end{example}

        \begin{proposition}[Recursive Structure]\label{recstruc}%recursive structutre
            Each boundary cycle is isomorphic to a product of moduli spaces of lower dimension.
            In particular, a boundary divisor of the type $D(A|B)$ with intersection node $x$ is isomorphic to the product $\overline{M}_{0,A\cup\{x\}}\times\overline{M}_{0,B\cup\{x\}}$.
            %justify dimension count, bdry cycle smooth and irreducible
        \end{proposition}
        The isomorphism in question simply separates the $A$- and $B$-twigs. This proposition justifies counting the dimension of the strata $\Sigma_\delta$ twig by twig, and also guarantees that boundary cycles are smooth and irreducible. 

        \begin{remark}
            Let $D(A|B)$ be a boundary divisor of $\overline{M}_{0,n}$ and consider the forgetful map $\epsilon:\overline{M}_{0,n+1}\rightarrow\overline{M}_{0,n}$ forgetting the last mark, say $q$. Then the pullback divisor is given by:
            $$\epsilon^*D(A|B)=D(A\cup\{q\}|B)+D(A|B\cup\{q\}).$$
        \end{remark}

        \begin{remark}[Prelude to Fundamental Boundary Relation]\label{bdryprel}%1.5.13 fund rel alpha
            For $n\ge4$ consider the forgetful map $\epsilon:\overline{M}_{0,n}\rightarrow\overline{M}_{0,4}$ (obtained by composition of $n-4$ maps forgetting one mark each) and write the set of marks of $\overline{M}_{0,4}$ as $\{i,j,k,l\}$. The three boundary divisors of $\overline{M}_{0,4}$ are $(ij|kl)$, $(ik|jl)$ and $(il|jk)$. Note that any two points in $\overline{M}_{0,4}\cong\pone$ are linearly equivalent. Pulling back these three equivalent divisors via $\epsilon$ therefore leads to the following equivalence:
            \begin{equation*}
                \sum_{\substack{i,j\in A\\k,l\in B}}D(A|B)\equiv\sum_{\substack{i,k\in A\\j,l\in B}}D(A|B)\equiv\sum_{\substack{i,l\in A\\j,k\in B}}D(A|B).
            \end{equation*}
        \end{remark}

        This concludes our initial discussion of curves. We now turn our attention to maps. This additional structure will allow us to exert control over the marked points using so called evaluation maps, which will be crucial in view of our enumerative goals.
        % here curves are present as images. looking at maps will allow us to involve marks which we need for enumerative goal.

        \begin{definition}%degree of a map
            We define the degree of a map $\mu:\pone\rightarrow\pr$ as the product of the transcendence degree of the field extension corresponding to $\mu$ and the degree of the image curve. In particular, if $\mu$ is birational, then its degree is the degree of its image curve.
        \end{definition}
        A degree $d$ map $\mu:\pone\rightarrow\pr$ is given up to a constant factor by $r+1$ binary forms of degree $d$ that do not vanish simultaneously at any point. 
        We denote the space formed by these selections of binary forms by $W(r,d)$. 

        \begin{proposition}%Wrd fine mod space, dim rd+r+d (2.1.1)
            The space $W(r,d)$ is a fine moduli space for maps $\pone\rightarrow\pr$ of degree $d$. It is of dimension $(r+1)(d+1)-1$.
        \end{proposition}
        Note that $(r+1)(d+1)$ is the number of coefficients of the $r+1$ binary coefficients, and we subtract $1$ because we only consider the binary forms up to a constant factor.
        % mention universal family? 2.1.1
        The universal family is given by the following diagram:
        \begin{center}
            \begin{tikzcd}
                W(r,d)\times\pone \arrow[rr,"\mu"] \arrow[dd,swap,"\pi"] && \pr \\\\
                W(r,d)
            \end{tikzcd}
        \end{center}
        where the horizontal arrow maps a pair $(\mu,x)$ to the point $\mu(x)$.

        Denote by $W^\circ(r,d)\subset W(r,d)$ the locus of immersions (morphisms with everywhere injective tangent map). It can be shown that it is open. %maps whose tangent map (def?) is injective at any point 
        Denote further by $W^*(r,d)\subset W(r,d)$ the (also open) locus of maps that are birational onto their image. Then we have inclusions $W^\circ(r,d)\subset W^*(r,d)\subset W(r,d)$, and for $d>0$ the complement of $W^*(r,d)$ consists of multiple covers. For $r\ge2$, $W^*(r,d)$ is dense in $W(r,d)$.

        \begin{remark}[Issues with $W(r,d)$]\label{issueswrd}%thee problems (and their solutions)
             There are three major issues rendering $W(r,d)$ ineffective as a tool for the description of families of rational curves. 
             
             Firstly, it may happen that curves in the same family are parametrized by a different $\pone$ (see Example 2.1.10 in \cite{invitation}). %2.1.10
             This problem can be solved by slightly altering the families under consideration. We do this via Definition \ref{famdefnew}.

             Secondly, reparametrizations of the same rational curve in $\pr$ yield distinct objects in $W(r,d)$. 
             This issue can be addressed by giving a variety structure to the quotient set $W(r,d)/\text{Aut}(\pone)$. More precisely, one can show (by comparing the function fields of the image curve and of $\pone$) that given a non-constant map $\mu:\pone\rightarrow\pr$, there are only finitely many automorphisms $\phi:\pone\rightarrow\pone$ with $\mu=\mu\circ\phi$. If $\mu$ is birational onto its image, there is precisely one such automorphism. This means that the locus $W^*(r,d)$ is precisely the subset of automorphism-free maps in $W(r,d)$. The intuition from Remark \ref{intuition} holds true in this case, as we shall see in Proposition \ref{fineandcoarse}. In order to further reduce the number of automorphisms, we shall introduce marks on our source curves. This will lead to the definition of the spaces $\m$ in Proposition \ref{m0n}.

             Finally, %non complete see exp2.2.1, rmk after 2.3.3 (cauchy seq) 2.2 1-par fams, th 2.3.2
             the space $W(r,d)$ is not complete (or compact) in the sense that the limits of certain families within $W(r,d)$ are not contained in $W(r,d)$. Looking at concrete one- and two-parameter families (e.g. a pencil of conics) indicates that in order to mend this, one must allow blow-ups of the base curve, which leads to trees whose twigs have degrees summing up to $d$ as sources of our maps. There are however some constraints, particularly on the twigs of degree zero (see section 2.2 of \cite{invitation} for details). This leads to the definition of Kontsevich-stable maps, and subsequently (in Theorem \ref{kontsspace}) to the spaces we shall work in throughout the rest of this thesis. These are the so called Kontsevich spaces $\mbar$, which serve as compactifications of $\m$.
        \end{remark}

        \begin{definition}\label{famdefnew}%tackles 1st concern
            A map of smooth rational curves is a map $\mu:C\rightarrow\pr$ with $C\cong\pone$.
            Two maps $\mu:C\rightarrow\pr$ and $\mu':C'\rightarrow\pr$ are isomorphic if there exists an isomorphism $\phi:C\rightarrow C'$ such that the following diagram commutes:
            
            \begin{center}
                \begin{tikzcd}
                    C \arrow[rr,"\phi"] \arrow[dr,swap,"\mu"] && C' \arrow[dl,"\mu'"]\\
                    &\pr
                \end{tikzcd}
            \end{center}

            Similarly, in the case $C=C'$, we say that both maps are automorphic if the map $\phi:C\rightarrow C$ in the diagram is an automorphism.
            
            A family of maps of smooth rational curves is a diagram of the following form:

            \begin{center}
                \begin{tikzcd}
                    \mathfrak{X} \arrow[rr,"\mu"] \arrow[dd,swap,"\pi"] && \pr \\\\
                    B
                \end{tikzcd}
            \end{center}
            
            such that $\pi:\mathfrak{X}\rightarrow B$ is a flat family whose geometric fibers $\mathfrak{X}_b$ are isomorphic to $\pone$ for all $b\in B$. Two such families $\pi$ and $\pi'$ are isomorphic if there exists an isomorphism $\phi:\mathfrak{X}\rightarrow\mathfrak{X}'$ making the following diagram commute:

            \begin{center}
                \begin{tikzcd}
                    &\pr\\
                    \mathfrak{X} \arrow[rr,"\phi"] \arrow[dr,swap,"\pi"] \arrow[ur,"\mu"] && \mathfrak{X}' \arrow[dl,"\pi'"] \arrow[ul,swap,"\mu'"]\\
                    &B
                \end{tikzcd}
            \end{center}

        Automorphisms are defined in the obvious way.
        \end{definition}
        
        Note that the maps $\mu_b:\mathfrak{X}_b\rightarrow\pr$ are maps from smooth rational curves. It can be shown that they all have the same degree. 

        \begin{proposition}\label{fineandcoarse}
            There exists a fine moduli space $M_{0,0}^*(\pr,d)\cong W^*(r,d)/\text{Aut}(\pone)$ for the classification of degree $d$ maps that are birational onto their image up to isomorphisms of $\pone$.

            Including the multiple covers from the complement of $W^*(r,d)$ yields a coarse moduli space $M_{0,0}(\pr,d)\cong W(r,d)/\text{Aut}(\pone)$.
        \end{proposition}
        The first zero in the subscripts again stands for the genus of the curves we are considering. The meaning of the second zero in the subscripts will become clear momentarily.
        The dimension of $M_{0,0}(\pr,d)$ is given by:
        \begin{equation}\label{dimm00}
            \dim(M_{0,0}(\pr,d))=\dim(W(r,d))-\dim(\text{Aut}(\pone))=rd+r+d-3.
        \end{equation} %argue with generic fiber of Wrd->M00prd

       Note that, since a linear map is always birational onto its image, we have that $M_{0,0}^*(\pr,1)\cong M_{0,0}(\pr,1)$. However, this is no longer true for $d>1$. Another way to eliminate isomorphisms is to add marks to the source curve and to force the automorphisms to respect this additional structure. Since our source curves are isomorphic to $\pone$, adding three or more marks eliminates all isomorphisms (but the trivial one). 

        \begin{proposition}\label{m0n}%2.1.16 m0n(...)
            Let $n\ge3$. Then there exists a fine moduli space 
            $$\m:=M_{0,n}\times W(r,d)$$
            for isomorphism classes of $n$-pointed maps $\pone\rightarrow\pr$ of degree $d$. It is a smooth variety.
        \end{proposition}
        The second subscript corresponds to the number of marks. The universal family for this new space can be constructed by combining those of the two factors.
        \begin{remark}
            Similarly, there exist fine moduli spaces $${M}^{\circ}_{0,n}(\pr,d) = {M}_{0,n}\times W^\circ(r,d),$$ 
            parametrizing immersions with smooth source curve, and 
            $${M}^{*}_{0,n}(\pr,d) = M_{0,n}\times W^*(r,d),$$ parametrizing automorphism-free maps with smooth source curve.
        \end{remark}

        \begin{definition}%n-pointed map
            An $n$-pointed map $(C;p_1,\dots,p_n;\mu)$ is a morphism $\mu:C\rightarrow\pr$, where $C$ is a tree of projective lines with $n$ distinct marks $p_i$ that are smooth points of $C$. 
            Two curves $\mu:C\rightarrow\pr$ and $\mu':C'\rightarrow\pr$ are isomorphic if there exists an isomorphism $\phi:C\rightarrow C'$ satisfying $\phi(p_i)=p_i'$ for all $1\le i\le n$ and making the first diagram from Definition \ref{famdefnew} commute. We can express this by the commutativity of the following diagram, writing $\bullet=\text{Spec}(\C)$ and defining sections $\sigma_i$ such that $\sigma_i(\bullet)=p_i$:
            \begin{center}
                \begin{tikzcd}
                   &\pr\\
                    C \arrow[rr,"\phi"] \arrow[dr,swap,shift right=.75ex,"\pi\textcolor{white}{'}"] \arrow[ur,"\mu\textcolor{white}{'}"] && C' \arrow[dl,shift left=.75ex,"\pi'"] \arrow[ul,swap,"\mu'"]\\
                    &\bullet \arrow[ul,swap,shift right=.75ex,"\sigma_i\textcolor{white}{'}",shorten >=1ex, shorten <=1ex] \arrow[ur,shift left=.75ex,"\sigma_i'",shorten >=1ex, shorten <=1ex]
                \end{tikzcd}
            \end{center}

            A family of $n$-pointed maps is a diagram of the following form:
            \begin{center}
                \begin{tikzcd}[column sep=huge,row sep=huge]
                   \mathfrak{X} \arrow[d,swap,shift right=.75ex,"\pi"] \arrow{r}{\mu} & \pr \\
                     B \arrow[u,swap,shift right=.75ex,"\sigma_i",shorten >=1ex, shorten <=1ex] 
                \end{tikzcd}
            \end{center}
            for a flat family $\pi:\mathfrak{X}\rightarrow B$ of trees of smooth rational curves with $n$ distinct sections $\sigma_i$ avoiding the singularities of the fibers of $\pi$. Isomorphisms are defined via the commutativity of the following diagram:
            % using a diagram like the one at the end of Definition \ref{famdefnew}, but this time including the sections $\sigma_i$ and $\sigma'_i$ for all $1\le i\le n$.
            \begin{center}
                \begin{tikzcd}
                   &\pr\\
                    \mathfrak{X} \arrow[rr,"\phi"] \arrow[dr,swap,shift right=.75ex,"\pi\textcolor{white}{'}"] \arrow[ur,"\mu\textcolor{white}{'}"] && \mathfrak{X}' \arrow[dl,shift left=.75ex,"\pi'"] \arrow[ul,swap,"\mu'"]\\
                    &B \arrow[ul,swap,shift right=.75ex,"\sigma_i\hphantom{'}\textcolor{white}{'}",shorten >=1ex, shorten <=1ex] \arrow[ur,shift left=.75ex,"\sigma_i'",shorten >=1ex, shorten <=1ex]
                    
                \end{tikzcd}
            \end{center}
        \end{definition}
        For the sake of brevity, we may simply write $\mu$ instead of $(C;p_1,\dots,p_n;\mu)$.

        \begin{definition}%konts-stability
            We say that an $n$-pointed map $\mu: C\rightarrow\pr$ is Kontsevich-stable if its degree zero twigs (i.e. those that are mapped to a point) are stable in the sense of $n$-pointed curves. 
        \end{definition}
        This means that each degree zero twig must contain three points that are either marks or nodes. The source curve $C$ is not necessarily stable as a pointed curve: a non-constant map $\pone\rightarrow\pr$ has a source curve made of a single positive degree twig and is therefore Kontsevich-stable. However, the source curve contains no marks or nodes and is therefore unstable as a pointed curve. From now on, we will often speak simply of 'stable' maps. The intended kind of stability can be deduced from the context.
        %combines trees and maps
        \begin{remark}
            An $n$-pointed map is Kontsevich-stable if and only if it has a finite number of automorphisms. This can be shown by examining what happens to twigs that are unstable as pointed curves and applying what we know about their automorphisms from Remark \ref{issueswrd}.
        \end{remark}
        
        In keeping with the intuition from Remark \ref{intuition}, there exists a coarse moduli space for these maps.
        \begin{theorem}\label{kontsspace}%coarse and fine moduli spaces (konts space)
            There exists a coarse moduli space $\mbar$ for the parametrization of isomorphism classes of Kontsevich-stable $n$-pointed maps of degree $d$ to $\pr$.

            This space is a projective normal irreducible variety. Furthermore, it is locally isomorphic to a quotient of a smooth variety by the action of a finite group. It contains a smooth open dense subvariety $\overline{M}_{0,n}^*(\pr,d)$ which is a fine moduli space for isomorphism-free maps.
        \end{theorem}
        For a more general version of this statement see Theorem 1 in \cite{FPnotes}. %thm1
        %TODO: proj implies separated and complete (recall: this is one of the advantages of proj spaces, separatedness as a kind of compactness)
        The dimension of this space is given by the following equation:
       \begin{equation}\label{dimmbar}%p2
                dim(\mbar)=rd+r+d+n-3. %(2.1.14)
        \end{equation}
        This uses equation (\ref{dimm00}) and the fact that each mark adds one dimension.

        \begin{definition}%evaluation maps
            Given a mark $p_i$, there is a natural map:
            \begin{align*}
                \nu_i:\mbar&\rightarrow\pr\\
                (C;p_1,\dots,p_n;\mu)&\mapsto\mu(p_i).
            \end{align*}
            We call these the evaluation maps. 
            %Here $(C;p_1,\dots,p_n;\mu)$ actually stands for an isomorphism class of $n$-pointed maps. 
            The evaluation maps are constant on isomorphism classes since by definition two isomorphic $n$-pointed curves satisfy $\mu(p_i)=\mu'(p_i')$ for all $i$.
        \end{definition}
        The evaluation maps are in fact flat morphisms. They will serve as a crucial tool when relating the geometry of $\pr$ with that of $\mbar$. As an example, note that the preimage of a hyperplane in $\pr$ under $\nu_i$ is the divisor of maps in $\mbar$ whose $i$-th mark gets mapped to that hyperplane.
        
        \begin{propdef}\label{forg}%forgetful map
            Given two sets of marks $B\subset A$, there is a natural forgetful map:
            \begin{equation*}
                \epsilon:\overline{M}_{0,A}(\pr,d)\rightarrow\overline{M}_{0,B}(\pr,d),
            \end{equation*}
            which forgets the marks in $A\setminus B$. It factors through forgetful maps forgetting one mark each. Special care must be taken to contract any degree $0$ twigs that would lead to Kontsevich-instability.

            For $n\ge3$, there is also a forgetful map $\mbar\rightarrow\overline M_{0,n}$ forgetting the map to $\pr$ and stabilizing the base curve. This map is a flat morphism.

            In particular, for $n\ge4$, there is a flat forgetful morphism $\mbar\rightarrow\overline M_{0,4}$.
        \end{propdef}

        %bdry, bdry cycle, bdry divisor
        We now look at the boundary $\mbar\setminus{M}_{0,n}(\pr,d)$ of $\mbar$. It parametrizes maps whose source curve is non-reduced. In addition to the distribution of marks that we encountered when discussing the boundary of $\overline{M}_{0,n}$, we will now have to distribute partial degrees.
        
        \begin{propdef} %D(A,B,dA,dB)
           A $d$-weighted partition of the set of marks $S=\{p_1,\dots,p_n\}$ consists of a partition $S=A\sqcup B$ together with a partition $d=d_A+d_B$, where $d_A$ and $d_B$ are non-negative integers.

           Given a $d$-weighted partition of $S$ such that $\abs{A}\ge2$ whenever $d_A=0$ and $\abs{B}\ge2$ whenever $d_B=0$, there is an associated irreducible divisor $D(A,B;d_A,d_B)$ which we call a boundary divisor. A general point of the boundary divisor represents a map $(C;p_1,\dots,p_n;\mu)$ whose domain is given by a tree with two twigs, i.e. $C=C_A\cup C_B$, such that the twigs $C_A$ and $C_B$ respectively carry the marks of $A$ and $B$, and such that the restrictions $\mu|_{C_A}$ and $\mu|_{C_B}$ are maps of degree $d_A$ and $d_B$, respectively. We illustrate this below.

           \definecolor{rvwvcq}{rgb}{0.08235294117647059,0.396078431372549,0.7529411764705882}
            \definecolor{wrwrwr}{rgb}{0.3803921568627451,0.3803921568627451,0.3803921568627451}
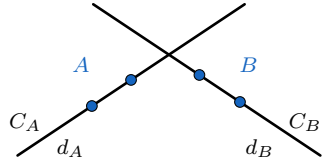
            \begin{tikzpicture}[line cap=round,line join=round,>=triangle 45,x=1cm,y=1cm]
                \clip(-5,0.5) rectangle (4.5,4); %size of rectangle containing image (given by 2pts)
                \draw [line width=1pt] (1,3)-- (4,1);
                \draw [line width=1pt] (3,3)-- (0,1);
                \begin{scriptsize}
                    \draw[color=black] (3.8,1.44) node {$C_{B}$};
                    \draw[color=black] (3.2,1.1) node {$d_{B}$};
                    
                    \draw [fill=rvwvcq] (2.4025391414130275,2.0649739057246483) circle (2pt);
                    % \draw[color=rvwvcq] (2.57,2.2890429805862818) node {$p_{1}$};
                    
                    \draw [fill=rvwvcq] (2.938269319967983,1.7078204533546781) circle (2pt);
                    \draw[color=rvwvcq] (3.05,2.2) node {$B$};
                    
                    \draw[color=black] (0.1,1.44) node {$C_{A}$};
                    \draw[color=black] (0.7,1.1) node {$d_{A}$};
                    
                    \draw [fill=rvwvcq] (1.4993503144908957,1.9995668763272638) circle (2pt);
                    %\draw[color=rvwvcq] (1.75,1.8) node {$m_{1}$};
                    
                    \draw [fill=rvwvcq] (0.9819670598590557,1.654644706572704) circle (2pt);
                    \draw[color=rvwvcq] (0.84,2.2) node {$A$};
                    
                    %\draw [fill=wrwrwr] (2,2.3333333333333335) circle (2pt);
                    %\draw[color=wrwrwr] (2,2.7) node {$x$};
                \end{scriptsize}
            \end{tikzpicture}
        \end{propdef}

        %2.7.2: what is finite quotient?
        \begin{proposition}\label{normalcrossingsmbar}%finite quotient = up to finitely many bdry divisors?
            Up to a finite quotient, the union of the boundary divisors in $\mbar$ is a divisor with normal crossings.
        \end{proposition}

        \begin{remark}[Recursive Structure]\label{recstruc}%2.7.3 chow ring  %gluing isomorphism.
           Consider the following fiber product:
           $$\overline{M}_{0,A\cup\{x\}}(\pr,d_A)\times_{\pr}\overline{M}_{0,B\cup\{x\}}(\pr,d_B),$$ 
           via the evaluation maps at the gluing mark $x$, i.e.: %\textcolor{red}{need twigs on lhs, so is the overline superfluous? mark always sits on a twig}
           \begin{align*}
               &\nu_{x_A}:\overline{M}_{0,A\cup\{x\}}(\pr,d_A)\rightarrow\pr,\\
               &\nu_{x_B}:\overline{M}_{0,B\cup\{x\}}(\pr,d_B)\rightarrow\pr.
           \end{align*} 
           Then there is a natural gluing map:
           $$\overline{M}_{0,A\cup\{x\}}(\pr,d_A)\times_{\pr}\overline{M}_{0,B\cup\{x\}}(\pr,d_B)\rightarrow D(A,B;d_A,d_B),$$
           gluing together a pair of twigs at the shared mark $x$. If $A\neq\emptyset$ and $B\neq\emptyset$, this map is an isomorphism; see Lemma 12 in \cite{FPnotes} for further details. 

           Denote by $\Delta$ the diagonal in $\pr\times\pr$. Then the fiber product can be seen as a subset of the usual product via the evaluation maps: $$(\nu_{x_A}\times\nu_{x_B})^{-1}(\Delta)\subset\overline{M}_{0,A\cup\{x\}}(\pr,d_A)\times\overline{M}_{0,B\cup\{x\}}(\pr,d_B).$$
           This allows us to examine intersections with the boundary divisor $D(A,B;d_A,d_B)$ in spaces of strictly smaller dimension, i.e. $\overline{M}_{0,A\cup\{x\}}(\pr,d_A)$ and $\overline{M}_{0,B\cup\{x\}}(\pr,d_B)$.
           
        \end{remark}
        %used in splitting lemma

        A special type of boundary divisor is obtained by pulling back the divisor $(ij|kl)$ from $\overline M_{0,4}$ via the flat forgetful morphism from Proposition-Definition \ref{forg}. We denote these divisors in $\mbar$ by $D(ij|kl)$. We can write them as formal sums:
        $$D(ij|kl)=\sum_{\substack{i,j\in A\\k,l\in B}} D(A,B;d_A,d_B),$$
        where we sum over all $d$-weighted partitions $A\sqcup B$ of $S=\{p_1,\dots,p_n\}$.
        Arguing as in Remark \ref{bdryprel}, we obtain the following theorem.

        \begin{theorem}[Fundamental Boundary Relation]\label{bdry}%TODO: see also FP notes intro!
        The following equivalence of divisors holds in $\mbar$:
            \begin{equation}\label{fundrel}
                \sum_{\substack{A\sqcup B=S\\i,j\in A\\k,l\in B\\d_A+d_B=d}}D(A,B;d_A,d_B)\;\equiv\sum_{\substack{A\sqcup B=S\\i,k\in A\\j,l\in B\\d_A+d_B=d}}D(A,B;d_A,d_B)\;\equiv\sum_{\substack{A\sqcup B=S\\i,l\in A\\j,k\in B\\d_A+d_B=d}}D(A,B;d_A,d_B).
            \end{equation}
        \end{theorem} %proof: M04bar isomorphic to P1, see if you can find junliangs explanation 
        
        We list a few more useful properties.

        \begin{proposition}[Compatibility]\label{compeval}%2.8.1, see also 1.5.9, and 2.7.3 for recursive structure
            %from easy props
            Given a partition $S=A\sqcup B$ such that $p_i\in A$, the following diagram commutes:
            \begin{center}
                \begin{tikzcd}
                    \overline{M}_{0,A\cup\{x\}}(\pr,d_A)\times_{\pr}\overline{M}_{0,B\cup\{x\}}(\pr,d_B) \arrow[r] \arrow[d]&D(A,B;d_A,d_B)\subset\mbar \arrow[d,"\nu_i"] \\
                    \overline{M}_{0,A\cup\{x\}}(\pr,d_A) \arrow[r, swap, "\nu_i"] & \pr
                \end{tikzcd}
            \end{center}
        \end{proposition}

        This proposition can be used to show the next.

        \begin{proposition}%proper intersection
            Let $\Gamma\subset\pr$ be a subvariety. Its inverse image $\nu_i^{-1}(\Gamma)\subset\mbar$ has proper intersection with each boundary divisor $D(A,B;d_A,d_B)$.
        \end{proposition}
        This means that if $\text{codim}_{\pr}(\Gamma)=k$, then: 
        $$\text{codim}_{\mbar}(\nu_i^{-1}(\Gamma)\cap D(A,B;d_A,d_B))=k+1.$$
        
        In case $d=0$, the entire source curve gets sent to a single point. By Kontsevich stability, we therefore have $n\ge3$. As such we can consider the forgetful morphism $\mbar\rightarrow\overline M_{0,n}$ mentioned previously and, of course, the evaluation maps $\nu_i$. As intuition suggests, one can then show that the product of these two maps yields an isomorphism, i.e. that giving a degree zero pointed map is equivalent to giving a pointed curve together with a point in $\pr$.
        
        \begin{proposition}[Degree 0]\label{deg0iso}%identification 2.8.5
            There exists an isomorphism between the following spaces:
            $$\overline{M}_{0,n}(\pr,0)\cong\overline{M}_{0,n}\times\pr.$$
        \end{proposition}

        This also shows that one can recover the spaces $\overline{M}_{0,n}$ from $\mbar$: just set $r=0$ in the proposition.

    \section{Moduli Spaces of Curves and Maps for  $\ponex$}\label{th1p1x}
         %mention both interpretations of ponex (via embedding) and of bidegree (d,e) maps
         %define bidegree, define W((e,d)),Wdeg,W* and mention 3 problems
         %explain horizontal and vertical rules MAKE SURE WHETHER BOOK HAS MISTAKE; SEE PROBLEMS NOTE
         %\textcolor{red}{SETTLE ON DEFINITION OF HORIZONTAL} %if we diverge from book, must exchange all hor/vert until now
         
         %explain terminology: (d,e)-curve means rational curve of bidegree (d,e)
         %rational map is of bidegree (d,e) if its img is a (d,e)-curve
         %3 equivalent ways of thinking about (d,e)-curves: V(f)xV(g), V(F), sum of fundamental classes. EXPLAIN ALL THREE PLUS RELATIONS
         
         %introcu^duce shorthand notation: mu ind^stead of (C;p1,...,pn;mu)
         %see last two pages of n3 note

         %LIST IMPORTANT CHANGES AND RESULTS P.R.A. THE ABOVE: DEF OF BIDEGREE, DIMENSION OF MBAR, ... CF NOTES

         We now consider maps into $\ponex$ or, equivalently, into the quadric surface in $\mathbb{P}^3$ via the Segre embedding. While the general spirit of the previous section remains the same, the objects in question are slightly different. In particular, curves now have bidegrees instead of degrees. On the other hand, the dimension of our space is now fixed at $\dim(\ponex)=2$. Despite these changes, all of the above statements have their analogous counterparts in the new setting. We mention only the most notable adaptations.
         
        %\textcolor{red}{have altered definition of hor/vert given in book, see problems note}

        We begin by defining the most basic types of curves in $\ponex$.
        \begin{definition}%horizontal and vertical rules (curves)
            We call horizontal rule a line of the form $\pone\times\{[y_0:y_1]\}\subset\ponex$, where $[y_0:y_1]$ is a point in $\pone$.

            A vertical rule is a line $\{[x_0,x_1]\}\times\pone\subset\ponex$ for $[x_0:x_1]\in\pone$.
        \end{definition}
        
        Let us now look at what other types of curves appear in $\ponex$.
        \begin{definition}%bidegree (curves)
            We say that a curve $C$ in $\ponex$ is of bidegree $(d,e)$ if it is a product of vanishing loci, i.e. $C=V(f)\times V(g)$, where $f(x_0,x_1)$ and $g(y_0,y_1)$ are homogeneous polyomials of degrees $d$ and $e$, respectively. 
            Equivalently, it can be described as the locus $\{([x_0:x_1],[y_0:y_1])\}$ of points in $\ponex$ that satisfy $F(x_0,x_1,y_0,y_1)=0$ where $F$ is a homogeneous polynomial of bidegree $(d,e)$, i.e. homogeneous of degree $d$ in $(x_0,x_1)$ and homogeneous of degree $e$ ind $(y_0,y_1)$.
            We may briefly refer to such a curve as a $(d,e)$-curve. 
            
            In particular, a horizontal rule is given by $V(0)\times V(g)$ with $\text{deg}(g)=1$ and is therefore a $(0,1)$-curve. Accordingly, a vertical rule is a $(1,0)$-curve.
            
        \end{definition}
        %linked by cohom to rules
        From a cohomological point of view, the class of a $(d,e)$-curve is given by $d\cdot V + e\cdot H$, where $V$ and $H$ respectively denote the classes of a vertical and a horizontal rule.

        Inspired by B\'ezout's Theorem and in order to make the upcoming connections between the $\ptwo$ and $\ponex$ cases more apparent, we develop some new notation. We introduce the product of two bidegrees $(d_A,e_A)$ and $(d_B,e_B)$ as the multiplicity of the intersection of two curves with these respective bidegrees: $$(d_A,e_A)\circ(d_B,e_B):=d_A e_B \cdot e_A d_B.$$
        In particular, we have $(d,e)\circ(1,0)=e$, $(d,e)\circ(0,1)=d$ and $(d,e)\circ(1,1)=d+e$.

        %now analogous spaces etc, name just mbarpone. konts stable analogous. eval, forg are the same. 
        %bidegree of maps is defined componentwise like pr case
        Using these definitions, one can replicate all the theory from the previous section. The maps now go to $\ponex$ instead of $\pr$, and their bidegree is defined componentwise like the degree of a map to $\pr$. With the same definition of Kontsevich-stability, one finds a coarse moduli space $\monebar$ parametrizing isomorphism classes of of Kontsevich-stable $n$-pointed maps of bidegree $(d,e)$ into $\ponex$. The dimension of this space is given by:
         \begin{equation}\label{dimmbarone}%see quantum notes
             dim\bigr(\monebar\bigl)=n+2d+2e-1,
        \end{equation}
        and evaluation and forgetful maps work exactly the same.  

        The boundary structure is again crucial. It now takes into account the bidegrees.

        \begin{definition}%d,e-valued partition
            A $(d,e)$-weighted partition of $S=\{p_1,\dots,p_n\}$ is a partition $S=A\sqcup B$ together with partitions $d=d_A+d_B$ and $e=e_A+e_B$ of both degrees into non-negative integers.
        \end{definition}

        Accordingly, we get boundary divisors $D\bigr(A,B;(d_A,e_A),(d_B,e_B)\bigl)$, which again satisfy a fundamental relation:
        \begin{align}\label{fundrelonex}
            \begin{split}
                \sum_{\substack{A\sqcup B=S\\i,j\in A\\k,l\in B\\d_A+d_B=d\\e_A+e_B=e}}D\bigr(A,B;(d_A,e_A),(d_B,e_B)\bigl)\;&\equiv\sum_{\substack{A\sqcup B=S\\i,k\in A\\j,l\in B\\d_A+d_B=d\\e_A+e_B=e}}D\bigr(A,B;(d_A,e_A),(d_B,e_B)\bigl)\\
                &\equiv\sum_{\substack{A\sqcup B=S\\i,l\in A\\j,k\in B\\d_A+d_B=d\\e_A+e_B=e}}D\bigr(A,B;(d_A,e_A),(d_B,e_B)\bigl).
            \end{split}
        \end{align}
        % $$\sum_{\substack{m_1,m_2\in A \\ p_1,p_2\in B\\d_A+d_B=d  \\e_A+e_B=e}}D(A,B;(d_A,e_A),(d_B,e_B)text{ and} \sum_{\substack{m_1,p_1\in A \\ m_2,p_2\in B \\d_A+d_B=d  \\e_A+e_B=e}}D\bigr(A,B;(d_A,e_A),(d_B,e_B)\bigl),$$
        We are now equipped to tackle the enumerative problem posed in the Introduction.

\chapter{Kontsevich's Formula: Classical Approach}\label{classkonts}

    In this chapter, we use the theory established in the previous ones to solve the moduli problem of classifying degree $d$ rational curves passing through finitely many points explained in the Introduction. We will at first compute some of the desired numbers directly in order to get a feeling for how to apply our theory. At the end, we prove a recursive relation giving us all the desired numbers. We do this first for the space $\ptwo$, then for $\ponex$.
    \section{Kontsevich's Formula for $\ptwo$}\label{kfp2}
        %TODO motivate Kontsevich formula historically 
        
        %discussion on why 3d-1 is the right amount of points
        We must first determine the necessary number of point conditions for fixing a finite number of irreducible curves of degree d in $\ptwo$. 
        
        Let us examine the space of irreducible plane curves of degree d. For $d=1$, such a curve is determined by the three coefficients $a$, $b$ and $c$ of the homogeneous polynomial $ax+by+cz$ up to scaling. In this case, the desired space is therefore simply $\ptwo$. For general $d$, we must count the number of coefficients of an homogeneous polynomial of degree d in $x$, $y$ and $z$. In other words, we must count the number of monomials $x^iy^jz^k$ with $i+j+k=d$. There are $d+1$ choices for $i$. Once $i$ is set, we can choose $j$ from $\{0,\dots,d-i\}$ (so $d-i+1$ possibilities), and the final coefficient must be $k=d-i-j$. So summing up the number of choices of $j$ for each $i$ yields the desired number: 
        $$\sum_{i=0}^d d-i+1=\sum_{k=1}^{d+1}k=\frac{(d+1)(d+2)}{2}=\frac{d(d+3)}{2}+1,$$
        where we substituted $k:=d+i-1$. The space of irreducible plane curves of degree d is thus $\mathbb{P}^{d(d+3)/2}$.

        Now the Genus Formula for nodal plane curves (see Theorem \ref{genus}) gives the genus $g$ of a plane curve as a function of $\delta$, the number of its nodes: 
        \begin{equation}\label{genform}
            g=\frac{(d-1)(d-2)}{2}-\delta.
        \end{equation}
        The curves we are interested in are rational curves, i.e. curves with $g=0$. So by \ref{genform}, we need $\delta={(d-1)(d-2)}/{2}$ nodes. Selecting a node in $\mathbb{P}^{d(d+3)/2}$ constitutes a codimension 1 condition (i.e. it allows the reformulation of one of the coefficients in terms of the others), so the space of rational curves of degree $d$ has dimension $d(d+3)/2-{(d-1)(d-2)}/{2}=3d-1.$

        In order to obtain a finite number of rational curves by adding point constraints, we must therefore choose exactly $3d-1$ general points in the plane (since the selection of a point is again a codimension 1 condition).
        
        We shall denote the number of rational plane curves of degree $d$ passing through $3d-1$ points in general position by $N_d$.

        \begin{remark}
            We say that two points $Q_1$, $Q_2$ in $\ptwo$ are in general position if $Q_1\neq Q_2$. Likewise, two lines $L_1$, $L_2$ in $\ptwo$ are said to be in general position if $L_1\neq L_2$. We refer to these conditions as genericity.
        \end{remark}

        The following fact was mentioned in the Introduction.
        \begin{fact}\label{n1} %reference Euclid, Whitehead's axioms (?)
            There is exactly $N_1=1$ line passing through 2 points in general position in $\ptwo$.
        \end{fact}
        %explain how this expands our definition of genericity
        Two points in general position thus correspond to a unique line. The points on this line are in this sense in a special positon in relation to the first two points. We extend our notion of genericity accordingly.

        \begin{remark}
            If more than $3d-1$ points lie on a degree $d$ rational curve, they are not in general position.
            A set of points and (one or two) lines is in general position if the points and the lines are in general position respectively amongst themselves, and none of the points lie on any of the lines.
        \end{remark}

        \begin{proposition}\label{n2} 
            There is exactly $N_2=1$ rational conic passing through 5 points in general position in $\ptwo$.
        \end{proposition}
        
        \begin{proof} %perhaps give this as is for motivation, then justify why it works in lemmas for n3.
            The central idea of the proof is to work in the space $\overline{M}_{0,6}(\ptwo,2)$ of six-pointed stable maps of degree 2 into $\ptwo$ and to exploit our knowledge of its boundary. We denote the six marks by $m_1$, $m_2$, $p_1$, \dots, $p_4$. Let $L_1$ and $L_2$ be two lines and $Q_1, \dots, Q_4$ four points in general position in $\ptwo$. We will refer to these as our data. Looking at stable maps $(C;m_1, m_2, p_1, \dots, p_4;\mu)$ sending $m_1$ to $L_1$, $m_2$ to $L_2$ and $p_i$ to $Q_i$ for $1\leq i\leq 4$ will grant us the necessary flexibility to recover the desired number $N_2$ in terms of the known number $N_1$. %from Fact \ref{n1}.

            More precisely, let $Y\subset\overline{M}_{0,6}(\ptwo,2)$ be the subvariety of the maps mentioned above: 
            $$Y:=\nu_{m_1}^{-1}(L_1) \cap \nu_{m_2}^{-1}(L_2) \cap \nu_{p_1}^{-1}(Q_1) \cap \dots \cap \nu_{p_4}^{-1}(Q_4).$$
            We will show (Lemmas \ref{lemma1} and \ref{lemma2}) that the genericity of the lines $L_i$ and points $Q_i$ guarantees certain desirable properties of $Y$: firstly it is a one-dimensional subvariety in the 11-dimensional projective variety $\overline{M}_{0,6}(\ptwo,2)$, and secondly the intersection of $Y$ with each boundary divisor is transversal and takes place in $\overline{M}_{0,6}^*(\ptwo,2)$. This will be discussed in more detail in Remark \ref{onY}.
            %mention SMOOTHNESS?
            
            %now replicate structure of n3 proof.
            %  Via the forgetful map $\overline{M}_{0,6}(\ptwo,2) \rightarrow \overline{M}_{0,\{m_1,m_2,p_1,p_2\}}$ (forgetting $p_3$ and $p_4$) and the fundamental boundary relation (\ref{bdry}), we obtain the following equivalence:
            
            % \begin{equation}\label{equivn2}
            %       Y\cap D(m_1,m_2|p_1,p_2)\equiv Y\cap D(m_1,p_1|m_2,p_2),
            %  \end{equation}
            %  where $D(m_1,m_2|p_1,p_2)$ and $D(m_1,p_1|m_2,p_2)$ represent sums of the form $$\sum_{\substack{m_1,m_2\in A \\ p_1,p_2\in B \\d_A+d_B=2}}D(A,B;d_A,d_B) \; \text{ and} \sum_{\substack{m_1,p_1\in A \\ m_2,p_2\in B \\d_A+d_B=2}}D(A,B;d_A,d_B),$$ respectively.
            %  \textcolor{red}{but an equivalence of what? the intersection has codimension n and is not a divisor...}

            Since $Y$ and $D(m_1,m_2|p_1,p_2)$ have complementary dimension and intersect transversally, they intersect in $\deg(Y)\cdot\deg(D(m_1,m_2|p_1,p_2))$ points by a generalized version of B\'ezout's Theorem (see Corollary 2.4 in \cite{eh2}). %eisenbud harris 2
            This obviously continues to hold if we replace $D(m_1,m_2|p_1,p_2)$ by $D(m_1,p_1|m_2,p_2)$. Consider the forgetful map $\overline{M}_{0,6}(\ptwo,2) \rightarrow \overline{M}_{0,\{m_1,m_2,p_1,p_2\}}$ forgetting $p_3$ and $p_4$. By the fundamental boundary relation (\ref{bdry}), we know that $D(m_1,m_2|p_1,p_2)\equiv D(m_1,p_1|m_2,p_2)$. Recall that rationally equivalent divisors have the same degree, since principal divisors have degree $0$. Combining these facts, we conclude that the intersections $Y\cap D(m_1,m_2|p_1,p_2)$ and $Y\cap D(m_1,p_1|m_2,p_2)$ contain the same number of points. We summarize this by writing:
            \begin{equation}\label{equivn2}
                 \abs{Y\cap D(m_1,m_2|p_1,p_2)} = \abs{Y\cap D(m_1,p_1|m_2,p_2)}.
            \end{equation}
            Recall that $D(m_1,m_2|p_1,p_2)$ and $D(m_1,p_1|m_2,p_2)$ represent sums of irreducible components of the forms $$\sum_{\substack{m_1,m_2\in A \\ p_1,p_2\in B \\d_A+d_B=2}}D(A,B;d_A,d_B) \; \text{ and} \sum_{\substack{m_1,p_1\in A \\ m_2,p_2\in B \\d_A+d_B=2}}D(A,B;d_A,d_B),$$ respectively.
            %https://encyclopediaofmath.org/wiki/Intersection_theory
            
            Counting up and subsequently equating the contributions of the left and the right hand side will yield the result.
             
            LHS:\\
            Let us begin by illustrating $D_{AB}:=D(A,B;d_A,d_B)$ before distributing the spare marks $p_3$ and $p_4$.
            
            \definecolor{rvwvcq}{rgb}{0.08235294117647059,0.396078431372549,0.7529411764705882}
            \definecolor{wrwrwr}{rgb}{0.3803921568627451,0.3803921568627451,0.3803921568627451}
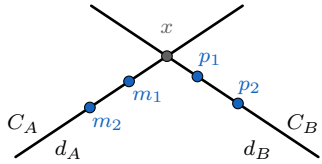
            \begin{tikzpicture}[line cap=round,line join=round,>=triangle 45,x=1cm,y=1cm]
                \clip(-5,0.5) rectangle (4.5,4); %size of rectangle containing image (given by 2pts)
                \draw [line width=1pt] (1,3)-- (4,1);
                \draw [line width=1pt] (3,3)-- (0,1);
                \begin{scriptsize}
                    \draw[color=black] (3.8,1.44) node {$C_{B}$};
                    \draw[color=black] (3.2,1.1) node {$d_{B}$};
                    
                    \draw [fill=rvwvcq] (2.4025391414130275,2.0649739057246483) circle (2pt);
                    \draw[color=rvwvcq] (2.57,2.2890429805862818) node {$p_{1}$};
                    
                    \draw [fill=rvwvcq] (2.938269319967983,1.7078204533546781) circle (2pt);
                    \draw[color=rvwvcq] (3.092541216255241,1.9233276303855245) node {$p_{2}$};
                    
                    \draw[color=black] (0.1,1.44) node {$C_{A}$};
                    \draw[color=black] (0.7,1.1) node {$d_{A}$};
                    
                    \draw [fill=rvwvcq] (1.4993503144908957,1.9995668763272638) circle (2pt);
                    \draw[color=rvwvcq] (1.75,1.8) node {$m_{1}$};
                    
                    \draw [fill=rvwvcq] (0.9819670598590557,1.654644706572704) circle (2pt);
                    \draw[color=rvwvcq] (1.2162624630513585,1.44) node {$m_{2}$};
                    
                    \draw [fill=wrwrwr] (2,2.3333333333333335) circle (2pt);
                    \draw[color=wrwrwr] (2,2.7) node {$x$};
                \end{scriptsize}
            \end{tikzpicture}

            Here $x$ denotes the intersection node of $C_A$ and $C_B$.
            Since we must have $d_A+d_B=2$, the only three possible pairs of degrees are $(d_A,d_B)\in\{(0,2),(1,1),(2,0)\}$. The two marks $p_3$ and $p_4$ remain to be distributed between $C_A$ and $C_B$. There are $\binom{2}{0}+\binom{2}{1}+\binom{2}{2}=4$ such distributions. Note that adding $k\in\{0,1,2\}$ marks to one twig is equivalent to adding the other $2-k$ marks to the other twig. 
            The number of irreducible components of the divisor $Y\cap D(m_1,m_2|p_1,p_2)$ is therefore equal to $3\cdot4=12$, i.e. there are 12 summands of the form $D_{AB}$ on the left hand side. 
            
            We now study the intersection of $Y$ with the divisors $D_{AB}$ by separately considering each of the three pairs of degrees $(d_A,d_B)$. In each case, we will determine the number of points in the intersection $Y\cap D(m_1,m_2|p_1,p_2)$ coming from that pair. In other words, we shall compute the number of stable maps subject to the constraints of being in $Y$ and having twigs $C_A$ and $C_B$ of respective degrees $d_A$ and $d_B$ with $d_A+d_B=2$. \\

            Ad $(2,0)$:
            No component with this pair of degrees can contribute to the intersection $Y\cap D(m_1,m_2|p_1,p_2)$. Indeed, $Q_1=\mu(p_1)=\mu(p_2)=Q_2$ would contradict the genericity of our data. \\
            
            Ad $(0,2)$:
            In this case we get $\{\mu(x)\}=\{\mu(m_1)\}=\{\mu(m_2)\}=L_1\cap L_2$. Unless both spare marks lie on $C_B$, $\mu(p_3)=Q_3$ or $\mu(p_4)=Q_4$ coincides with $L_1\cap L_2$, contradicting the genericity of the data (indeed, any $Q_i$ lying on $L_1$ or $L_2$ would suffice for a contradiction). 
            Thus only the divisor with $p_1$, \dots, $p_4$ on $C_B$, illustrated below, can yield a contribution. 

            
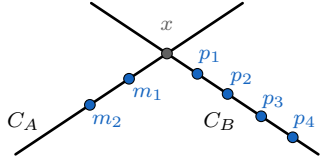
\begin{tikzpicture}[line cap=round,line join=round,>=triangle 45,x=1cm,y=1cm]
            \clip(-5,0.5) rectangle (4.5,4); %size of rectangle containing image (given by 2pts)
            \draw [line width=1pt] (1,3)-- (4,1);
            \draw [line width=1pt] (3,3)-- (0,1);
            \begin{scriptsize}
                \draw[color=black] (2.7,1.44) node {$C_{B}$};
                
                \draw[color=black] (0.1,1.44) node {$C_{A}$};
                
                \draw [fill=wrwrwr] (2,2.3333333333333335) circle (2pt);
                \draw[color=wrwrwr] (2,2.7) node {$x$};
                
                \draw [fill=rvwvcq] (0.9819670598590557,1.654644706572704) circle (2pt);
                \draw[color=rvwvcq] (1.2162624630513585,1.44) node {$m_{2}$};
                
                \draw [fill=rvwvcq] (1.4993503144908957,1.9995668763272638) circle (2pt);
                \draw[color=rvwvcq] (1.75,1.8) node {$m_{1}$};
                
                \draw [fill=rvwvcq] (2.4025391414130275,2.0649739057246483) circle (2pt);
                \draw[color=rvwvcq] (2.57,2.2890429805862818) node {$p_{1}$};
                
                \draw [fill=rvwvcq] (2.801758442831036,1.7988277047793095) circle (2pt);
                \draw[color=rvwvcq] (2.99,2) node {$p_{2}$};
                
                \draw [fill=rvwvcq] (3.245069807567101,1.503286794955266) circle (2pt);
                \draw[color=rvwvcq] (3.41,1.7) node {$p_{3}$};
                
                \draw [fill=rvwvcq] (3.6617003957218293,1.2255330695187805) circle (2pt);
                \draw[color=rvwvcq] (3.82,1.43) node {$p_{4}$};
            \end{scriptsize}
            \end{tikzpicture}
            
            The twig $C_B$ gets mapped to a conic passing through the five points $Q_1$, \dots, $Q_4$ as well as through the point $\mu(C_A)$. Thus the contribution from this divisor is equal to the number of conics through five points in general position, precisely the number $N_2$ we are looking for.
            
            This argument will be formalized in Lemmas \ref{lemma1} and \ref{lemma2}. \\

            Ad $(1,1)$: If one or both of the spare marks were to lie on $C_B$, then three or four of the $Q_i$ would lie on the line $\mu(C_B)$. Since this would contradict the genericity of our data, only the divisor with $m_1$, $m_2$, $p_3$ and $p_4$ on $C_A$ can contribute to the sum. 
            
            For the same reason, the image lines of $C_A$ and $C_B$ must be distinct. Furthermore, by Fact \ref{n1}, $\mu(C_A)$ is uniquely defined by $Q_3$ and $Q_4$, as is $\mu(C_B)$ by $Q_1$ and $Q_2$. This means that there is precisely $N_1\cdot N_1=1$ way of choosing the image lines of $C_A$ and $C_B$. %this is where we use N1
            Similarly, the images of the marks $m_1$ and $m_2$ must respectively lie in $\mu(C_A)\cap L_1$ and $\mu(C_A)\cap L_2$. Both of these intersections contain only one point, since $\mu(C_A) = L_i$ would imply $Q_3,Q_4\in L_i$ for $i\in\{1,2\}$ and thus contradict the genericity of the data.
            
            The final special point, namely the node $x$, must be mapped to the unique intersection point of $\mu(C_A)$ and $\mu(C_B)$, i.e. there is only one way of 'gluing' these two image lines together.
            
            Therefore, since the positions of all the special points are fixed, there is exactly one stable map fulfilling all the necessary properties, and the contribution from this pair is equal to 1.
            \\
            \\
            As discussed above, these are all possible pairs of degrees. The left hand side of equation \ref{equivn2} thus yields the number $N_2+1$.
            \\
            \\
            RHS:\\
            As for the left hand side, we begin with an illustration of the possible $D_{AB}$ before distribution of the spare marks $p_3$ and $p_4$.
            
            \definecolor{rvwvcq}{rgb}{0.08235294117647059,0.396078431372549,0.7529411764705882}
            \definecolor{wrwrwr}{rgb}{0.3803921568627451,0.3803921568627451,0.3803921568627451}
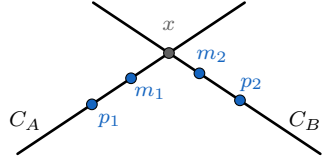
            \begin{tikzpicture}[line cap=round,line join=round,>=triangle 45,x=1cm,y=1cm]
                \clip(-5,0.5) rectangle (4.5,4); %size of rectangle containing image (given by 2pts)
                \draw [line width=1pt] (1,3)-- (4,1);
                \draw [line width=1pt] (3,3)-- (0,1);
                \begin{scriptsize}
                    \draw[color=black] (3.8,1.44) node {$C_{B}$};
                    
                    \draw [fill=rvwvcq] (2.4025391414130275,2.0649739057246483) circle (2pt);
                    \draw[color=rvwvcq] (2.57,2.2890429805862818) node {$m_2$};
                    
                    \draw [fill=rvwvcq] (2.938269319967983,1.7078204533546781) circle (2pt);
                    \draw[color=rvwvcq] (3.092541216255241,1.9233276303855245) node {$p_{2}$};
                    
                    \draw[color=black] (0.1,1.44) node {$C_{A}$};
                    
                    \draw [fill=rvwvcq] (1.4993503144908957,1.9995668763272638) circle (2pt);
                    \draw[color=rvwvcq] (1.75,1.8) node {$m_{1}$};
                    
                    \draw [fill=rvwvcq] (0.9819670598590557,1.654644706572704) circle (2pt);
                    \draw[color=rvwvcq] (1.2162624630513585,1.44) node {$p_1$};
                    
                    \draw [fill=wrwrwr] (2,2.3333333333333335) circle (2pt);
                    \draw[color=wrwrwr] (2,2.7) node {$x$};
                \end{scriptsize}
            \end{tikzpicture}

            The possible $(d_A,d_B)$, the spare marks to distribute and thus the number of irreducible components are the same as for the left hand side. We proceed to an analogous examination of the divisors $D_{AB}$ for each pair of degrees.

            Ad $(0,2)$, $(2,0)$: If $(d_A,d_B)=(0,2)$, then $Q_1=\mu(p_1)=\mu(m_1)\in L_1$. Similarly, if $(d_A,d_B)=(2,0)$, then $Q_2\in L_2$. In both cases we would get a contradiction to the genericity of our data. Consequently neither case contributes to the right hand side.\\

            Ad $(1,1)$: If either twig contained both spare marks (and thereby the other neither), the image line of that twig would contain three of the $Q_i$'s. Since this would contradict the genericity of our data, only components with $p_3$ on one twig and $p_4$ on the other can yield a contribution to the intersection $Y\cap D(m_1,p_1|m_2,p_2)$. Clearly, there are 2 such distributions.
            
            In both cases, each twig is mapped to a line containing two $Q_i$'s. By Fact \ref{n1}, these lines are uniquely determined. Similarly, $m_1$ must be mapped to the unique intersection point of $\mu(C_A)$ and $L_1$ and $m_2$ must be mapped to the unique intersection point of $\mu(C_B)$ and $L_2$. There is thus only $N_1\cdot N_1=1$ choice of the image lines in either case.

            In each case, the intersection node $x$ can only be mapped to the unique intersection of the image lines of $C_A$ and $C_B$. 
            \\
            \\
            The right hand side of the equation thus yields the number 2.
            Combining the results of both counts, equation \ref{equivn2} implies that $N_2+1=2$, i.e. $N_2=1$, as claimed.
            
        \end{proof}
        
        %we now formalize the arguments from n2
        %explain roadmap. explain why we require sum of codims = dim mbar. define Y, m_1,f,p_2. explain C_A, d_A,...
        %explain the definition of genericity in dependence of N_d

        Let us now establish a formal basis for the more subtle arguments in the above proof. They will then be used more explicitly in the proofs of the more general statements to come. %We do so in greater generality than immediately necessary and will thus be able to use the same results once we replace $\ptwo$ by $\ponex$.
        We shall work with the space $\overline{M}:=\overline{M}_{0,n}(X,d)$ of n-pointed stable maps of degree $d$ into $X:=\pr$ for $r\geq2$. %a variety $X$, where $X$ can be thought of as as a stand-in for $\ptwo$.  %not $\ponex$, bc that need bidegrees. 
        We denote the $n$ marks by $p_1$, \dots, $p_n$. 

        Let $X^n=X\times\dots\times X$ be the $n$-fold product of $X$ with itself and $\tau_i:X^n\rightarrow X$ the $i$-th projection for $1\leq i \leq n$. Consider $n$ irreducible %bc we need pts and lines
        subvarieties $\Gamma_i\subset X$ with respective codimensions $k_i$, and denote their product by 
        \begin{equation}\label{defgammabar}
            \underline\Gamma=\Gamma_1\times\dots\times\Gamma_n=\cap\tau_i^{-1}(\Gamma_i)\subset X^n.
        \end{equation}
        %\textcolor{red}{codim gamma sing gamma neq 0?}
        Similarly, let $\underline\nu:\overline{M}\rightarrow X^n$ be the map induced by the $n$ evaluation maps $\nu_i$.

        We are interested in the locus of maps $\mu$ that send the mark $p_i$ to the subvariety $\Gamma_i\subset X$ for each $1\leq i \leq n$, which in our new notation is given by the scheme-theoretic intersection $$\underline\nu^{-1}(\underline\Gamma)=\nu_1^{-1}(\Gamma_1)\cap\dots\cap\nu_n^{-1}(\Gamma_n).$$ 
        Note that by flatness of the evaluation maps $\nu_i$, we have $\text{codim}_{\overline{M}}(\nu_i^{-1}(\Gamma_i))=\text{codim}_X(\Gamma_i)=k_i$. However, this flatness property does not necessarily translate to the product map $\underline\nu$ (see example 2.5.4 in \cite{invitation}), so the above locus may not be of codimension $\sum k_i$.

        From now on, we will examine the case $\sum k_i=\text{dim} (\overline{M})$, since this is precisely the setting we need for the proofs of this section. In fact, if $r=2$, $n=3d-1$ and all the $\Gamma_i$ are points (i.e. $k_i=2$ for $1\leq i \leq n$), then 
        $$\sum k_i=2(3d-1) = \text{dim} (\overline{M}_{0,3d-1}(\ptwo,d))=\text{dim} (\overline{M}),$$
        where we used equation (\ref{dimmbar}).
        
        Intuitively, given general enough $\Gamma_i$ with $\sum k_i=\text{dim} (\overline{M})$, the intersection of their preimages should have dimension equal to $0$. As we shall see in the following lemma, this is ineed the case.

        \begin{lemma}\label{lemma1}
            The locus $\underline\nu^{-1}(\underline\Gamma)$ is smooth and enumerable: if $\sum k_i=\text{dim} (\overline{M})$, it consists of finitely many reduced %reduced: does this mean reduced source curve?
            %since \C is alg. closed, reduced means that the point does not lie in the intersection of two irred. components: https://math.stackexchange.com/questions/1604506/classifying-non-reduced-points-in-noetherian-schemes
            %here this should follow from dim0 and smooth
            points (i.e. points that do not lie in the intersection of several irreducible components) and is contained in $M^*:={M}^*_{0,n}(X,d)\subset \overline{M},$ the locus of automorphism-free maps with smooth source.%note that in the proofs, the sources are reducible: this is because we look only at DAB, not at how it looks like once it intersects Y.
        \end{lemma}

        \begin{proof}
            Repeated use of Kleiman's Theorem on the transversality of the general translate, recalled in \ref{kleiman}, will yield the result. In the language of Kleiman's theorem we set $G=\text{Aut}(X)^n$, and different choices of $A$, $B$, $C$, $f$ and $g$ will lead to the different parts of the statement. Recall the diagram from the statement of the theorem:
            \begin{center}
                \begin{tikzcd}
                && C \arrow[dd,"g"]\\\\
                B\arrow[rr,"f"] && A
                \end{tikzcd}
            \end{center}
            \begin{itemize}
                \item Set $A=X^n$, $B=\underline\Gamma$ and $C=\overline M \setminus M^*$. Let $f$ be the inclusion $B\hookrightarrow A$ and let $g$ be the map $\underline \nu: C\rightarrow A$. Assume that
                $$\text{dim}(\underline\Gamma^{\sigma}\times_{X^n} (\overline M \setminus M^*))=\text{dim} (\underline\Gamma) + \text{dim} (\overline M \setminus M^*) - \text{dim} (X^n).$$
                Together with
                $$\text{dim} (X^n) - \text{dim} (\underline\Gamma) \geq \text{codim}_{X^n}(\underline\Gamma)=\sum_i \text{codim}_X(\Gamma_i)=\text{dim}(\overline M),$$
                this would imply that 
                $$\text{dim}(\underline\Gamma^{\sigma}\times_{X^n} (\overline M \setminus M^*))\leq \text{dim} (\overline M \setminus M^*) - \text{dim}(\overline M) < 0,$$
                since $M^*$ is dense in $\overline M$. As this is impossible, Kleiman's Theorem implies that there exists a dense subset $V_1\subset G$ such that for every $\sigma\in V_1$, the set
                $$\underline\Gamma^{\sigma}\times_{X^n} (\overline M \setminus M^*) = \underline\nu^{-1}(\underline\Gamma^{\sigma})\subset \overline M \setminus M^*$$ is empty, i.e. $\underline\nu^{-1}(\underline\Gamma^{\sigma})\subset M^*.$ Given our choice of $G$ and the generality of $\underline\Gamma$, we can thus assume that  $\underline\nu^{-1}(\underline\Gamma)\subset M^*.$
                    
                \item Set $A=X^n$, $B=\text{Sing}(\underline\Gamma)$ and $C=M^*$ (the choice of $C$ being inspired by the previous point). Again, let $f$ be the inclusion $B\hookrightarrow A$ and $g$ be the map $\underline \nu: C\rightarrow A$. Arguing as for the previous setup, we find that
                $$\text{dim}(\text{Sing}(\underline\Gamma)^{\sigma}\times_{X^n} M^*)\leq \text{dim} (M^*) - \text{codim}_{X^n}(\text{Sing}(\underline\Gamma)) < 0, $$
                where we used $\text{dim}(M^*)\leq \text{dim}(\underline M)$ and $\text{codim}_{X^n}(\text{Sing}(\underline\Gamma))>\text{codim}_{X^n}(\underline\Gamma)$ (this assumes that $\text{codim}_{\underline\Gamma}(\text{Sing}(\underline\Gamma))\ge1$, which will always be the case in our applications).
                Since this is impossible, Kleiman's Theorem again implies the existence of a dense subset $V_2\subset G$ such that for every $\sigma\in V_2$, the set $\underline\nu^{-1}(\text{Sing}(\underline\Gamma)^{\sigma})$ is empty. As above, we may therefore assume that $\underline\Gamma$ is smooth. 
                
                \item Set $A=X^n$, $B=\underline\Gamma \setminus \text{Sing}(\underline\Gamma)$ (due to the previous point) and $C=M^*$. The maps $f$ and $g$ are again given by the inclusion $B\hookrightarrow A$ and $\underline \nu$, respectively. Note that $A$, $B$, and $C$ are smooth. This time, the equality from Kleiman's theorem yields the following:
                $$\text{dim}(\nu^{-1}(B^{\sigma})) \leq \text{dim} (M^*) - \text{codim}_{X^n}(B)) \leq 0,$$
                because $\text{dim}(M^*)\leq\text{dim}(\overline{M})$ and $\text{codim}_{X^n}(B)\geq\text{codim}_{X^n}(\underline\Gamma)=\text{dim} (\overline{M}).$
                Since $(\nu^{-1}(B^{\sigma})$ is non-empty, Kleiman's Theorem implies the existence of a dense subset $V_3\subset G$ such that for every $\sigma\in V_3$, the final inequality above is actually an equality, i.e. $(\nu^{-1}(B^{\sigma})$ has dimension $0$ and is smooth.
            \end{itemize}
            Combining the above results (i.e. taking $\sigma\in V_1\cap V_2\cap V_3$), and invoking again our choice of $G$ and the generality of $\underline\Gamma$, we conclude.
            
        \end{proof}
        
        \begin{remark}\label{onY}
            Note that in the proof of Proposition \ref{n2}, the variety $Y$ does not satisfy the condition $\sum k_i=\text{dim} (\overline{M})$, but instead satisfies $10=\sum k_i=\text{dim} (\overline{M})-1$. However, Kleiman's Theorem again yields the necessary properties. We provide a sketch of the procedure, omitting those parts pertaining to containment in $M^*$ and smoothness discussed in the proof of Lemma \ref{lemma1}.
            
            Let $D$ be a placeholder for either boundary divisor in equation \ref{equivn2}. Arguing as in the proof of Lemma \ref{lemma1}, we find:
            $$0\leq\text{dim}(Y\cap D)\leq\text{dim}(Y)-\text{codim}(D),$$
            where we used that $Y\cap D = Y\times_{\overline M} D$. Since $D$ is a divisor, its codimension is equal to $1$. 
            Another application of Kleiman's Theorem leads to 
            $$\text{dim}(Y) \leq \text{dim}(\overline M)-\text{codim}(\underline\Gamma)=1,$$
            where we used $Y=\underline\nu^{-1}(\underline\Gamma)$.
            Plugging this into the first equation and using that $Y\cap D$ is non-empty, we conclude that $\text{dim}(Y)=1$ and $\text{dim}(Y\cap D)=0$.
        \end{remark}
        
        %is it too early to mention $N_d_A$ or $N_d_B$?
        In the proof of Proposition \ref{n2}, when examining both sides of equation \ref{equivn2}, we counted stable maps mapping $d_A$ (resp. $d_B$) marks, say $p_i$, to corresponding $3 d_A -1$ (resp. $3 d_B -1$) points, say $\Gamma_i$. However, in order to recover $N_{d_A}$ or $N_{d_B}$, we should count rational curves of degree $d_A$ (resp. $d_B$) passing through the points $\Gamma_i$, irrespective of marks or maps. It is not clear that counting both types of objects should yield the same result: if a rational curve meets a $\Gamma_i$ more than once (or with multiplicity greater than 1), then there are several corresponding stable maps sending $p_i$ to $\Gamma_i$. In fact, each intersection of the curve with $\Gamma_i$ corresponds to a possible choice of position for the mark $p_i$ for the stable map. The following lemma will show that such a situation does not arise, hence we can indeed recover $N_{d_A}$ or $N_{d_B}$ by counting stable maps.
        
        \begin{lemma}\label{lemma2}
            If $d>0$ and $\sum k_i=\text{dim} (\overline{M})$, then for every map $\mu\in\underline\nu^{-1}(\underline\Gamma)$ and for every $1\leq i \leq n$, we have $\mu^{-1}(\mu(p_i))=\{p_i\}$ with multiplicity $1$. %\textcolor{red}{this needs d greater than 0: then the preimg is finite and speaking of reduced (points) makes sense.}
            
            Furthermore, if $r=2$, $n=3d-1$ and the $\Gamma_i=Q_i$ are points in general position, then the number of degree $d$ stable maps sending $p_i$ to $Q_i$ for all $1\leq i \leq n$ is equal to $N_d$.
        \end{lemma}

        \begin{remark}
            In the proof of Proposition \ref{n2}, Lemma \ref{lemma2} justifies counting maps instead of curves and legitimizes the appearance of $N_2$ and $N_1$ in our counts of the contributions from each side of equation \ref{equivn2}.
        \end{remark}

        \begin{proof}
            For the first statement, arguing as in the proof of Lemma \ref{lemma1}, we find that $\nu^{-1}(\underline\Gamma)$ lies in the dense subset ${M}^{\circ}_n:={M}^{\circ}_{0,n}(\pr,d)\subset \overline M$ of immersions with smooth source curve and consists of finitely many reduced points. This implies %(\textcolor{red}{how? bc img has only one irreducible component? or rather bc the source curve has only one irr comp.?})
            that for every $1\leq i \leq n$, the subvariety $\mu^{-1}(\mu(p_i))$ is reduced.
            
            Consider the subsets $J_i:=\{\mu \in {M}^{\circ}_n |\; \exists\; q \neq p_i: q\in\mu^{-1}(\mu(p_i))\} \subset {M}^{\circ}_n$ and $Q_{i,j}:=\{\mu\in M:={M}_{0,n}(\pr,d)|\mu(p_i)=\mu(p_j)\} \subset M$ for $1\leq i,j\leq n$. 
            By adding an extra mark $p_0$, we can look at the space ${M}^{\circ}_{n+1}:={M}^{\circ}_{0,n+1}(\pr,d)$ and the forgetful map 
            $\epsilon: {M}^{\circ}_{n+1} \rightarrow {M}^{\circ}_{n}$ (forgetting $p_0$).
            Combining the following two points yields the first statement.
            
            \begin{itemize}
                \item We have $\text{codim}(J_i)>0$. Indeed, consider the subsets $Q_{i,j}$ defined above. Here it is enough to consider $M$ instead of $\overline M$ due to Lemma \ref{lemma1}. Clearly $\epsilon$ maps $Q_{i,0}$ into $J_i$. Moreover, it does so surjectively: for any $\mu\in J_i$ with $p_i\neq  q\in\mu^{-1}(\mu(p_i))$, add the extra mark $p_0=q$ to $\mu$. The resulting map in ${M}_{0,n+1}(\pr,d)$ gets mapped to $\mu$ by $\epsilon$. Thus $J_i=\epsilon(Q_{i,0})$. Hence $\text{codim}(Q_{i,0})=r$ (see \ref{qijcodim}) implies that $\text{codim}(J_i)\geq r-1\geq1>0$ (since by forgetting a mark, $\epsilon$ effectively reduces the dimension by 1). 
                \item If $\text{codim}(J_i)>0$, then $\mu^{-1}(\mu(p_i))=\{p_i\}$. Indeed, applying Kleiman's Theorem to the data $A={M}^{\circ}_n$, $B=J_i$, $C=Q_{i,0}$, $f=\text{inc}:J_i\xhookrightarrow{}{M}^{\circ}_n$ and $g=\epsilon|_{Q_{i,0}}$, we get the following inequality: 
                $$\text{dim}((\epsilon|_{Q_{i,0}})^{-1}) \leq \text{dim}(Q_{i,0}) - \text{codim}_{{M}^{\circ}_n}(J_i).$$
                But this is a contradiction unless $\text{codim}_{{M}^{\circ}_n}(J_i)=0$, since $(\epsilon|_{Q_{i,0}})^{-1}(J_i)=Q_{i,0}$.
            \end{itemize}

            Now the second statement follows immediately, since we've shown that $\mu^{-1}(Q_i)=\mu^{-1}(\mu(p_i))=\{p_i\}$ and the intersection has multiplicity $1$.
            
        \end{proof}

        A sequel of sorts to these Lemmas can be found in Lemma \ref{355}. We omit it here since it will only be applied in section \ref{th2quantum}.

        For additional intuition on these ideas see Proposition \ref{n3}, where we prove that $N_3=12$.
        
        \begin{theorem}[Kontsevich's Formula for $\ptwo$]\label{nd}
            For any $d\geq1$, the number $N_d$ of rational degree $d$ curves passing through $3d-1$ points in general position in $\ptwo$ satisfies the following recursive equation:
                \begin{align}\label{kontsformp2}
                    \begin{split}
                        N_d \; \; & + \sum_{\substack{d_A+d_B=d \\ d_A\geq1,d_B\geq1}} \binom{3d-4}{3d_A-1} d_A^2 N_{d_A} \cdot d_A d_B \cdot N_{d_B} \\
                        &=\sum_{\substack{d_A+d_B=d \\ d_A\geq1,d_B\geq1}} \binom{3d-4}{3d_A-2} d_A N_{d_A} \cdot d_A d_B \cdot d_B N_{d_B}.
                    \end{split}
                \end{align}
            Since $N_1=1$ by Fact \ref{n1}, we can thus compute any $N_d$.
        \end{theorem}

        \begin{proof} %analogous to n3
            The overall structure of the proof is analogous to that of the proofs of Propositions \ref{n2} and \ref{n3}. 
            We now work in $\overline{M}_{0,3d}(\ptwo,d)$ with marks $m_1$, $m_2$, $p_1$, \dots, $p_{3d-2}$ and consider as our data two lines $L_1$ and $L_2$ and $3d-2$ points $Q_1, \dots, Q_{3d-2}$ in general position in $\ptwo$.
            We denote by $Y$ the subvariety of $\overline{M}_{0,3d}(\ptwo,d)$ consisting of stable maps $(C;m_1, m_2, p_1, \dots, p_{3d-2};\mu)$ sending $m_1$ to $L_1$, $m_2$ to $L_2$ and $p_i$ to $Q_i$ for $1\leq i\leq 3d-2$, that is:
            $$Y=\nu_{m_1}^{-1}(L_1) \cap \nu_{m_2}^{-1}(L_2) \cap \nu_{p_1}^{-1}(Q_1) \cap \dots \cap \nu_{p_{3d-2}}^{-1}(Q_{3d-2}).$$
            By Lemmas \ref{lemma1} and \ref{lemma2} and Remark \ref{onY}, $Y$ is a curve and its intersection with each boundary divisor is transversal and takes place in $\overline{M}_{0,3d}^*(\ptwo,d)$.
            The fundamental boundary relation (\ref{bdry}) again yields the equality:
            \begin{equation}\label{equivnd}
                \abs{Y\cap D(m_1,m_2|p_1,p_2)} = \abs{Y\cap D(m_1,p_1|m_2,p_2)},
            \end{equation}
            where $D(m_1,m_2|p_1,p_2)$ and $D(m_1,p_1|m_2,p_2)$ represent sums of the form $$\sum_{\substack{m_1,m_2\in A \\ p_1,p_2\in B \\d_A+d_B=d}}D(A,B;d_A,d_B) \; \text{ and} \sum_{\substack{m_1,p_1\in A \\ m_2,p_2\in B \\d_A+d_B=d}}D(A,B;d_A,d_B),$$ respectively.
            We proceed to an examination of the right and left hand sides of equality  \ref{equivnd}. The $3d-4$ spare marks $p_3$, \dots, $p_{3d-2}$ remain to be distributed between $C_A$ and $C_B$.\\
            
            LHS:\\
            Ad $(d,0)$:
            This pair does not contribute to the intersection as $Q_1=\mu(p_1)=\mu(p_2)=Q_2$ would contradict the genericity of our data. \\
        
            Ad $(0,d)$: %by rmk in notes, total codimension is correct.
            In this case we have $\{\mu(m_1)\}=\{\mu(m_2)\}=L_1\cap L_2$. We apply Lemmas \ref{lemma1} and \ref{lemma2} with ${\Gamma}_1=\mu(m_1)$, $\Gamma_2=Q_1$, \dots, $\Gamma_{3d-1}=Q_{3d-2}$ and deduce a contribution of $N_d$, the number we seek. \\

            Ad $(d_A,d_B)$ with $d_A,d_B>0$:
            We examine first the distribution of the spare $p_i$'s, then the positions of the $p_i$'s and $m_i$'s, and finally the position of the intersection node $x$.
            \begin{itemize}
                %number and distribution of marks
                \item We must add $3d_A-1$ of the spare marks to $C_A$ and $3d_B-3$ to $C_B$ (note that these are all the spare marks: $3d_A-1+3d_B-3=3(d_A+d_B)-4$ and $d=d_A+d_B$). Indeed, any other distribution would contradict the genericity of the data: putting more than $3d_A-1$ marks on $C_A$ (or equivalently, fewer than $3d_B-3$ on $C_B$) would force more than $3d_A-1$ of the $Q_i$'s onto the degree $d_A$ curve $\mu(C_A)$. Similarly, putting more than $3d_B-3$ marks on $C_B$ (or equivalently, fewer than $3d_A-1$ on $C_B$) would force more than $3d_B-1$ of the $Q_i$'s onto the degree $d_B$ curve $\mu(C_B)$.
                There are thus $\binom{3d-4}{3d_A-1}$ ways to distribute the spare marks onto $C_A$ and $C_B$.
                %positions of marks
                \item The positions of the points $p_1$, \dots, $p_{3d-2}$ on their respective twigs are now uniquely determined. Indeed, for each $i\in\{1, \dots,3d-2\}$, the point $p_i$ is the only element of $\mu^{-1}(Q_i)$ by Lemma \ref{lemma2}, and the order of the $p_i$'s is fixed because the maps $C_A \rightarrow \ptwo$ and $C_B \rightarrow \ptwo$ are each birational onto their respective image. 
                Setting $\Gamma_i=p_i$ for the $3d_A-1$ spare $p_i$ of $C_A$ and applying Lemmas \ref{lemma1} and \ref{lemma2} yields a contribution of $N_{d_A}$. The analogous procedure on the $3d_B-1$ marks of $C_B$ yields a contribution of $N_{d_B}$. 

                The position of the mark $m_1$ is determined by the intersection $\mu(C_A)\cap L_1$. By B\'ezout's Theorem, there are $d_A\cdot1=d_A$ possible positions for $m_1$. The same argument applies to $m_2$, where $\mu(C_A)\cap L_2$ yields another contribution of $d_A$, so the total contribution from the positions of $m_1$ and $m_2$ is equal to $d_A^2$.
                %gluing
                \item Finally, we account for the different ways of 'gluing' the images of $C_A$ and $C_B$ together: by B\'ezout's Theorem, there are $d_A\cdot d_B$ possible choices for $\mu(x)$.
            \end{itemize} 
            The contribution from $(d_A,d_B)$ with $d_A,d_B>0$ thus looks as follows: 
            $$\binom{3d-4}{3d_A-1} d_A^2 N_{d_A}\cdot d_A d_B \cdot N_{d_B}.$$
            The binomial coefficient comes from the distribution of the spare marks, the factor $d_A^2 N_{d_A}$ comes from $C_A$, the factor $d_A d_B$ comes from gluing together both twigs,  and the factor $N_{d_B}$ comes from $C_B$.

            Summing over all the above pairs of degrees, we obtain the following total contribution from the left hand side:
            $$N_d \; \; + \sum_{\substack{d_A+d_B=d \\ d_A\geq1,d_B\geq1}} \binom{3d-4}{3d_A-1} d_A^2 N_{d_A} \cdot d_A d_B \cdot N_{d_B}.$$\\

            RHS:\\
            Ad $(0,d)$, $(d,0)$: If $(d_A,d_B)=(0,d)$, then $Q_1=\mu(p_1)=\mu(m_1)\in L_1$. Similarly, if $(d_A,d_B)=(d,0)$, then $Q_2\in L_2$. Both conclusions contradict the genericity of the data, so neither case contributes to the right hand side.\\
            
            Ad $(d_A,d_B)$ with $d_A,d_B>0$:
            \begin{itemize}
                %number and distribution of marks
                \item We must add $3d_A-2$ of the spare marks to $C_A$ and $3d_B-2$ to $C_B$ (note that these are all the spare marks: $3d_A-2+3d_B-2=3(d_A+d_B)-4$ and $d=d_A+d_B$). Again, any other distribution would contradict the genericity of the data: putting more than $3d_A-2$ marks on $C_A$ (or equivalently, less than $3d_B-2$ on $C_B$) would force more than $3d_A-1$ of the $Q_i$'s onto the degree $d_A$ curve $\mu(C_A)$. The argument for $C_B$ is symmetric.
                There are thus $\binom{3d-4}{3d_A-2}$ ways to distribute the spare marks onto $C_A$ and $C_B$.
                %positions of marks
                \item The positions of the points $p_1$, \dots, $p_{3d-2}$ on their respective twigs are now uniquely determined, just as on the left hand side.
                Setting $\Gamma_i=p_i$ for the $3d_A-1$ marks $p_i$ of $C_A$ and applying Lemmas \ref{lemma1} and \ref{lemma2} yields a contribution of $N_{d_A}$. The analogous procedure on the $3d_B-1$ marks of $C_B$ yields a contribution of $N_{d_B}$. 

                The positions of $m_1$ and $m_2$ are determined by the intersections $\mu(C_A)\cap L_1$ and $\mu(C_B)\cap L_2$, respectively. By B\'ezout's Theorem, we get respective contributions of $d_A\cdot1=d_A$ and $d_B\cdot1=d_B$.
                %gluing
                \item As on the left hand side, there are $d_A\cdot d_B$ ways of 'gluing' the images of $C_A$ and $C_B$ together.
            \end{itemize} 
            On the right hand side, the contribution from $(d_A,d_B)$ with $d_A,d_B>0$ is therefore the following:
            $$\binom{3d-4}{3d_A-2} d_A N_{d_A}\cdot d_A d_B \cdot d_B N_{d_B}.$$
            Similarly to the left hand side, the binomial coefficient comes from the distribution of the spare marks, the factor $d_A N_{d_A}$ comes from $C_A$, the factor $d_A d_B$ comes from gluing together both twigs,  and the factor $d_B N_{d_B}$ comes from $C_B$.

            Summing over all pairs of degrees as above yields the following total contribution from the right hand side:
            $$\sum_{\substack{d_A+d_B=d \\ d_A\geq1,d_B\geq1}} \binom{3d-4}{3d_A-2} d_A N_{d_A}\cdot d_A d_B \cdot d_B N_{d_B}.$$
            
            Equating the total contributions from both sides of equality \ref{equivnd} yields the desired formula.
            
        \end{proof}

        We collect the values of $N_d$ for $1\le d \le 12$ in a table taken from \cite{Itzy}.
              
        \begin{table}[h!]
            \centering
            \begin{tabular}{|c|c|c|}
                \hline
                d  & $3d-1$ & $N_d$                        \\ \hline
                1  & 2    & 1                           \\
                2  & 5    & 1                           \\
                3  & 8    & 12                          \\
                4  & 11   & 620                         \\
                5  & 14   & 87304                       \\
                6  & 17   & 26312976                    \\
                7  & 20   & 14616808192                 \\
                8  & 23   & 13525751027392              \\
                9  & 26   & 19385778269260800           \\
                10 & 29   & 40739017561997799680        \\
                11 & 32   & 120278021410937387514880    \\
                12 & 35   & 482113680618029292368686080 \\ \hline
            \end{tabular}
        \end{table}

    \section{Kontsevich's Formula for $\ponex$}\label{kfp1x}
    %structure as above. avoid repetition. 
    %NB: rules, bidegree etc already explained in theory 1.
    
    %discussion on why we need 2d+2e-1 pts
    %first as for ptwo. then mention shafarevic aliter (and cite book).
        In $\ponex$, the presence of bidegrees adds a new dimension to the problem. As a result, some new and interesting properties of the desired numbers arise.
        As in the $\pr$ case, we first determine the necessary number of point conditions for fixing a finite number of irreducible curves of bidegree (d,e) in $\ponex$. We will always assume that $d+e\ge1$.
        
        Such a curve is determined by a pair of homogeneous polynomials $f(x,y)$ and $g(x,y)$ of respective degrees $d$ and $e$. We can write out $f$ as $a_0 x^d + a_1 x^{d-1}y +\dots+a_{d-1}xy^{d-1}+a_d y^d$, a polynomial $d+1$ coefficients. The same of course works for $g$, this time with $e+1$ coefficients. There are thus a total of $(d+1)(e+1)$ total coefficients up to scaling, so the space of irreducible curves of bidegree $(d,e)$ in $\ponex$ is of dimension ${(d+1)(e+1)-1}$. 
        %TODO perhaps remark that this is equal to dim M_0,0(p1xp1,(d,e))

        But for such a curve $C$ of genus $g_C$, it holds that:
        $$(d+1)(e+1)-1=(d-1)(e-1)+2(d+e)-1=g_C+2d+2e-1,$$
        see \ref{genusponex}. The dimension of the subspace of rational (i.e. $g_C=0$) curves is thus equal to $2d+2e-1$, so we need $2d+2e-1$ codimension $1$ conditions to obtain a finite number of such curves. 
    
        Hence we shall define $N_{(d,e)}$ as the number of rational curves of bidegree $(d,e)$ passing through $2d+2e-1$ points in general position in $\ponex$.
        
        \begin{remark}
            If more than $2d+2e-1$ points lie on a bidegree $(d,e)$ rational curve, then they are not in general position.
        \end{remark}

        The following fact is clear from the definition of horizontal and vertical rules.
        \begin{fact}\label{n10} %reference Euclid?
            There are exactly $N_{(0,1)}=1$ horizontal rule and $N_{(1,0)}=1$ vertical rule passing through a general point in $\ponex$.
        \end{fact}

        \begin{lemma}\label{sym}
            For all $d,e$ such that $d+e\geq1$, it holds that $N_{(d,e)}=N_{(e,d)}$.
        \end{lemma}

        \begin{proof}
            By definition, $N_{(d,e)}$ is the number of rational bidegree $(d,e)$ curves passing through $2d+2e-1$ points in general position in $\ponex$. 
            
            Let $C$ be such a $(d,e)$-curve given by the pair of homogeneous polynomials $(f(x,y),g(x,y))$, with $\text{deg}(f)=d$, $\text{deg}(g)=e$.
            Define $\Tilde{C}$ as the curve given by the pair $(g(x,y),f(x,y))$. 
            
            Clearly, $Q:=([a:b],[c:d])\in C$ if and only if $\Tilde{Q}:=([c:d],[a:b])\in \Tilde{C}$.
            Thus the number of $(d,e)$-curves passing through general points $Q_1$, \dots, $Q_{2d+2e-1}$ is equal to the number of $(e,d)$-curves passing through the corresponding general points $\Tilde{Q}_1$, \dots, $\Tilde{Q}_{2d+2e-1}$.
            
            % $F(x,y,z,w)$, with $F$ a homogeneous bidegree $(d,e)$ polynomial. 
            % Define $\Tilde{C}$ as the curve given by 
        \end{proof}

        Consider again the language of Lemmas \ref{lemma1} and \ref{lemma2} and their proofs from the $\ptwo$ case. With some reassignments, i.e. $\overline M := \overline{M}_{0,2d+2e}(\ponex,(d,e))$, $X:=\ponex$ etc., we arrive at analogous statements for the $\ponex$ case. Note that,  similarly to before, $\sum k_i=\text{dim} (\overline{M})$ holds if $n=2d+2e-1$ and the $\Gamma_i$ are points (i.e. $\text{codim}(\Gamma_i)=k_i=2$ for $1\leq i \leq n$). Indeed:
        $$\sum k_i=2(2d+2e-1) = \text{dim} (\overline{M}_{0,2d-2e}(\ptwo,(d,e))=\text{dim} (\overline{M}),$$
        where we used \ref{dimmbarone}.     
        
        Instead of repeating the statements here, we shall keep referring to the original Lemmas \ref{lemma1} and \ref{lemma2}, keeping the necessary adjustments of terminology in mind. 
        For instance, after adjustment, the second part of Lemma \ref{lemma2} states that if $n=2d+2e-1$ and the $\Gamma_i=Q_i$ are points in general position, then the number of bidegree $(d,e)$ stable maps sending $p_i$ to $Q_i$ for all $1\leq i \leq n$ is equal to $N_{(d,e)}$.

        \begin{proposition}\label{nd0} 
            For any $d>1$, it holds that $N_{(0,d)}=0$. 
        \end{proposition}

        \begin{proof}
            By definition, $N_{(0,d)}$ is the number of bidegree $(0,d)$ curves through $2d-1$ points in general position in $\ponex$. %\textcolor{red}{what is (0,d)? is interpretation of (d,e)-curves via zero sets even valid? what is the link? SUSPICION. verify and elaborate the next statement}%understand the link between classes and zero sets of polynomials
            But a bidegree $(0,d)$ curve (a horizontal rule of multiplicity $d$) is a selection of $d$ horizontal rules. These are determined uniquely by $d$ general points (their 'heights' on different vertical rules). Since $2d-1>d$ if and only if $d>1$, $2d-1$ points on a bidegree $(0,d)$ curve cannot be in general position for $d>1$. Thus $N_{(0,d)}=0$ as claimed.  
            
        \end{proof}

        \begin{theorem}[Kontsevich's Formula for $\ponex$]\label{nde}
                For any $d,e$ with $d+e\geq1$, the number $N_{(d,e)}$ of rational bidegree $(d,e)$ curves passing through $2d+2e-1$ points in general position in $\ponex$ satisfies the following recursive equation:
                    \begin{align}\label{kontsponex}
                        \begin{split}
                            N_{(d,e)} \; \; & + \sum_{\substack{d_A+d_B=d \\ e_A+e_B=e}} \binom{2d+2e-4}{2d_A+2e_A-1} d_A e_A N_{(d_A,e_A)} \cdot (d_A e_B + e_A d_B) \cdot N_{(d_B,e_B)} \\
                            &=\sum_{\substack{d_A+d_B=d \\ e_A+e_B=e}} \binom{2d+2e-4}{2d_A+2e_A-2} d_A N_{(d_A,e_A)} \cdot (d_A e_B + e_A d_B) \cdot e_B N_{(d_B,e_B)},
                        \end{split}
                    \end{align}
                    where the sums are over pairs $(d,e)$ such that $d+e>0$.
                Since $N_{(1,0)}=N_{(0,1)}=1$ by Fact \ref{n10} and $N_{(d,0)}=N_{(0,d)}=0$ for $d>1$ by Proposition \ref{nd0}, we can thus compute any $N_{(d,e)}$.
        \end{theorem}
        Note that since $N_{(d,0)}=N_{(0,d)}=0$ for $d>1$, one could also restrict the sums to pairs $(d,e)$ with $d\ne0$ and $e\ne0$ and require only the two initial values $N_{(1,0)}=N_{(0,1)}=1$.

        \begin{remark}
            More precisely, let us arrange the sought after numbers in an arbitrarily large matrix:
            $$\begin{pmatrix}
                \times & N_{(0,1)} & N_{(0,2)} & \dots\\
                N_{(1,0)} & N_{(1,1)} & N_{(1,2)} & \dots\\
                N_{(2,0)} & N_{(2,1)} & N_{(2,2)} & \dots\\
                \dots& \dots& \dots& \dots
            \end{pmatrix}.
            $$
            Note that $N_{(0,0)}$ is not defined, hence the $\times$ in the upper left entry. The sums 
            %of type $\sum_{\substack{d_A+d_B=d \\ e_A+e_B=e}}$ 
            in formula \ref{kontsponex} are over pairs with $(d_A,e_A)\neq(0,0)\neq(d_B,e_B)$, but even if this condition was erroneously omitted, the factors $d_A$ (in the case $(d_A,e_A)=(0,0)$) and $(d_A e_B + e_A d_B)$ (in the case $(d_B,e_B)=(0,0)$) would annihilate any summands containing the undefined number $N_{(0,0)}$.
            
            The aim of Theorem \ref{nde} is to recursively describe all of the entries of this matrix. By Lemma \ref{sym}, it would be sufficient to determine the diagonal and either the upper or lower triangular part. 
            
            In order to find any unknown number $N_{(d,e)}$, formula \ref{kontsponex} requires the numbers $N_{(d_A,e_A)}$ and $N_{(d_B,e_B)}$ for all $(d_A,e_A)$ and $(d_B,e_B)$ such that $d_A+d_B=d$ and $e_A+e_B=e$, i.e. the following submatrix: 
             $$\begin{pmatrix}
                \times & N_{(0,1)} & \dots &  \dots\\
                N_{(1,0)} & N_{(1,1)}  & \dots &  \dots\\
                \dots& \dots& \dots & \dots \\
                \dots& \dots& \dots& N_{(d,e)}
            \end{pmatrix}
            $$
            Due to the nature of the recursion, it is therefore sufficient to know the first line and first column of the initial matrix, and these are identical by Lemma \ref{sym}. Our knowledge of $N_{(0,d)}$ for any $d\geq0$ is thus enough to compute any $N_{(d,e)}$. 
        \end{remark}
    
        \begin{proof} 
            We now work in $\overline{M}_{0,2d+2e}(\ponex,(d,e))$ with marks $m_1$, $m_2$, $p_1$, \dots, $p_{2d+2e-2}$ and consider as our data a horizontal rule $L_1$, a vertical rule $L_2$ and $2d+2e-2$ points $Q_1, \dots, Q_{2d+2e-2}$ in general position in $\ponex$ (lying on neither of the rules).
            Denote by $Y$ the subvariety of $\overline{M}_{0,2d+2e}(\ponex,(d,e))$ consisting of stable maps $(C;m_1, m_2, p_1, \dots, p_{2d+2e-2};\mu)$ sending $m_1$ to $L_1$, $m_2$ to $L_2$ and $p_i$ to $Q_i$ for $1\leq i\leq 2d+2e-2$:
            $$Y=\nu_{m_1}^{-1}(L_1) \cap \nu_{m_2}^{-1}(L_2) \cap \nu_{p_1}^{-1}(Q_1) \cap \dots \cap \nu_{p_{2d+2e-2}}^{-1}(Q_{2d+2e-2}).$$
            By Lemmas \ref{lemma1} and \ref{lemma2}, $Y$ is a curve and its intersection with each boundary divisor is transversal and takes place in $\overline{M}_{0,2d+2e}^*(\ponex,(d,e))$.
            The fundamental boundary relation (\ref{bdry}) yields the equality :
            \begin{equation}\label{equivnde}
                Y\cap D(m_1,m_2|p_1,p_2)\equiv Y\cap D(m_1,p_1|m_2,p_2),
            \end{equation}
            where $D(m_1,m_2|p_1,p_2)$ and $D(m_1,p_1|m_2,p_2)$ represent sums of the form 
            $$\sum_{\substack{m_1,m_2\in A \\ p_1,p_2\in B\\d_A+d_B=d  \\e_A+e_B=e}}D\bigr(A,B;(d_A,e_A),(d_B,e_B)\bigl) \; \text{ and} \sum_{\substack{m_1,p_1\in A \\ m_2,p_2\in B \\d_A+d_B=d  \\e_A+e_B=e}}D\bigr(A,B;(d_A,e_A),(d_B,e_B)\bigl),$$
            respectively.
            We proceed to an examination of the right and left hand sides of equality  \ref{equivnde}. The $2d+2e-4$ spare marks $p_3$, \dots, $p_{2d+2e-2}$ remain to be distributed between $C_A$ and $C_B$.\\
            
            LHS:\\
            Ad $(d_A,e_A)=(0,0)$, $(d_B,e_B)=(d,e)$:
            All spare marks must go on $C_B$, otherwise there exists $3\leq i\leq2d+2e-3$ with $Q_i\in L_1\cap L_2$, contradicting genericity. Applying Lemmas \ref{lemma1} and \ref{lemma2} to the points ${\Gamma}_1=\mu(m_1)$, $\Gamma_2=Q_1$, \dots, $\Gamma_{2d+2e-1}=Q_{2d+2e-1}$ yields a contribution of $N_{(d,e)}$, the number we seek.
            \\

            Ad any other $(d_A,d_B)$, $(e_A,e_B)$ with $d_A+d_B=d$, $e_A+e_B=e$:
            % We examine first the distribution of the spare $p_i$'s, then the positions of the $pi$'s and $m_i$'s, and finally the position of the intersection node $x$.
            \begin{itemize}
                %number and distribution of marks
                \item We must add $2d_A+2e_A-1$ of the spare marks to $C_A$ (since it has no $p_i$'s) and $2d_B+2e_B-3$ to $C_B$ (since it has two $p_i$'s). Any other distribution contradicts the genericity of the data. These are all the spare marks: $2d_A+2e_A -1 + 2d_B+2e_B-3=2(d_A+d_B)+2(e_A+e_B)-4$.
                There are thus $\binom{2d+2e-4}{2d_A+2e_A-1}$ ways to distribute the spare marks onto $C_A$ and $C_B$.
                %positions of marks
                \item Once distributed, the positions of the points $p_1$, \dots, $p_{2d+2e-2}$ on their respective twigs are uniquely determined. Indeed, for each $i\in\{1, \dots,2d+2e-2\}$, the point $p_i$ is the only element of $\mu^{-1}(Q_i)$ by Lemma \ref{lemma2}, and the order of the $p_i$'s is fixed because the maps $C_A \rightarrow \ponex$ and $C_B \rightarrow \ponex$ are each birational onto their respective image. 
                Setting $\Gamma_i=p_i$ for the $2d_A+2e_A-1$ marks $p_i$ of $C_A$ and applying Lemmas \ref{lemma1} and \ref{lemma2} yields a contribution of $N_{(d_A,e_A)}$. The analogous procedure on the $2d_B+2e_B-1$ marks of $C_B$ yields a contribution of $N_{(d_B,e_B)}$. 
    
                The position of the mark $m_1$ is determined by the intersection $\mu(C_A)\cap L_1$. By B\'ezout's Theorem, there are $(d_A,e_A)\circ(0,1)=d_A$ possible positions for $m_1$. For $m_2$, $\mu(C_A)\cap L_2$ yields a factor $(d_A,e_A)\circ(1,0)=e_A$, so the total contribution from the positions of $m_1$ and $m_2$ is equal to $d_A e_A$.
                
                %gluing
                \item Finally, we account for the different ways of 'gluing' the images of $C_A$ and $C_B$ together: by B\'ezout's Theorem, there are $(d_A,e_A)\circ(d_B,e_B)=d_A e_B + e_A d_B$ possible choices for $\mu(x)$.
            \end{itemize} 
            The contribution from these $(d_A,d_B)$, $(e_A,e_B)$ is thus equal to the following: 
            $$\binom{2d+2e-4}{2d_A+2e_A-1} d_A e_A N_{(d_A,e_A)} \cdot (d_A e_B + e_A d_B) \cdot N_{(d_B,e_B)}.$$
            The binomial coefficient comes from the distribution of the spare marks, the factor $d_A e_A N_{(d_A,e_A)}$ comes from $C_A$, the factor $d_A e_B + e_A d_B$ comes from gluing together both twigs,  and the factor $N_{(d_B,e_B)}$ comes from $C_B$.

            Summing over all possible pairs of bidegrees, we obtain the claimed total contribution from the left hand side.\\

            RHS:\\          
            Ad $(d_A,d_B)$, $(e_A,e_B)$ with $d_A+d_B=d$, $e_A+e_B=e$:
            \begin{itemize}
                %number and distribution of marks
                \item We must add $2d_A+2e_A-2$ of the spare marks to $C_A$ and $2d_B+2e_B-2$ to $C_B$. Again, any other distribution would contradict the genericity of the data and these are all the spare marks.
                As such there are $\binom{2d+2e-4}{2d_A+2e_A-2}$ possible distributions.
                %positions of marks
                \item After distribution, the positions of the points $p_1$, \dots, $p_{2d+2e-2}$ on their respective twigs are uniquely determined, just as on the left hand side.
                Applying Lemmas \ref{lemma1} and \ref{lemma2} yields a contribution of $N_{(d_A,e_A)}$ on $C_A$ and one of $N_{(d_B,e_B)}$ on $C_B$. 

                The positions of $m_1$ and $m_2$ are determined by the intersections $\mu(C_A)\cap L_1$ and $\mu(C_B)\cap L_2$, respectively. By B\'ezout's Theorem, we get respective factors of  $(d_A,e_A)\circ(0,1)=d_A$ and  $(d_B,e_B)\circ(1,0)=e_B$.
                % %gluing
                \item Just like on the left hand side, there are $(d_A,e_A)\circ(d_B,e_B)=d_A e_B + e_A d_B$ ways of 'gluing' the images of $C_A$ and $C_B$ together.
            \end{itemize} 
            On the right hand side, the contribution from such a pair of bidegrees $(d_A,d_B)$, $(e_A,e_B)$ is thus given by the following product: 
            $$\binom{2d+2e-4}{2d_A+2e_A-2} d_A N_{(d_A,e_A)} \cdot (d_A e_B + e_A d_B) \cdot e_B N_{(d_B,e_B)}.$$
            Similarly to the left hand side, the binomial coefficient comes from the distribution of the spare marks, the factor $d_A N_{(d_A,e_A)}$ comes from $C_A$, the factor $d_A e_B + e_A d_B$ comes from gluing together both twigs,  and the factor $e_B N_{(d_B,e_B)}$ comes from $C_B$.
            
            Summing over all pairs of bidegrees and equating the total contributions from both sides of equality  \ref{equivnd} yields the desired formula.
            %this does not fit the page
            % \begin{align}
            %     \begin{split}
            %         N_{(d,e)} \; \; & + \sum_{\substack{d_A+d_B=d \\ e_A+e_B=e}} \binom{2d+2e-4}{2d_A+2e_A-1} (d_A,e_A)\circ(0,1) \cdot (d_A,e_A)\circ(1,0) N_{(d_A,e_A)} \cdot (d_A,e_A)\circ(d_B,e_B) \cdot N_{(d_B,e_B)} \\
            %         &=\sum_{\substack{d_A+d_B=d \\ e_A+e_B=e}} \binom{2d+2e-4}{2d_A+2e_A-2} (d_A,e_A)\circ(0,1) N_{(d_A,e_A)} \cdot (d_A,e_A)\circ(d_B,e_B) \cdot (d_B,e_B)\circ(1,0) N_{(d_B,e_B)},
            %     \end{split}
            % \end{align}
            % which strongly resembles the formula for the $\ptwo$ case from Theorem \ref{nd} and after writing out the products of bidegrees becomes precisely the claimed formula.
            
        \end{proof}

        As an example, we calculate the numbers $N_{(d,1)}$ for all $d\geq0$. Of course we could also have deduced them 'directly' by proceeding as in the above proof.
        
        \begin{corollary}\label{nd1} 
            For any $d\geq0$, there is exactly $N_{(d,1)}=1$ rational $(d,1)$-curve passing through $2d+1$ points in general position in $\ponex$.
        \end{corollary}
        
        \begin{proof}
            For $d=0$, this is true by Fact \ref{n10} and Lemma \ref{sym}. 
            We proceed by induction on $d$.

            Pick an arbitrary $d>1$. Then by Theorem \ref{nde}, we have:
            \begin{align*}
                \begin{split}
                    N_{(d,1)} \; \; & + \sum_{\substack{d_A+d_B=d \\ e_A+e_B=1}} \binom{2d-2}{2d_A+2e_A-1} d_A e_A N_{(d_A,e_A)} \cdot (d_A e_B + e_A d_B) \cdot N_{(d_B,e_B)} \\
                    &=\sum_{\substack{d_A+d_B=d \\ e_A+e_B=1}} \binom{2d-2}{2d_A+2e_A-2} d_A N_{(d_A,e_A)} \cdot (d_A e_B + e_A d_B) \cdot e_B N_{(d_B,e_B)}.
                \end{split}
            \end{align*}
            We examine the sums on both sides.
            Note that in order to satisfy $e_A+e_B=1$, either $e_A=1$ and $e_B=0$ or vice versa. The right hand sum has a factor $e_A$, so only the summands with $e_A=1$ and $e_B=0$ survive here. On the left hand side, there is a factor $e_B$, so only the summands with $e_A=0$ and $e_B=1$ remain. Keeping this in mind, we get the following: 
            \begin{align*}
                \begin{split}
                    N_{(d,1)} \; \; & + \sum_{d_A+d_B=d} \binom{2d-2}{2d_A+1} d_A N_{(d_A,1)} \cdot d_B \cdot N_{(d_B,0)} \\
                    &=\sum_{d_A+d_B=d} \binom{2d-2}{2d_A-2} d_A N_{(d_A,0)} \cdot d_A \cdot N_{(d_B,1)}.
                \end{split}
            \end{align*}
            By the induction hypothesis, the factors $N_{(d_A,1)}$ and $N_{(d_B,1)}$ are both equal to one. Moreover, using Fact \ref{n10} and Proposition \ref{nd0}, we deduce that $N_{(d_B,0)}=\delta_{d_B,1}$ on the left hand side and $N_{(d_A,0)}=\delta_{d_A,1}$ on the right, where $\delta_{i,j}$ is the usual Kronecker delta. Plugging this into the equation yields:
            $$N_{(d,1)} + \binom{2d-2}{2d-1} = \binom{2d-2}{0}.$$
            Evaluating the binomial coefficients we conclude that $N_{(d,1)} = 1$ for all $d\geq0$.
            
        \end{proof}

        For the sake of illustration, we borrow another table from \cite{Itzy}, this time representing the numbers $N_{(d,e)}$ for $0\le d\le 3$.

        \begin{center}
            \renewcommand\arraystretch{1.3}
            \setlength\doublerulesep{0pt}
            \begin{tabular}{|c||*{4}{c|}}
                \hline\backslashbox{d}{e}  & 0        & 1 & 2  & 3 \\
                \hline\hline
                0 & $\times$ & 1 & 0  & 0   \\ 
                \hline
                1 & 1        & 1 & 1  & 1 \\ 
                \hline
                2 & 0        & 1 & 12 & 96 \\ 
                \hline
               3 & 0        & 1 & 96 & 3510 \\ 
                \hline
                
            \end{tabular}
        \end{center}

\chapter{Extending the Theory}\label{th2}%titles to be adapted
%motivate. basically the PD of our previous approach!
%coefficients in \Q
%do as much as possible over general X, with pr and ponex as special cases.
In this section, we generalize the numbers $N_d$ and $N_{d,e}$ into the theory of Gromov-Witten invariants. 
This will allow us to involve useful techniques from combinatorics and cohomological algebra, some of which we shall recall momentarily. The more general machinery will lead to the definition of a new structure, the so called quantum cohomology ring.

    \section{Sequences and Generating Functions}
        %combinatorics: series and recursions
        We briefly discuss some handy combinatorial tools.
        We begin by introducing generating functions associated to sequences of numbers. These are data structures that will allow us to express properties of and relationships between the numbers of a sequence in a concise and elegant way.
        \begin{definition}%generating function
            Let $(N_k)_{k=0}^\infty$ be a sequence of numbers. The exponential generating function of this sequence is the following formal power series:
            $$F(x)=\sum_{k=0}^\infty\frac{x^k}{k!}N_k.$$
            
        \end{definition}
        Note that the degree $k$ term recovers the number $N_k$.
        
        \begin{definition}\label{prodrule}
        The product of two exponential power series $\sum f_n \frac{x^n}{n!}$ and $\sum g_n \frac{x^n}{n!}$ is defined as follows:
            \begin{equation*}
                \left( \sum_{n\geq0}f_n \frac{x^n}{n!}\right)\cdot \left(\sum_{n\geq0} g_n \frac{x^n}{n!}\right)=\sum_{n\geq0}\left( \sum_{i=0}^n \binom{n}{i}f_i \cdot g_{n-i}\right)\frac{x^n}{n!}.
            \end{equation*}
            Note that the coefficient on the left hand side can be rewritten in the following way:
            $$\sum_{i=0}^n \binom{n}{i}f_i \cdot g_{n-i}=\sum_{i+k=n}\frac{n!}{k!\cdot i!}f_i\cdot g_k.$$
        \end{definition}

        \begin{definition}
            The formal derivative of an exponential power series $F(x)=\sum\frac{x^k}{k!}N_k$ is defined as the following power series:
            $$\frac{d}{dx}F:=\sum_{k=1}^\infty k\frac{x^{k-1}}{k!}N_k.$$
        \end{definition}

        The main takeaway of this section is that recursive relations of the numbers in the sequence correspond to differential equations of the generating functions, as the next proposition demonstrates.
        \begin{proposition}\label{genfctder}%5.1.2
            Denote by $F_x:=\frac{d}{dx}F$ the formal derivative of the generating function $F(x)=\sum\frac{x^k}{k!}N_k$. Then $F_x$ is the generating function of the sequence $(N_{k+1})_{k=0}^\infty$, i.e.
            $F_x=\sum_{k=0}^\infty\frac{x^k}{k!}N_{k+1}.$
        \end{proposition}
        The proof consists of slightly rewriting the summands and shifting the index.
        
        \begin{example}[The Drosophilia of Combinatorics]%fibonacci
            This can of course be applied to the Fibonacci numbers. These are given by the sequence $(N_k)_{k=0}^\infty$ with initial condition $N_0=N_1=1$ and the recursive relation $N_{k+2}=N_{k+1}+N_k$ for $k\ge0$. 
            This relation yields an equality of sequences, namely $(N_{k+2})_{k=0}^\infty=(N_{k+1}+N_k)_{k=0}^\infty$.
            
            Let $F=\sum\frac{x^k}{k!}N_k$ be the generating function of the Fibonacci numbers. By the above proposition, the equation of sequences is equivalent to the differential equation $F_{xx}=F_x+F$ for $F$. The initial condition for the sequence translates to an initial condition for the differential equation: $F(0)=F_x(0)=1$.
        \end{example}
        
        \section{Homology and Cohomology}
        %define cup product as PD of intersection product (not always possible, see red box in 4.1.2 book: not poss. for mbar bc singular.
        %define intersection product in particular
        % \textcolor{red}{may have written some Aupstar where there should be Adownstar. WATCH OUT}
        
        We will work with coefficients in $\Q$. Details can be found in \cite{Fulton}.

        For a smooth variety $X$ (e.g. $X=\ptwo$ or $X=\ponex$) of dimension $n$, the Chow group $A_*(X)=\bigoplus_{k\ge0}A_k(X)$ of cycle classes modulo rational equivalence is in fact a ring. Here $A_k(X)$ denotes the group of rational equivalence classes of algebraic cycles of dimension $k$. The product in the ring $A_*(X)$ is the intersection product. 

        We also work with the intersection group $A^*(X)=\bigoplus_{k\ge0}A^k(X)$. The Poincar\'e duality isomorphism $A^k(X)\rightarrow A_{n-d}(X)$, sending a class $c$ to its cap product with the fundamental class $c\cap[X]$, then induces a ring structure on $A^*(X)$, and we call the Poincar\'e dual of the intersection product the cup product. It is denoted by $\cup$. 

        For the spaces under consideration here, there are isomorphisms $A^*(X)\cong H^*(X)$ and $A_*(X)\cong H_*(X)$, both with a doubling of degrees (i.e. $A_k(X)\cong H_{2k}(X)$ and $A^k(X)\cong H^{2k}(X)$). We shall therefore also refer to $A^*(X)$ and $A_*(X)$ as the cohomology and homology rings of $X$, respectively.

        The coarse moduli spaces $\overline M$ are singular varieties. As such, the above notion of an intersection product of cycle classes does not exist for them. It is however still possible to define a cohomology ring with a cup product, see for example 4.1.2 in \cite{invitation}. The classes in $A^*(\overline M)$ that we encounter will be those obtained via pullbacks from $X$.

        The following is an important property of the cap product.
        %projection formula for chow rings https://mathoverflow.net/questions/67228/where-do-all-these-projection-formulas-come-from
        \begin{proposition}[Naturality of $\cap$/Projection Formula]
            Given classes $\beta\in A_*(X)$ and $c\in A^*(Y)$ and a proper morphism $f:X\rightarrow Y$, we have the following equality:
            \begin{equation}\label{projform}
                f_*(f^*(c)\cap \beta) = c\cap f_*(\beta),
            \end{equation}
            see section 3.1 in \cite{fulton2}.
            If $f:X\rightarrow Y$ is a generically finite map of degree $d$, then we have $f_*([X])=d[Y]$. %\textcolor{red}{find source}: https://en.wikipedia.org/wiki/Algebraic_cycle#:~:text=In%20mathematics%2C%20an%20algebraic%20cycle,directly%20accessible%20by%20algebraic%20methods.
            %property like for cap product for closed oriented manifolds
            %or is it https://math.stackexchange.com/questions/3780366/fs-to-x-is-generically-finite-of-degree-d-then-f-fd-dd  ??
        \end{proposition}
        Here $f_*:A_*(X)\rightarrow A_*(Y)$ and $f^*:A^*(Y)\rightarrow A^*(X)$ respectively denote the pushforward (or direct image) and pullback maps induced by $f$. They both preserve degrees. Later we also encounter a pushforward map in cohomology. For $f$ this would be defined as $\text{PD}^{-1}\circ f_*\circ \text{PD}$, where $\text{PD}$ is the Poincar\'e duality isomorphism.

        Let us take a look at our spaces of interest. The ring $A^*(\pr)$ is generated by the classes $\{h^0,\dots,h^r\}$, where $h^i$ corresponds to a generic subvariety of codimension $i$. In particular, $h^0$ is the fundamental class and $h^r$ is the point class.
        Similarly, the ring $A^*(\ponex)$ is generated by classes $\{T_0, T_1, T_2, T_3\}$, where $T_0$ is the fundamental class, $T_1=V$ is the class of a vertical rule, $T_2=H$ is the class of a horizontal rule, and $T_3$ is the point class. 

        We introduce some further notation.

        %see p34 of invit notes and problem notes
        %define integral notation as in FP notes
            % Hodge integrals are defined to be the top intersection products of the ψi and
            % λj classes in Mg,n, see "on hodge integrals" doc
            %panda in hodge integr. paper: integral in GW def is a TAUTOLOGICAL INTEGRAL
            %FP-notes 0.2: "If X is complete, and c is a class in the ring A∗X = L in AkX, we denote by R AdX, and β is a class β c the degree of the class of the zero cycle obtained by
                % evaluating ck on β, where ck is the component of c in AkX. When V is a closed, pure dimensional subvariety of X, we write
                % R V c instead of R [V] c. We use cup
                % product notation ∪ for the product in A∗X."
                
        \begin{definition}
            Given classes $\beta\in A_k(X)$ and $c=\sum c_d\in \bigoplus A^d(X)=A^*(X)$, we denote by $\int_\beta c$ the degree of the zero cycle $c_k\cap\beta\in A_0(X)$. Here the degree refers to the sum of the multiplicities of the points of the subvariety corresponding to the zero cycle.
            For a closed, pure dimensional subvariety $V\subset X$, we may write $\int_V c$ instead of $\int_{[V]} c$.

            Given a class $\sigma\in A_0(X)$, we denote its degree by $\int\sigma$. In particular, if $c\in A^{\dim V}$, we can write $\int_V c=\int c\cap[V]$ .
        \end{definition}

        %\textcolor{red}{issue: greek letter gamma later used for COhomology}
        %connection between intersection and integration (see eg pf of reconstruction)

        %equivalent divisors represent the same class (by def...)
        
        %MUST explain p112 from book
        %DEF PUSHFWD AND PULLBACK. CF P38 OF INVIT. NOTES

    \section{Quantum Formalism}\label{th2quantum}%THE? or: QUANTUM OF FORMALISM
        %for \pone, \pr and \ponex: (IMPORTANT: different for $r\geq 2$ $\pone$ -> after relevant lemma, explain why it does not hold for p1) 
        % degree 0 invariants are =1
        %I_1(hr*hr)=1
        
        In this section, we finally take our first steps in Gromov Witten theory and introduce Gromov-Witten invariants, the quantum product and the quantum cohomology ring. As mentioned in the Introduction, these objects originated in quantum field theory, see for example \cite{Witten:1990hr}. There, GW-invariants correspond to $n$-point correlation functions.
        As before, we begin with the case of $\pr$ and then move on to $\ponex$.
        
        \subsection{Gromov Witten Invariants and Quantum Cohomology for $\pr$}\label{th2qpr}%mention that this section is strongly inspired by \cite{invitiation}
            %We examine the properties of the invariants $I_d(\gamma_1\cdots\gamma_n)$ in $\pr$.

            % INTRODUCE SMALL AND LARGE GAMMA AND EXPLAIN PD. INTRODUCE UNDERLINE GAMMA NOTATION. REMOVE REDUNDANCIES FROM PROP GWMEANINGPR.
            We first revisit the notation from Lemma \ref{lemma1}. Let $\Gamma_1,\dots,\Gamma_n\subset\pr$ be general subvarieties with $[\Gamma_i]=\gamma_i\cap[\pr]$ for all $i$, i.e. such that $\gamma_i\in A^*(X)$ corresponds to $\Gamma_i\in A_*(X)$ via Poincar\'e duality for $X=\pr$. The class corresponding to $[\underline\Gamma]\in A_*(X^n)$ is then given by $\underline\gamma:=\bigcup\tau_i^*(\gamma_i)\in A^*(X^n)$, see (\ref{defgammabar}). 

            Intuitively, the upcoming theory will replicate the arguments from Chapter \ref{classkonts}, but this time on the other side of Poincar\'e duality, allowing us to work with cohomological algebra techniques. In this vein, consider the following product of cohomology classes:
            $$\underline\nu^*(\underline\gamma)=\underline\nu^*(\bigcup\tau_i^*(\gamma_i))=\bigcup\nu_i^*(\gamma_i)\in A^*(\overline M).$$
            The following lemma shows that our new notation is sensible.
            
            \begin{lemma}\label{interslemma}%4.1.3
                The number of maps in the intersection $\underline\nu^{-1}(\underline\Gamma)$ is given by:
                \begin{equation}\label{enumlemma}
                    \int[\underline\nu^{-1}(\underline\Gamma)]=\int\underline\nu^*(\underline\gamma)\cap[\overline M].
                \end{equation}
            \end{lemma}
            We omit the proof. See Lemma 4.1.3 in \cite{invitation} for further details.

            Inspired by this lemma, we define a new kind of invariant.
            \begin{definition}
                We call Gromov-Witten invariant (or GW-invariant) of degree $d$ associated to the cohomology classes $\gamma_1,\dots,\gamma_n\in A^*(\pr)$ the following number:
                \begin{equation*}
                    I_d(\gamma_1\cdots\gamma_n):=\int_{\overline M}\underline\nu^*(\underline\gamma).
                \end{equation*}
            \end{definition}
            We write $\gamma_1\cdots\gamma_n$ instead of $\gamma_1,\dots,\gamma_n$ since $I_d(\gamma_1\cdots\gamma_n)$ is invariant under permutations of the $\gamma_i$ by commutativity of the cup product. The GW-invariants are $\Q$-linear in each entry because both integration and pullback respect sums.
            This allows us to consider inputs $h^i$ from the basis $\{h^0,\dots,h^r\}$ of $A^*(\pr)$ instead of general $\gamma_i$ without any loss of generality.
            The following remark is another important observation.
            \begin{remark}[Dimension Constraint for $\pr$]\label{dimconsprrmk}
                %\textcolor{red}{reference in future chapters}
                In order for $I_d(\gamma_1\cdots\gamma_n)$ to be non-zero, the following condition must be met:
                \begin{equation}\label{dimconspr}
                    \sum_i\text{codim}(\gamma_i)=\text{dim}(\overline M)=rd+r+d+n-3.
                \end{equation}
                This is because we assumed the classes $\gamma_i$ to be homogeneous. Indeed, we have $\underline\nu^*(\underline\gamma)=\bigcup\nu_i^*(\gamma_i)\in A^{\sum\text{codim}(\gamma_i)}(\overline M)$ by flatness of the evaluation maps $\nu_i$ (and $[\overline M]\in A_{\text{dim}(\overline M)}(\overline M)$).
            \end{remark}
            
            %\textcolor{red}{this and other statements could be phrased for general X and applied to ponex case directly}
            The following proposition provides an enumerative interpretation of these invariants.
            \begin{proposition}\label{gwmeaningpr} %4.1.5, need lemma 4.1.3
                Let $r\ge2$ and let $\gamma_1,\dots,\gamma_n\in A^*(\pr)$ be homogeneous classes with $\text{codim}(\gamma_i)\geq2$ for all $i$ and $\sum_{i=1}^n\text{codim}(\gamma_i)=\text{dim}(\mbar)$. 
                Then $I_d(\gamma_1\cdots\gamma_n)$ is the number of rational degree $d$ curves incident to all the subvarieties $\Gamma_1$, \dots, $\Gamma_n$.
            \end{proposition}
 
            \begin{proof} %need lemma 3.5.5
                By Lemma \ref{interslemma}, $I_d(\gamma_1\cdots\gamma_n)$ is equal to the number of n-pointed degree $d$ stable maps $\mu:\pone\rightarrow\pr$ such that $\mu(p_i)\in\Gamma_i$ for $1\leq i\leq n$. As in the discussion leading up to Lemma \ref{lemma2}, we would like to show that counting these maps is equivalent to counting rational curves. Again, we show that there is no freedom in the choice for the marks.
                Consider some $\Gamma_i$. By Lemma \ref{355} (here we need $\text{codim}(\gamma_i)\geq2$), each image curve $\mu(\pone)$ intersects $\Gamma_i$ exactly once, namely in $\mu(p_i)$. Lemma \ref{lemma2} further guarantees that the preimage of this intersection point also consists of a single point, which therefore is the only choice for the mark $p_i$. 
                Thus also in this case, counting maps with marks is the same as counting curves without marks.
                
            \end{proof}

            \begin{corollary}\label{gwndpr}
                In $\ptwo$, we have 
                $I_d(h^2\cdots h^2)=N_d$ for $3d-1$ factors $h^2$.
            \end{corollary}
            
            \begin{lemma}\label{d=0inv}%4.2.1 mapping to a pt
                Unless $n=3$ and $\sum\text{codim}(\gamma_i)=r$, we have $I_0(\gamma_1\cdots\gamma_n)=0$.
                If these conditions are met, then the invariant can be expressed as follows:
                $$I_0(\gamma_1\cdot\gamma_2\cdot\gamma_3)=\int (\gamma_1\cup\gamma_2\cup\gamma_3)\cap[\pr]=1.$$
                %\textcolor{red}{prev prop implies that I0 is equal to the number of points in the intersection of gamm1,gamma2,gamma3.. by gen. bezout (cor 2.4 eh2), number of pts in inters=intersection number is product of degrees of the Gammai. wlog the gammai are generators hi. maybe these have degree 1?}%check if the cases where we use =1 satisfy deg gammai is 1 for all i. not really... se pf of recursion algo. 
            \end{lemma}

            \begin{proof}%see also FP p35 (I)
                Note first that by definition, a constant map is unstable if it has less than three marks. This eliminates the option $n<3$.
                The main idea now is to exploit the isomorphism  $\overline{M}_{0,n}(\pr,0)\cong\overline{M}_{0,n}\times\pr$ from Proposition \ref{deg0iso}. Under this isomorphism, the evaluation maps coincide with the projection $\pi:\overline{M}_{0,n}\times\pr\rightarrow\pr$. We compute:
                \begin{align*}
                    I_0(\gamma_1\cdots\gamma_n) &= \int_{[\overline M]}\nu_1^*(\gamma_1)\cup\cdots\cup\nu_n^*(\gamma_n)\\
                    &= \int_{[\overline{M}_{0,n}\times\pr]}\pi^*(\gamma_1\cup\cdots\cup\gamma_n)\\
                    &= \int_{\pi_*[\overline{M}_{0,n}\times\pr]}\gamma_1\cup\cdots\cup\gamma_n,
                \end{align*}
                where we used the projection formula (\ref{projform}) in the final step. %\textcolor{red}{why no} $p_*$ \textcolor{red}{? argument would work for n=3 because then the projection is 1 to 1, but we have not deduced that yet. proper pushforward???}%https://stacks.math.columbia.edu/tag/02R3
                Note that in the second to last line, we integrate the following class:
                $$\pi^*(\gamma_1\cup\cdots\cup\gamma_n)\cap[\overline{M}_{0,n}\times\pr]\in A_0(\overline{M}_{0,n}\times\pr),$$
                and not the class one would expect from the formulation of the projection formula, that is:
                $$\pi_*\bigl(\pi^*(\gamma_1\cup\cdots\cup\gamma_n)\cap[\overline{M}_{0,n}\times\pr]\bigr).$$
                But this is justified: proper pushforward respects rational equivalence, and rationally equivalent zero cycles have the same degree, which is exactly what the integral detects.
                
                Now unless $n=3$, we have $\dim(\overline{M}_{0,n})>0$. But then $\pi_*[\overline{M}_{0,n}\times\pr]=0$, because the fibers of $\pi$ have positive codimension. In particular there is a positive integer $k$ such that $\pi_*[\overline{M}_{0,n}\times\pr]\in A_{r+k}(\pr)=0$. %since the induced homom. in homology goes from A_k -> A_k and Pr only has non zero A_k up to k=r, we would land in a trivial A_k of too high dimension.
                Thus the only remaining case is $n=3$, where $\overline{M}_{0,n}\times\pr\cong\pr$ and thus $I_0(\gamma_1\cdot\gamma_2\cdot\gamma_3)=\int_{\pr}\gamma_1\cup\gamma_2\cup\gamma_3$.
                Recall that w.l.o.g., we may assume that the input classes are generators, i.e. that we are dealing with $I_0(h^i\cdot h^j\cdot h^k)$ with $i+j+k=r$. Then the product $h^i\cup h^j\cup h^k=h^r$ is the point class, and thus the invariant is equal to 1.

            \end{proof}

            The following is a generalization of Fact \ref{n1}.
            \begin{lemma}\label{2ptin}%4.2.2, uses r \geq 2
                For $r\geq2$ and $n<3$, the only non-zero GW-invariant is $$I_1(h^r\cdot h^r)=1.$$ 
            \end{lemma}
            The case $r=1$ will be discussed in detail in Proposition \ref{p1invs}.
            \begin{proof}
               By what we've seen in the proof of Lemma \ref{d=0inv}, we may assume that $d>0$. With the dimension constraint (\ref{dimconspr}), we have:
               $$nr\geq\sum\codim\gamma_i)=\dim\mbar=rd+r+d+n-3\geq2r+n-2.$$
               Since we assumed that $r\geq2$, this is impossible unless $n=2$.
               In the latter case, the above inequalities become equalities and thus $\gamma_1=\gamma_2=h^r$. The only way for the dimension constraint to be satisfied in this setting is if $d=1$. By Proposition \ref{gwmeaningpr} (and the Introduction), this invariant is equal to 1. %\textcolor{red}{prove that 2 pts uniquely determine line, i.e. kill all coeffs. should follow from lemma gwmeaningpr with argument as in classical Kontsevich} actually, that is an axiom of projective geometry
                
            \end{proof}

            Consider the following commutative diagram: %\textcolor{red}{replace by mbarn, mbarn1 then reuse in ponex?}
            \begin{equation}\label{diag}
                \begin{tikzcd}
                    \overline{M}_{0,n+1}(\pr,d) \arrow{rr}{\tilde\nu_i} \arrow[swap]{dd}{\epsilon} && \pr \\ \\
                    \mbar \arrow[swap]{rruu}{\nu_i} &&
                \end{tikzcd}
            \end{equation}
            where $\epsilon$ is the forgetful map and $\tilde\nu_i$ and $\nu_i$ are the evaluation maps for the first $n$ marks of $\overline{M}_{0,n+1}(\pr,d)$ (respectively $\mbar$). Note that the map $\epsilon$ is not defined if $n<3$ and $d=0$, since there are no degree $0$ stable maps with less than three marks. The commutativity implies that $\tilde\nu_i^*(\gamma_i)=\epsilon^*\nu_i^*(\gamma_i)$ in $A^*(\overline{M}_{0,n+1}(\pr,d))$.

            \begin{lemma}\label{fundclassinv} %4.2.3 why is this called string equation
                Consider the fundamental class $1=h^0\in A^0(\pr)$. If $\gamma_i=1$ for some $1\leq i\leq n$, then $I_d(\gamma_1\cdots\gamma_n)=0$ unless $d=0$ and $n=3$. 
            \end{lemma}

            \begin{proof}
                We work in $\overline{M}_{0,n+1}(\pr,d)$. The case $n<2$ was handled in Lemma \ref{d=0inv}. Assume now that $n>2$ and $d>0$. Let w.l.o.g. the fundamental class be the final input, i.e. $\gamma_{n+1}=1$. Note that $\nu_{n+1}^*(1)=1\in\overline{M}_{0,n+1}(\pr,d)$. Using the projection formula (\ref{projform}) for the map $\epsilon$ in diagram (\ref{diag}), we compute:
                \begin{align*}
                    I_d(\gamma_1\cdots\gamma_n\cdot1) 
                    & = \int_{[\overline{M}_{0,n+1}(\pr,d)]}\underline{\tilde\nu}^*(\underline\gamma)\cup\tilde\nu_{n+1}^*(1)\\
                    & =\int\underline{\tilde\nu}^*(\underline\gamma)\cap[\overline{M}_{0,n+1}(\pr,d)]\\
                    & =\int\underline\nu^*(\underline\gamma)\cap\epsilon_*[\overline{M}_{0,n+1}(\pr,d)].%by projform. but where is \epsilon_*?
                \end{align*}
                But, by an argument analogous to that in the proof of Lemma \ref{d=0inv}, the class $\epsilon_*[\overline{M}_{0,n+1}(\pr,d)]$ is zero. This concludes the proof. 
                
            \end{proof}

            The following lemma allows us to handle codimension $1$ classes and is therefore also referred to as 'divisor equation'.
            \begin{lemma}\label{divisoreq} %4.2.4
                Consider the hyperplane class $h=h^1\in A^1(\pr)$. If $\gamma_i=h$ for some $1\leq i\leq n+1$, say $\gamma_{n+1}=h$, then: 
                $$I_d(\gamma_1\cdots\gamma_n\cdot h)=d\cdot I_d(\gamma_1\cdots\gamma_n).$$ 
            \end{lemma}

            \begin{proof}
                Denote by $H\subset\pr$ the hyperplane satisfying $\tilde\nu_{n+1}^*(h)\cap[\overline{M}_{0,n+1}(\pr,d)]=[\tilde\nu_{n+1}^{-1}(H)]$. The locus $\tilde\nu_{n+1}^{-1}(H)$ consists of those maps in $\overline{M}_{0,n+1}(\pr,d)$ whose mark $p_{n+1}$ gets mapped to $H$. Restricting the forgetful map $\epsilon$ from diagram (\ref{diag}) to $\tilde\nu_{n+1}^{-1}(H)$ yields a generically finite map $\epsilon|$ of degree $d$. Indeed, the image curve of a general map $\mu\in\overline{M}_{0,n}(\pr,d)$ intersects $H$ in $d$ points, say $Q_1,\dots,Q_d$. The mark $p_{n+1}$ can be placed upon any of the $d$ points $\mu^{-1}(Q_i)$, yielding a total of $d$ maps in $\epsilon|^{-1}(\mu)$. %if a p_i already gets sent to H, rename. Is this what generality means?
                
                As when proving the previous lemmas, we shall compute the relevant invariant using the projection formula (\ref{projform}):
                \begin{align*}
                    I_d(\gamma_1\cdots\gamma_n\cdot h) 
                    & = \int\underline{\tilde\nu}^*(\underline\gamma)\cup\tilde\nu_{n+1}^*(h)\cap{[\overline{M}_{0,n+1}(\pr,d)]}\\
                    & =\int\underline{\tilde\nu}^*(\underline\gamma)\cap[\tilde\nu_{n+1}^{-1}(H)]\\
                    & =\int\underline\nu^*(\underline\gamma)\cap\epsilon_*[\tilde\nu_{n+1}^{-1}(H)]\\
                    & =\int\underline\nu^*(\underline\gamma)\cap d[\mbar]\\
                    & =d\cdot I_d(\gamma_1\cdots\gamma_n),%by projform. but where is \epsilon_*?
                \end{align*}
                where we used the compatibility between the cup and cap products in the first step, see for example Proposition 4.17 in \cite{greg}.
                
            \end{proof}

            Recall the notation $D_{AB}=D(A,B;d_A,d_B)$ from section \ref{kfp2}. As before, denote the gluing mark between the two twigs by $x$ and consider the diagonal $\Delta\subset\pr\times\pr$ as in Remark \ref{recstruc}. Then we can once more consider the inclusion $\iota:D_{AB}\hookrightarrow \overline{M}_A\times\overline{M}_B$ by writing $D_{AB}=(\nu_{x_A}\times\nu_{x_B})^{-1}(\Delta)$, where $\overline{M}_A=\overline{M}_{0,A\cup\{x\}}(\pr,d_A)$ and $\overline{M}_B=\overline{M}_{0,B\cup\{x\}}(\pr,d_B)$.

            \begin{lemma}[Splitting Lemma for $\pr$]\label{splitting}
                Let $\alpha:D_{AB}\hookrightarrow\overline M$ be the inclusion. Then for all $\gamma_1,\dots,\gamma_n\in A^*(\pr)$, the following holds in $A^*(\overline M_A\times \overline M_B)$:
                $$\iota_*\alpha^*\underline\nu^*(\underline\gamma)=\sum_{e+f=r}\bigl(\bigcup_{a\in A}\nu_a^*(\gamma_a)\cup\nu_{x_A}^*(h^e)\bigr)\cdot\bigl(\bigcup_{b\in B}\nu_b^*(\gamma_b)\cup\nu_{x_B}^*(h^f)\bigr).$$
                In particular, we get the following formula:
                $$\int_{D_{AB}} \nu_1^*(\gamma_1)\cup\dots\cup\nu_n^*(\gamma_n) = \sum_{e+f=r} I_{d_A}(\prod_{a\in A}\gamma_a\cdot h^e)\cdot I_{d_B}(\prod_{b\in B}\gamma_b\cdot h^f).$$
            \end{lemma}
            We omit the proof. The main idea is to exploit the compatibility between the evaluation maps $\nu$ and the recursive structure of the boundary divisor $D_{AB}$ (Proposition \ref{compeval}). 
            The second formula is then obtained from the first by integration. See Lemma 4.3.2 and Corollary 4.3.3 in \cite{invitation}.

            The following theorem generalizes Theorem \ref{nd} to $\pr$. Again, the only assumption is that there exists a unique line through two general points. This time, we do not give an explicit recursion, but only prove the computability of the invariants. In cases where we need explicit expressions for invariants, we shall make use of the relevant parts of the arguments of this proof, see for example Proposition \ref{3ptpr}. 
            \begin{theorem}[Reconstruction for $\pr$ with $r\geq2$]\label{reconpr}
                Given that $I_1(h^r\cdot h^r)=1$, all genus $0$ GW-invariants for $\pr$ can be computed recursively.
            \end{theorem}

            \begin{proof}
                We will prove this theorem via a threefold recursion. Let $I_d(\gamma_1\cdots\gamma_n)$ be the invariant to be computed. Our goal is to reduce this computation to the knowledge of  $I_1(h^r\cdot h^r)=1$. Concretely, this means reducing it to the computation of invariants of lower degree or with fewer marks.

                Thanks to Lemmas \ref{fundclassinv} and \ref{divisoreq}, we may assume that $\codim \gamma_i)\geq2$ for $1\leq i\leq n$. We may also assume w.l.o.g. that $\codim\gamma_1)\geq\dots\geq\codim\gamma_n)$.

                We must make some preparations for the recursions. Due to the first reduction above, we may rewrite the final class as $\gamma_n=\lambda_1\cup\lambda_2$ with $\codim\lambda_1)\leq\codim\lambda_2)<\codim\gamma_n).$
                Similarly to the proof of Theorem \ref{nd}, consider the space $\overline{M}_{0,n+1}(\pr,d)$ with marks $m_1,m_2,p_1,\dots,p_{n-1}$. In that proof, we considered the subvariety of $\overline M$ of maps sending $m_i$ to $L_i$ and $p_i$ to $Q_i$, where $L_1,L_2$ where lines and $Q_1,\dots,Q_{n-1}$ where points. We then compared its intersections with the equivalent boundary divisors $D(m_1,m_2|p_1,p_2)\equiv D(m_1,p_1|m_2,p_2)$. 
                Note how this setup mirrors the constellation of our classes $\lambda_1,\lambda_2$ with smaller codimensions than $\gamma_1,\dots,\gamma_{n-1}$. Our approach in this setting is essentially the Poincar\'e dual of that in the proof of Theorem \ref{nd}. 
                In particular, we now take the following class:
                $$\nu_{m_1}^*(\lambda_1)\cup\nu_{m_2}^*(\lambda_2)\cup\nu_{p_1}^*(\gamma_1)\cup\dots\cup\nu_{p_{n-1}}^*(\gamma_{n-1}),$$
                and we integrate it over the two equivalent boundary divisors. Applying the Splitting Lemma \ref{splitting} to each component leads to the following equation:
                \begin{equation}\label{reconpreq} %reformat?
                    \sum_{\substack{A\cup B:\\m_1,m_2\in A \\ p_1,p_2\in B }} \sum_{e+f=r} I_{d_A}(\prod_{a\in A}\gamma_a\cdot h^e)\cdot I_{d_B}(\prod_{b\in B}\gamma_b\cdot h^f) =\sum_{\substack{A\cup B:\\m_1,p_1\in A \\ m_2,p_2\in B }} \sum_{e+f=r} I_{d_A}(\prod_{a\in A}\gamma_a\cdot h^e)\cdot I_{d_B}(\prod_{b\in B}\gamma_b\cdot h^f),
                \end{equation}
                where the external sums are over partitions $A\cup B=\{m_1,m_2,p_1,\dots,p_{n-1}\}$ with $d_A+d_B=d$. The invariants $I_{d_A}$ and $I_{d_B}$ are calculated in the spaces of curves with marks $A\cup \{x\}$ and $B\cup \{x\}$ respectively, and the classes $h^e$ and $h^f$ correspond to the gluing mark $x$. %\textcolor{red}{explain more?}
                
                We are now ready to describe the recursion algorithm.
                If $d_A>0$ and $d_B>0$, then we also have $d_B<d$ and $d_A<d$. We can pass these invariants to a first recursion in $d$. It remains to show that this recursion terminates, i.e. that the other cases can be reduced to $I_1(h^r\cdot h^r)$. %This is the first recursion. 

                If $d_A=0$ or $d_B=0$, we know that $I_0$ can only bear three marks and the sum of their codimensions must be equal to $r$ by Lemma \ref{d=0inv}. Writing $c_i:=\codim\lambda_i)$ and $b_i:=\codim\gamma_i)$, the relevant terms for this case are:
                \begin{align*}
                    &I_0(\lambda_1\cdot\lambda_2\cdot h^{r-c_1-c_2})\cdot I_d(\gamma_1\cdot\gamma_2\cdots\gamma_{n-1}\cdot h^{c_1+c_2}),\\
                    &I_d(\lambda_1\cdot\lambda_2\cdot\gamma_3\cdots\gamma_{n-1}\cdot h^{b_1+b_2})\cdot I_0(\gamma_1\cdot\gamma_2\cdot h^{r-b_1-b_2}),
                \end{align*}
                on the left hand side and
                \begin{align*}
                    &I_0(\lambda_1\cdot\gamma_1\cdot h^{r-c_1-b_1})\cdot I_d(\lambda_2\cdot\gamma_2\cdots\gamma_{n-1}\cdot h^{b_1+c_1}),\\
                    &I_d(\lambda_1\cdot\gamma_1\cdot\gamma_3\cdots\gamma_{n-1}\cdot h^{c_2+b_2})\cdot I_0(\lambda_2\cdot\gamma_2\cdot h^{r-c_2-b_2}),
                \end{align*}
                on the right hand side. The $I_0$ invariants are equal to $1$ by Lemma \ref{d=0inv}. The first term is therefore equal to the desired invariant $I_d(\gamma_1\cdots\gamma_n)$. The class with the lowest codimension in the remaining three $I_d$ is in each case either $\lambda_1$ or $\lambda_2$, both of which have strictly smaller codimension than $\gamma_n$. We can thus pass these to a second recursion in $\codim\tilde{\gamma}_n)$, where $\Tilde{\gamma}_n$ stands for the final input class of $I_d$ after reordering the classes up to decreasing codimension.
                Once this recursion arrives at $\codim\tilde{\gamma}_n)=1$, we can remove the final class from $I_d$ using Lemma \ref{d=0inv} and pass $I_d$ to a third and final recursion in $n$. 
                Once this recursion reaches $n=2$, it recovers either $I_1(h^r\cdot h^r)=1$ or $0$, the only invariants with less than three points according to Lemma \ref{2ptin}. 
                Thus the third, second and first recursions terminate.% \textcolor{red}{explain in more detail, cf algoko.}
                
            \end{proof}
            
            \begin{remark}
                Since we used Lemma \ref{2ptin} at the end of the proof, we had to assume $r\geq2$. The $r=1$ case is much simpler, in fact we will determine all GW-invariants of $\pone$ explicitly in Proposition \ref{p1invs}.
            \end{remark}

            We would now like to explore the GW-invariants $I_d(\gamma_1\cdots\gamma_n)$ further by encoding them into the language of (exponential) generating functions. To this end, we first get rid of the parameter $d$.

            \begin{propdef}[Collected GW-Invariants]\label{collinvpr}
                We define the collected GW-invariant associated to the cohomology classes $\gamma_1,\dots,\gamma_n\in A^*(\pr)$ as follows: 
                $$I(\gamma_1\cdots\gamma_n):=\sum_{d=0}^\infty I_d(\gamma_1\cdots\gamma_n).$$
                It is a finite sum. More precisely, at most one summand is non-zero and there is thus no loss of information in considering the collected invariant $I(\gamma_1\cdots\gamma_n)$.
            \end{propdef}

            \begin{proof}
                Let's examine the summands on the right hand side using the dimension constraint (\ref{dimcons}). The latter implies that $I_d(\gamma_1\cdots\gamma_n)=0$ unless $\sum \text{codim}(\gamma_i)=\dim(\mbar)=rd+r+d+n-3$. Thus only the summand satisfying $$d=\frac{\sum \text{codim}(\gamma_i)-r-n+3}{r+1}$$ survives.
                
            \end{proof}

            As previously mentioned, the linearity of the invariants allows us to assume w.l.o.g. that the input classes are given by basis elements $h^i\in A^*(\pr)$. The collected invariants are therefore determined by the number of occurences of each class $h^i$, so a general invariant looks like
            $I((h^0)^{\bullet a_0}\cdots (h^r)^{\bullet a_r}).$ We now define the generating function for these numbers.
            %potentials
            
            \begin{definition}[Potentials]
                Using formal variables $x_0,\dots,x_r$, we define the GW-potential as follows:
                $$\Phi(x_0,\dots,x_r):=\sum_{a_0,\dots,a_r}\frac{x_0^{a_0}\cdots x_r^{a_r}}{a_0!\cdots a_r!}I((h^0)^{\bullet a_0}\cdots (h^r)^{\bullet a_r}).$$
                Using multi-index notation with $\mathbf{x}:=(x_0,\dots,x_r)$ and $\mathbf{h}^{\mathbf{a}}:=(h^0)^{\bullet a_0}\cdots (h^r)^{\bullet a_r}$, we may more compactly rewrite this as $\Phi(\mathbf{x})=\sum_{\mathbf{a}}\frac{\mathbf{x}^{\mathbf{a}}}{\mathbf{a}!}I(\mathbf{h}^{\mathbf{a}})$, where $\mathbf{a}\in(\Z_{\geq0})^{r+1}$.
            \end{definition}

            In what follows, it will be useful to decompose $\Phi$ into a part of degree $0$ and a part of positive degree. We first define $I_+:=\sum_{d>0}I_d$ and can thus write $I=I_0+I_+$. We then define the classical potential $\Phi^{\text{cl}}$ with the same formula as $\Phi$, but with only $I_0$ instead of $I$.
            %$$\Phi^{\text{cl}}:=(x_0,\dots,x_r):=\sum_{a_0,\dots,a_r}\frac{x_0^{a_0}\cdots x_r^{a_r}}{a_0!\cdots a_r!}I_0((h^0)^{\bullet a_0}\cdots (h^r)^{\bullet a_r}).$$
            By Lemma \ref{d=0inv}, only invariants $I_0$ with three input classes whose codimensions sum up to $r$ survive, which yields the following simplified formula for the classical potential:
            $$\Phi^{\text{cl}}(\mathbf{x})=\sum_{i+j+k=r}\frac{x_ix_jx_k}{3!}I_0(h^i\cdot h^j\cdot h^k),$$
            where the factor $3!$ compensates for the different possible permutations of $i,j,k$ which all yield the same summand.
            
            We call the remaining part of $\Phi$ the quantum potential and denote it by $\Gamma:=\Phi-\Phi^{\text{cl}}$. Clearly, the explicit expression of $\Gamma$ is identical to that of $\Phi$, except that it sums up $I_+$ instead of $I$. The terminology is justified by Proposition \ref{classprodpr} and Definition \ref{qprodpr} below.
            
            Note that by Proposition \ref{genfctder}, we have $\Phi_i:=\frac{\delta}{\delta x_i}\Phi=\sum_{\mathbf{a}}\frac{\mathbf{x}^\mathbf{a}}{\mathbf{a}!}I(\mathbf{h}^\mathbf{a}\cdot h^i).$ %textcolor{red}{expand?}
            In particular, we have the following:
            \begin{equation}\label{phiijk}
                \Phi_{ijk}(\mathbf{x})=\sum_{\mathbf{a}}\frac{\mathbf{x}^\mathbf{a}}{\mathbf{a}!}I(\mathbf{h}^\mathbf{a}\cdot h^i \cdot h^j\cdot h^k).
            \end{equation}

            %5.2.6 "classical product"
            \begin{proposition}\label{classprodpr}
                The structure constants of the cup product in $A^*(\pr)$ are given by the third derivatives of the classical potential. More precisely, for arbitrary classes $h^i,h^j\in A^*(\pr)$, it holds that:
                \begin{equation*}
                    h^i\cup h^j = \sum_{e+f=r}\Phi_{ije}^{\text{cl}}h^f.
                \end{equation*}
                We shall therefore also refer to the cup product as the classical product.
            \end{proposition}

            \begin{proof}
                Consider the basis $\{h^0,\dots,h^r\}$ of $A^*(\pr)$. By definition we have $h^i\cup h^j=h^{i+j}\in A^{i+j}(\pr)$. 
                % We deduce the following: SEEMS TO BE UNNECESSARY!
                % \begin{equation}\label{inegcases}
                %     \int_{\pr} h^i\cup h^j=
                %     \begin{cases}
                %         1, \text{ if } i+j=r;\\
                %         0, \text{ otherwise}.
                %     \end{cases}
                % \end{equation}
                % \textcolor{red}{why exactly do Hi, Hj intersect in exactly one pt?}
                From Lemma \ref{d=0inv}, we know that $I_0(h^i\cdot h^j \cdot h^e)$ %=\int_{\pr} (h^i\cup h^j\cup h^e)$ and that this expression 
                is zero unless $i+j=r-e$. We can therefore rewrite $h^{i+j}$ as $\sum_e I_0(h^i\cdot h^j \cdot h^e)h^{r-e}$. Replacing $r-e$ by $f$ yields the desired expression.
                
            \end{proof}
            
            %\textcolor{red}{introduce coord free form here?}
            \begin{definition}[Quantum Product]\label{qprodpr}
                We define the quantum product of two classes $h^i,h^j\in A^*(\pr)$ as:
                $$h^i*h^j:=\sum_{e+f=r}\Phi_{ije}h^f\in \Q[[\mathbf{x}]] \otimes_{\Z}A^*(\pr),$$
                with $e,f\in\Z_{\ge0}$. We define the quantum product on all of $\Q[[\mathbf{x}]] \otimes_{\Z}A^*(\pr)$ as the $\Q[[\mathbf{x}]]$-linear extension of this definition.
            \end{definition}

            By Proposition \ref{classprodpr}, we can now decompose the quantum product as follows:
            $$h^i*h^j=h^i\cup h^j+\sum_{e+f=r}\Gamma_{ije}h^f.$$
            Intuitively, the quantum product of two classes encodes 'fuzzy' intersections between the corresponding subvarieties. Indeed, the classical product is the Poincar\' dual of the intersection product, and its structure constants come from the classical potential. The quantum product retains all this information, but due to the additional structure constants coming from the quantum potential, it also takes into account the number of rational curves of positive degree incident to both subvarieties, which on the other side of Poincar\'e duality leads to a less strict notion of intersection.

            \begin{remark}[Commutativity of the Quantum Product]
                Since changing the order of the partial derivatives has no effect on $\Phi_{ijk}$, the quantum product is commutative.
            \end{remark}

            \begin{lemma}%5.2.9
                The identity element for the quantum product is given by the fundamental class $h^0$.
            \end{lemma}

            \begin{proof}
                By definition, we have $h^0*h^i=\sum_{e+f=r}\sum_\mathbf{a}\frac{\mathbf{x}^\mathbf{a}}{\mathbf{a}!}I(\mathbf{h}^{\bullet\mathbf{a}}\cdot h^0 \cdot h^i\cdot h^e)h^f$ for all $i$. But by Lemma \ref{fundclassinv}, only the invariant $I_0(h^0\cdot h^i \cdot h^e)$ in the sum $I$ survives. The product reduces to:
                $$h^0*h^i=\sum_{e+f=r}I_0(h^0\cdot h^i \cdot h^e)h^f,$$
                where all summands except the one in $e=r-i$ vanish. By Lemma \ref{d=0inv}, this summand is just $h^i$.
                
            \end{proof}

            The following proposition will be proven in a way similar to Kontsevich's formula for $\ptwo$ (Theorem \ref{nd}). Indeed, it also relies on the fundamental equivalence (\ref{fundrel}) and counts contributions from two equivalent boundary divisors. Instead of intersections, we now work with cup products. Following a by now familiar pattern, this argument could be described as the Poincar\'e dual of that in the proof of Theorem \ref{nd}. Fittingly, we shall deduce that theorem from this proposition in section \ref{konassp2}.

            \begin{proposition}[Associativity of the Quantum Product in $\pr$]\label{asspr}
            %associativity of * (for pr) 
                For arbitrary classes $h^i,h^j,h^k\in A^*(\pr)$, the following equality holds:
                $$(h^i*h^j)*h^k=h^i*(h^j*h^k).$$
            \end{proposition}
            
            \begin{proof}
                We will first examine the associativity condition from the statement more precisely. We will then show that it can be deduced from the fundamental equivalence (\ref{fundrel}).

                By definition of the quantum product, we have:
                \begin{align*}
                    (h^i*h^j)*h^k= &\;(\sum_{e+f=r}\Phi_{ije}h^f)*h^k=\sum_{e+f=r}\sum_{l+m=r}\Phi_{ije}\Phi_{fkl}h^m,  \\
                    h^i*(h^j*h^k)= &\;h^i*(\sum_{e+f=r}\Phi_{jke}h^f)=\sum_{e+f=r}\sum_{l+m=r}\Phi_{jke}\Phi_{ifl}h^m.
                \end{align*}              
                Due to the linear independence of the $h^m$, comparing coefficients yields the following equations for any $i,j,k,l$:
                \begin{equation*}
                    \sum_{e+f=r}\Phi_{ije}\Phi_{fkl}=\sum_{e+f=r}\Phi_{jke}\Phi_{ifl}.
                \end{equation*}
                 These are also called the \textit{WDVV equations} after Witten, Dijkgraaf, Verlinde and Verlinde. As we will see below, they serve as a quantum counterpart to the classical fundamental relation (\ref{fundrel}). In order to be able to translate this differential equation for the generating functions $\Phi$ into a recursion for the invariants $I$ via the product rule (\ref{prodrule}), we define $\gamma:=\sum_{i=0}^rx_ih^i$ and rewrite $\Phi_{ijk}$ in a coordinate-free form: 
                
                 $$\Phi_{ijk}(\mathbf{x})=\sum_\mathbf{a}\frac{\mathbf{x}^\mathbf{a}}{\mathbf{a}!}I(\mathbf{h}^{\bullet\mathbf{a}}\cdot h^i \cdot h^j\cdot h^k)=I(\text{exp}(\gamma)\cdot h^i \cdot h^j\cdot h^k)=\sum_{n\geq0}\frac{1}{n!}I(\gamma^{\bullet n}\cdot h^i \cdot h^j\cdot h^k),$$
                where the right hand side can be thought of as the restriction of a power series in one variable
                $\bigl(\sum_{n\geq0}\frac{y^n}{n!}I(\gamma^{\bullet n}\cdot h^i \cdot h^j\cdot h^k)\bigr)|_{y=1}$. The second equality uses the following:
                $$\text{exp}(\gamma)=\text{exp}(\sum x_ih^i)=\prod\text{exp}(x_ih^i)=\prod\sum_{a_i\geq0}\frac{x_i^{a_i}}{a_i!}(h^i)^{a_i}=\sum_{\mathbf{a}}\frac{\mathbf{x}^\mathbf{a}}{\mathbf{a}!}\mathbf{h}^{\bullet\mathbf{a}},$$ which also leads to $I(\text{exp}(\gamma))=\Phi(\mathbf{x})$. 
                
                This trick allows us to deduce the following recursion: 
                 \begin{align}\label{wdvvpr}
                    \begin{split}
                        &\sum_{e+f=r}\;\sum_{n_A+n_B=n}\binom{n}{n_A}I(\gamma^{\bullet n_A}\cdot h^i\cdot h^j\cdot h^e)I(\gamma^{\bullet n_B}\cdot h^f\cdot h^k\cdot h^l)\\
                        = & \sum_{e+f=r}\;\sum_{n_A+n_B=n}\binom{n}{n_A}I(\gamma^{\bullet n_A}\cdot h^j\cdot h^k\cdot h^e)I(\gamma^{\bullet n_B}\cdot h^i\cdot h^f\cdot h^l).
                    \end{split}
                 \end{align}
                We now show the associativity of $*$ by proving this recursion via the fundamental equivalence (\ref{fundrel}). Since the equivalence involves four named marks, we shall work in the space $\overline{M}_{0,n+4}(\pr,d)$ for arbitrary $n$ and $d$, with four marks called $p_1$, $p_2$, $p_3$, $p_4$ and $n$ marks that shall remain unnamed.
                By the fundamental relation (\ref{fundrel}), we have:
                $$D(p_1p_2|p_3p_4)\equiv D(p_2p_3|p_1p_4).$$
                Note how the marks are permuted in the same way as the indices $i$, $j$, $k$, $l$ in (\ref{wdvvpr}). We exploit this by considering the pullbacks via the evaluation maps of the corresponding classes, namely $\nu_1^{-1}(h^i)$, $\nu_2^{-1}(h^j)$, $\nu_3^{-1}(h^k)$ and $\nu_4^{-1}(h^l)$. We further denote by $\underline\nu^{*}(\underline \gamma)$ the cup product of the pullbacks of $n$ copies of $\gamma$ via the remaining $n$ evaluation maps. Integrating the cup product of these pullbacks over the equivalent boundary divisors yields:
                \begin{align}\label{wdvvintpr}
                    \begin{split}
                        &\int_{D(p_1p_2|p_3p_4)}\underline\nu^{*}(\underline \gamma)\cup\nu_1^{-1}(h^i) \cup\nu_2^{-1}(h^j) \cup \nu_3^{-1}(h^k)\cup\nu_4^{-1}(h^l)\\
                        = & \int_{D(p_2p_3|p_1p_4)}\underline\nu^{*}(\underline \gamma)\cup\nu_1^{-1}(h^i) \cup\nu_2^{-1}(h^j) \cup \nu_3^{-1}(h^k)\cup\nu_4^{-1}(h^l).
                    \end{split}
                 \end{align}
                 Like in the proof of Theorem \ref{nd}, we examine the contributions from each summand (i.e. component) making up the boundary divisors. The situation on the left hand side is depicted below (with \textit{e.m.} standing for \textit{extra marks}).

                \definecolor{rvwvcq}{rgb}{0.08235294117647059,0.396078431372549,0.7529411764705882}
                \definecolor{wrwrwr}{rgb}{0.3803921568627451,0.3803921568627451,0.3803921568627451}
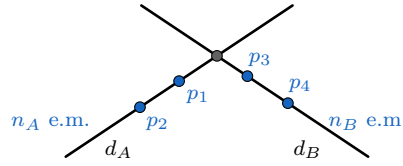
                \begin{tikzpicture}[line cap=round,line join=round,>=triangle 45,x=1cm,y=1cm]
                    \clip(-5,0.5) rectangle (4.5,4); %size of rectangle containing image (given by 2pts)
                    \draw [line width=1pt] (1,3)-- (4,1);
                    \draw [line width=1pt] (3,3)-- (0,1);
                        \begin{scriptsize}
                            \draw[color=black] (3.2,1.1) node {$d_{B}$};
                            \draw[color=rvwvcq] (4,1.44) node {$n_B$ e.m.};
                            
                            \draw [fill=rvwvcq] (2.4025391414130275,2.0649739057246483) circle (2pt);
                            \draw[color=rvwvcq] (2.57,2.2890429805862818) node {$p_{3}$};
                            
                            \draw [fill=rvwvcq] (2.938269319967983,1.7078204533546781) circle (2pt);
                            \draw[color=rvwvcq] (3.102541216255241,1.9433276303855245) node {$p_{4}$};
                            
                            \draw[color=black] (0.7,1.1) node {$d_{A}$};
                            \draw[color=rvwvcq] (-0.2,1.44) node {$n_A$ e.m.};
                            
                            \draw [fill=rvwvcq] (1.4993503144908957,1.9995668763272638) circle (2pt);
                            \draw[color=rvwvcq] (1.75,1.8) node {$p_{1}$};
                            
                            \draw [fill=rvwvcq] (0.9819670598590557,1.654644706572704) circle (2pt);
                            \draw[color=rvwvcq] (1.2162624630513585,1.44) node {$p_{2}$};
                            
                            \draw [fill=wrwrwr] (2,2.3333333333333335) circle (2pt);
                            \draw[color=wrwrwr] (2,2.7) node {};
                        \end{scriptsize}
                \end{tikzpicture}
                
                On both sides of (\ref{wdvvintpr}), there are $\binom{n}{n_A}$ ways of distributing the marks and we must have $d_A+d_B=d$. Luckily, the remainder of the count is simpler than in the proof of Kontsevich's formula: all of the $n$ extra classes are given by $\gamma$, so each component yields the same contribution.
                
                Applying the Splitting Lemma \ref{splitting} to the components of the boundary divisors in (\ref{wdvvintpr}) therefore leads to the following equation:
                \begin{align*}
                    \begin{split}
                        &\sum_{d_A+d_B=d}\;\sum_{n_A+n_B=n}\binom{n}{n_A}\sum_{e+f=r}I_{d_A}(\gamma^{\bullet n_A}\cdot h^i\cdot h^j\cdot h^e)I_{d_B}(\gamma^{\bullet n_B}\cdot h^f\cdot h^k\cdot h^l)\\
                        = &\sum_{d_A+d_B=d}\;\sum_{n_A+n_B=n}\binom{n}{n_A}\sum_{e+f=r}I_{d_A}(\gamma^{\bullet n_A}\cdot h^j\cdot h^k\cdot h^e)I_{d_B}(\gamma^{\bullet n_B}\cdot h^i\cdot h^f\cdot h^l).
                    \end{split}
                 \end{align*}
                Finally, summing over all $d$ eliminates the first sum on both sides and recovers the collected invariants $I$, which leads to the desired equation (\ref{wdvvpr}).
                
            \end{proof}

        \begin{remark}[Big and small Quantum Cohomology]
            The quantum product determines a new ring structure in $\Q[[\mathbf{x}]] \otimes_{\Z}A^*(\pr)$. We call this ring the quantum cohomology ring of $\pr$ and we denote it by $Q^*(\pr)$.

            A simpler yet useful variation of this quantum cohomology is obtained by setting all variables in $\Phi_{ijk}(\mathbf{x})$ to zero except for those corresponding to codimension $1$ classes, i.e. $x_1$. We can rewrite $x:=x_1$ and simplify the $\Phi_{ijk}$ using Lemma \ref{divisoreq}:
            \begin{align*}
                \Phi_{ijk}(x) & = \sum_{n=0}^\infty \frac{x^n}{n!}\sum_{d\ge0} I_d(h^{\bullet n}\cdot h^i\cdot h^j\cdot h^k) \\
                & = \sum_{n=0}^\infty \frac{x^n}{n!}\sum_{d\ge0}d^n\cdot I_d(h^i\cdot h^j\cdot h^k).
            \end{align*}
            We call the product given by these structure constants the small quantum product. By the above calculation, it involves only three-point invariants. By the dimension constraint (\ref{dimcons}), $I_d(h^i\cdot h^j\cdot h^k)$ is zero unless $3\ge i+j+k=rd+r+d$, which can only be achieved if $d\in\{0,1\}$. For the small quantum product we therefore have the simplification $I=I_0+I_+=I_0+I_1$. The ring given by this product is called the small quantum cohomology ring, and we shall compute it explicitly in section \ref{comprings}.
        \end{remark}

        \subsection{Gromov Witten Invariants and Quantum Cohomology for $\ponex$} %corresponding statements from "adapting the theory"
            Let us now examine the corresponding objects and statements in $\ponex$.
           
            %generators, GW ins etc are defined in the notes in the calculation of the small quantum ring of ponex on pages 10,11,12

            %big quantum ring can be reduced to small quantum ring: in this sense equiv. to Kontsevich formula

            We adopt the notation from the beginning of section \ref{th2qpr} with the obvious adjustments.
            Lemma \ref{interslemma} continues to hold in this setting.

            \begin{definition}
                We define the Gromov-Witten invariant (or GW-invariant) of bidegree $(d,e)$ associated to the cohomology classes $\gamma_1,\dots,\gamma_n\in A^*(\ponex)$ as follows:
                \begin{equation*}
                    I_{(d,e)}(\gamma_1\cdots\gamma_n):=\int_{\overline M}\underline\nu^*(\underline\gamma).
                \end{equation*}
            \end{definition}

            \begin{remark}[Dimension Constraint for $\ponex$]\label{dimconsponexrmk}
                This time, the constraint for $I_d(\gamma_1\cdots\gamma_n)$ to be non-zero looks as follows:
                \begin{equation}\label{dimconsponex}
                    \sum_i\text{codim}(\gamma_i)=\text{dim}(\monebar)=n+2d+2e-1.
                \end{equation}
            \end{remark}

            The following proposition corresponds to Corollary \ref{gwndpr} and is derived like Proposition \ref{gwmeaningpr}. In particular, if all the $\Gamma_i$ are points, Lemma \ref{lemma2} is strong enough and one needn't invoke an analogon of Lemma \ref{355}. %\textcolor{red}{shown likd gwmeaningpr, and it needs only analogon of 3.5.3, not 3.5.5, because all classes are points.}.
            \begin{proposition}\label{gwndeponex} %4.1.6' corresponds to cor for p2
                For $\ponex$, we have
                $I_{(d,e)}(T_3 \cdots T_3)=N_{(d,e)}$ for $2d+2e-1$ factors $T_3$.
            \end{proposition}

            % The invariants are interpreted similarly as in the $\pr$ case. \textcolor{red}{ADD STATEMENT AND ARGUMENT}
            % \begin{remark}\label{gwmeaningponex}
            %     We can make some further enumerative interpretations in some cases involving more than point classes.
            %     Indeed, let $n_1,n_2,n_3\in\Z_{\ge0}$
            %     Note that for these rules, we cannot simply extract the divisor classes using the divisor equation (Lemma \ref{divisoreqponex}). That case will be handled in Proposition \ref{ruleinvs}.
            %     %     Let $\gamma_1,\dots,\gamma_n\in A^*(\ponex)$ be homogeneous classes with \textcolor{red}{argue why codim constr unnecessary, cf pone. we must be able to allow codim greater equal 1, else problem.}
            %     %     %$\text{codim}(\gamma_i)\geq2$ for all $i$
            %     %     and $\sum_{i=1}^n\text{codim}(\gamma_i)=\text{dim}\bigl(\monebar\bigr)$. 
            %     %     Then $I_{(d,e)}(\gamma_1\cdots\gamma_n)$ is the number of rational bidegree $(d,e)$ curves incident to all the subvarieties $\Gamma_1$, \dots, $\Gamma_n\subset\ponex$.
            %     This extends the divisor equation for $\ponex$.
            % \end{remark}      

            \begin{lemma}\label{00inv}%4.2.1'
                $I_{(0,0)}=0$ unless $n=3$ and $\sum \text{codim}(\gamma_i)=2$.

                In this case we have:
                $$I_{(0,0)}(T_i\cdot T_j\cdot T_k)=\int (T_i\cup T_j\cup T_k)\cap[\ponex].$$
                Since $T_1\cup T_1=T_2\cup T_2=0$, this further implies that $I_{(0,0)}(T_i\cdot T_j\cdot T_k)=0$ if it contains two identical classes from $A^1(\ponex).$ 
            \end{lemma}
            
            \begin{proof}
                The proof is analogous to that of Lemma \ref{d=0inv}. This time we have $\sum \text{codim}(\gamma_i)=2=\dim\bigl(\overline{M}_{0,3}(\ponex,(0,0))\bigr)$.
                
            \end{proof}

            \begin{lemma}\label{1ptin}%4.2.2'
                For $n=1$, it holds that $I_{(d,e)}(\gamma_1)=0$ unless $(d,e)\in\{(0,1),(1,0)\}$ and $\gamma_1=T_3$, in which case we have $I_{(d,e)}(T_3)=1$.
            \end{lemma}
            
            \begin{proof}
                By Lemma \ref{00inv}, we may assume that $d>0$ or $e>0$. We treat the former case, the latter is handled analogously. The dimension constraint (\ref{dimconsponex}) leads to the following inequality:
                $$2\geq\codim \gamma_1)=\dim\bigl(\overline{M}_{0,1}(\ponex,(d,e))\bigr)=2(d+e)\geq2e+2.$$
                This implies that $e=0$ and $\codim \gamma_1)=2$, and further that $d+e=1$. The invariant is therefore $I_{(1,0)}(T_3)=1$ by Proposition \ref{gwndeponex}. The case $e>0$ yields $I_{(0,1)}(T_3)$.
                
            \end{proof} 
            \begin{remark}\label{2ptinv0}
                In the same way, one can show that for $n=2$, only invariants with  $(d,e)\in\{(0,1),(1,0)\}$ containing $T_3$ and either $T_1$ or $T_2$ can be non-zero.
            \end{remark}

            \begin{lemma}\label{fundclassinvponex}%4.2.3'
                %$I_{(d,e)}(\cdots T_0 \cdots)=0$ unless $(d,e)=(0,0)$ and $n=3$.
                Consider the fundamental class $1=T_0\in A^0(\ponex)$. If $\gamma_i=1$ for some $1\leq i\leq n$, then $I_{(d,e)}(\gamma_1\cdots\gamma_n)=0$ unless $(d,e)=(0,0)$ and $n=3$.
            \end{lemma}

            \begin{proof}
                The proof is analogous to that of Lemma \ref{fundclassinv}.
            \end{proof}

            \begin{lemma}[Divisor Equation for $\ponex$]\label{divisoreqponex} %4.2.4'
                Consider the rule classes $T_1,T_2\in A^1(\ponex)$.% If $d+e>0$, \textcolor{red}{is this enough for the degree argument? would look like I+, Gamma in pr }
                We can extract these classes from $I_{(d,e)}$ as follows: 
                
                If $e>0$ and $\gamma_i=T_1$ for some $1\leq i\leq n+1$, say $\gamma_{n+1}=T_1$, then: 
                $$I_{(d,e)}(\gamma_1\cdots\gamma_n\cdot T_1)=e\cdot I_{(d,e)}(\gamma_1\cdots\gamma_n).$$
                If $d>0$ and $\gamma_{n+1}=T_2$, then: 
                $$I_{(d,e)}(\gamma_1\cdots\gamma_n\cdot T_2)=d\cdot I_{(d,e)}(\gamma_1\cdots\gamma_n).$$
            \end{lemma}

            \begin{proof}
                We discuss the extraction of $T_1$, the other case is analogous. 
                The proof is analogous to that of Lemma \ref{divisoreq}, with a small twist due to the presence of bidegrees. In particular, if  $H\subset\ponex$ is the hyperplane satisfying $\tilde\nu_{n+1}^*(T_1)\cap[\overline{M}_{0,n+1}(\ponex,(d,e))]=[\tilde\nu_{n+1}^{-1}(H)]$, then restricting the forgetful map $\epsilon$ to $\tilde\nu_{n+1}^{-1}(H)$ yields a generically finite map $\epsilon|$ of degree $e$. Indeed, B\'ezout's theorem implies that the image curve of a general map $\mu\in\overline{M}_{0,n}(\ponex,(d,e))$ intersects $H$ in $(1,0)\circ(d,e)=e$ points. The rest of the argument is the same as in the proof of Lemma \ref{divisoreq}.
                
            \end{proof}

            % \begin{remark}\label{2ptin}%4.2.2''
            %     Using Lemma \ref{divisoreqponex}, Remark \ref{2ptinv0} and Proposition \ref{gwndeponex}, we can deduce that the only non-zero invariants with $n=2$ are $I_{(0,1)}(T_3)=I_{(1,0)}(T_3)=1$.
            % \end{remark}
            % \textcolor{red}{this remark is wrong with stronger restrictions in divisor equation}

            % \begin{lemma}[Splitting Lemma for $\ponex$]\label{splittingponex}
            %     Let $\alpha:D\hookrightarrow\overline M$ and $\iota:D\hookrightarrow \overline M_A\times \overline M_B$ be the inclusions. \textcolor{red}{must explain iota, gluing mark} Then for all $\gamma_1,\dots,\gamma_n\in A^*(\ponex)$, the following holds in $A^*(\overline M_A\times \overline M_B)$: 
            %     $$\iota_*\alpha^*\underline\nu^*(\underline\gamma)=\sum_{\codim T_e)+\codim T_f)=2}(\bigcup_{a\in A}\nu_a^*(\gamma_a)\cdot\nu_{x_A}^*(T_e))\times(\bigcup_{b\in B}\nu_b^*(\gamma_b)\cdot\nu_{x_B}^*(T_f)).$$
            %     In particular, we get the following integral formula:
            %     $$\int_D \nu_1^*(\gamma_1)\cup\cdots\cup\nu_n^*(\gamma_n) = \sum_{\codim T_e)+\codim T_f)=2} I_{(d_A,e_A)}(\prod_{a\in A}\gamma_a\cdot T_e)\cdot I_{(d_B,e_B)}(\prod_{b\in B}\gamma_b\cdot T_f).$$
            % \end{lemma}

            \begin{theorem}[Reconstruction for $\ponex$]\label{reconponex}
                Given that $I_{(1,0)}(T_3)=I_{(0,1)}(T_3)=1$, all genus $0$ GW-invariants for $\ponex$ can be computed recursively.
            \end{theorem}
            The idea for the proof is the same as for Theorem \ref{reconpr}. This proof requires an adapted version of the splitting lemma and the first recursion is over bidegrees instead of degrees. 

            % \begin{proof}
            %     %probably analogous. what are the three recursions? only the first one over (d,e) instead of just d should be more complicated. How to order the bidegrees? answer given by splitting lemma for ponex which gives the relevant equation. probably reducing e or d (equiv by symmetry) and passing to recursion.
            % \end{proof}

            %potentials
            \begin{propdef}[Collected GW-Invariants]
                The expression $$I(\gamma_1\cdots\gamma_n):=\sum_{\substack{(d,e)\\d+e>0}} I_{(d,e)}(\gamma_1\cdots\gamma_n)$$ is a finite sum. %this time more than one summand is non-zero. problem?
            \end{propdef}

            \begin{proof}
                The dimension constraint (\ref{dimconsponex}) implies that $I_{(d,e)}(\gamma_1\cdots\gamma_n)=0$ unless $\sum_i\text{codim}(\gamma_i)=\text{dim}(\monebar)=n+2d+2e-1$. Thus only the summands satisfying $$d+e=\frac{\sum \text{codim}(\gamma_i)-n+1}{2}$$ survive.  
            \end{proof}

            Again, the linearity of the invariants allows us to assume that the input classes are given by basis elements $T_i\in A^*(\ponex)$. The collected invariants are thus w.l.o.g. given by $I((T_0)^{\bullet a_0}\cdots (T_3)^{\bullet a_3}).$
            We again define the corresponding generating function.

            \begin{definition}[Potentials]
                Using formal variables $x_0,x_1,x_2,x_3$, we define the GW-potential for $\ponex$:
                $$\Phi(x_0,x_1,x_2,x_3):=\sum_{(a_0,a_1,a_2,a_3)}\frac{x_0^{a_0}x_1^{a_1}x_2^{a_2} x_3^{a_3}}{a_0!a_1!a_2!a_3!}I((T_0)^{\bullet a_0}(T_1)^{\bullet a_1}(T_2)^{\bullet a_2} (T_3)^{\bullet a_3}).$$
                Rewriting $\mathbf{x}:=(x_0,x_1,x_2,x_3)$ and $\mathbf{T}^{\mathbf{a}}:=(T_0)^{\bullet a_0}(T_1)^{\bullet a_1}(T_2)^{\bullet a_2} (T_3)^{\bullet a_3}$, we can write this as $\Phi(\mathbf{x})=\sum_{\mathbf{a}}\frac{\mathbf{x}^{\mathbf{a}}}{\mathbf{a}!}I(\mathbf{T}^{\mathbf{a}})$.
            \end{definition}
           
            We decompose $I$ into $I_{(0,0)}+I_+$. We define the classical potential $\Phi^{\text{cl}}$ and the quantum potential $\Gamma$ accordingly.
            By Lemma \ref{00inv}, the classical potential is given by:
            $$\Phi^{\text{cl}}(\mathbf{x})=\sum_{(i,j,k)}\frac{x_ix_jx_k}{3!}I_{(0,0)}(T_i\cdot T_j\cdot T_k),$$
            where the sum of the codimensions of the inputs is equal to $2$ and the factor $3!$ again compensates for the different permutations of $i,j,k$.

            As before, \ref{genfctder} implies that $\Phi_i=\frac{\delta}{\delta x_i}\Phi=\sum_{\mathbf{a}}\frac{\mathbf{x}^\mathbf{a}}{\mathbf{a}!}I(\mathbf{h}^\mathbf{a}\cdot h^i).$ 
            In particular, we again deduce that:
            \begin{equation}\label{phiijkponex}
                \Phi_{ijk}(\mathbf{x})=\sum_{\mathbf{a}}\frac{\mathbf{x}^\mathbf{a}}{\mathbf{a}!}I(\mathbf{T}^\mathbf{a}\cdot T_i \cdot T_j\cdot T_k).
            \end{equation}

            \begin{proposition}\label{classprodponex}     %5.2.6'
                The structure constants of the cup product in $A^*(\ponex)$ are given by the third derivatives of the classical potential. More precisely, for arbitrary classes $T_i,T_j\in A^*(\ponex)$, it holds that:
                \begin{equation*}
                    T_i\cup T_j = \sum_{f=0}^3I_{(0,0)}(T_i\cdot T_j\cdot T_{3-f})T_f=\sum_{k+f=3}\Phi_{ijk}^{\text{cl}}T_f.
                \end{equation*}
                We shall therefore also refer to this cup product as the classical product.
            \end{proposition}

            \begin{proof}
                This can be verified explicitly using Lemma \ref{00inv}. Note that unlike in the proof of Proposition \ref{classprodpr}, we must be careful since the codimension of $T_i$ is not necessarily equal to $i$.
                
                Consider the expression $I_{(0,0)}(T_i\cdot T_j\cdot T_{3-f})T_f$. In order for the coefficient $I_{(0,0)}(T_i\cdot T_j\cdot T_{3-f})$ not to vanish, the sum of the codimensions of the three inputs must be equal to $2$ and there cannot be two identical codimension $1$ classes. This is enough to verify the equality in each case.
                
            \end{proof}

            \begin{definition}[Quantum Product]
                The quantum product of two classes $T_i,T_j\in A^*(\ponex)$ is defined as:
                $$T_i*T_j:=\sum_{e+f=3}\Phi_{ije}T_f\in \Q[[\mathbf{x}]] \otimes_{\Z}A^*(\ponex),$$
                with $e,f\in\Z_{\ge0}$. We define the quantum product on all of $\Q[[\mathbf{x}]] \otimes_{\Z}A^*(\ponex)$ by $\Q[[\mathbf{x}]]$-linear extension.
            \end{definition}

            By Proposition \ref{classprodponex}, we can decompose the quantum product as follows:
            $$T_i*T_j=T_i\cup T_j+\sum_{e+f=3}\Gamma_{ije}h^f.$$
            
            %5.2.8' omit(?)

            \begin{lemma}%5.2.9'
                The identity element for the quantum product is given by the fundamental class $T_0$.
            \end{lemma}

            \begin{proof}
                For all $i$, the definition of the quantum product says that $T_0*T_i=\sum_{e+f=3}\Phi_{0ie}T_f$, where 
                $\Phi_{0ie}=\sum_\mathbf{a}\frac{\mathbf{x}^\mathbf{a}}{\mathbf{a}!}I(\mathbf{T}^{\mathbf{a}}\cdot T_0 \cdot T_i\cdot T_e)$. Using Lemma \ref{fundclassinvponex}, we deduce that $\Phi_{0ie}=I_{(0,0)}(T_0 \cdot T_i\cdot T_e)$ and can therefore rewrite the initial product as follows:
                $$T_0*T_i=\sum_{e+f=3}I_{(0,0)}(T_0 \cdot T_i\cdot T_e)T_f,$$
                where $\codim T_i)+\codim T_e)=2$. As in the proof of Proposition \ref{classprodponex}, these conditions are enough to verify that $T_0*T_i=T_i$ for $0\leq i\leq 3$.
                
            \end{proof}

            \begin{proposition}[Associativity of the Quantum Product in $\ponex$]\label{assponex}
            For arbitrary classes $T_i,T_j,T_k\in A^*(\ponex)$, the following equality holds:
                $$(T_i*T_j)*T_k=T_i*(T_j*T_k).$$
            \end{proposition}

            The proof is analogous to the $\pr$ case, see Proposition \ref{asspr}. It also relies on an adapted splitting lemma for $\ponex$.

            \begin{remark}[Big and small Quantum Cohomology for $\ponex$]
                The quantum product again determines a quantum cohomology ring which we denote by $Q^*(\ponex)$.
                The small quantum cohomology is obtained by setting the variables $x_1$ and $x_2$ corresponding to the divisor classes $T_1$ and $T_2$ to zero. We compute the small quantum cohomology in detail in Proposition \ref{qringponex}.

                % A simpler yet useful variation of this quantum cohomology is obtained by setting all variables in $\Phi_{ijk}(\mathbf{x})$ to zero except for those corresponding to codimension $1$ classes, i.e. $x_1$. We can rewrite $x:=x_1$ and simplify the $\Phi_{ijk}$ using Lemma \ref{divisoreq}:
                % \begin{align*}
                %     \Phi_{ijk}(x) & = \sum_{n=0}^\infty \frac{x^n}{n!}\sum_{d\ge0} I_d(h^{\bullet n}\cdot h^i\cdot h^j\cdot h^k) \\
                %     & = \sum_{n=0}^\infty \frac{x^n}{n!}\sum_{d\ge0}d^n\cdot I_d(h^i\cdot h^j\cdot h^k).
                % \end{align*}
                % We call the product given by these structure constants the small quantum product. By the above calculation, it involves only three-point invariants. By the dimension constraint (\ref{dimcons}), $I_d(h^i\cdot h^j\cdot h^k)$ is zero unless $3\ge i+j+k=rd+r+d$, which can only be achieved if $d\in\{0,1\}$. For the small quantum product we therefore have the simplification $I=I_0+I_+=I_0+I_1$. The ring given by this product is called the small quantum cohomology ring, and we shall compute it explicitly in Section \ref{comprings}.
            \end{remark}
\chapter{Computations}\label{comp} %titles to be adapted
In order to familiarize ourselves with the newly introduced theory, we perform some explicit calculations. 
    \section{GW-Invariants}\label{compgwi}
    We begin by computing some GW-invariants in our various target spaces. 

        \subsection{All GW-invariants of $\pone$} %separatre because of reason mentioned in theory 2
        We previously alluded to the fact that the $\pone$ case is especially simple. We now justify this claim.
            \begin{proposition}\label{p1invs}
                The only non-zero GW-invariants for $\pone$ are $I_0(h^1\cdot h^0\cdot h^0)=I_1((h^1)^{\bullet n})=1$.
            \end{proposition}
    
            \begin{proof}
                For $d=0$, we must have $n=3$ and $\sum\text{codim}(\gamma_i)=r=1$ by Lemmas \ref{d=0inv} and \ref{fundclassinv}. The only such invariant is $I_0(h^1\cdot h^0\cdot h^0)$.

                For $d>0$, we write $h:=h^1$ and denote by $n$ the number of marks. By definition, we must have:
                $$n\geq \sum\text{codim}(\gamma_i)=\text{dim}(\overline{M}_{0,n}(\pone,d)=2d-2+n\geq n.$$
                Thus $2d-2+n=n=\sum\text{codim}(\gamma_i)$, so $d=1$ and $\gamma_i=h$ for all $1\leq i\leq n$. The only non-zero invariants of degree $d=1$ in $\pone$ are thus $I_1(h^{\bullet n})=1$ for any $n\geq1$. They take the value 1 because there is precisely one line through $n$ points in $\pone$, namely $\pone$ itself. 
                
                Note that this last argument uses an adapted version of Proposition \ref{gwmeaningpr} for the case where $r=1$ and the $\Gamma_i$ are points, the proof of which is a simpler version of that of the original proposition, since it no longer needs Lemma \ref{qijcodim} and therefore no lower bound on the codimensions of the $\Gamma_i$.
         
            \end{proof}

            \begin{proposition}
                Let $D_{AB}=D(A,B;d_A,d_B)\subset \overline{M}_{0,n}(\pone,2)$ be a boundary divisor. Then for any tuple $(\gamma_1,\dots,\gamma_n)$ of cohomology classes, it holds that:
                $$\int_{D_{AB}} \nu_1^*(\gamma_1)\cup\dots\cup\nu_n^*(\gamma_n)=0.$$
            \end{proposition}
    
            \begin{proof}
                An application of Lemma \ref{splitting} will yield the result (since all degree $2$ invariants in $\pone$ are zero, we must only consider the case $d_A=1=d_B$): 
                \begin{align*}
                    & \int_{D_{AB}} \nu_1^*(\gamma_1)\cup\dots\cup\nu_n^*(\gamma_n)\\
                     = & \sum_{e+f=1} I_1(\prod_{a\in A}\gamma_a\cdot h^e)\cdot I_1(\prod_{b\in B}\gamma_b\cdot h^f) \\
                     = & I_1(\prod_{a\in A}\gamma_a\cdot h^1)\cdot I_1(\prod_{b\in B}\gamma_b\cdot h^0) + I_1(\prod_{a\in A}\gamma_a\cdot h^0)\cdot I_1(\prod_{b\in B}\gamma_b\cdot h^1).
                \end{align*}
                But by Proposition \ref{p1invs}, we have $I_1(\prod_{b\in B}\gamma_b\cdot h^0)=I_1(\prod_{a\in A}\gamma_a\cdot h^0)=0$. This concludes the proof.
                
            \end{proof}

        \subsection{Three-Point Invariants of $\pr$} % of $\pr$ for $r\geq 2$.  excs4.1-4.4, applications of recon algo

            Consider now a three-point invariant $I_d(h^i\cdot h^j \cdot h^k)$ in $\pr$. By definition, it holds that: 
            $$3r\geq i+j+k=\text{dim}(\overline{M}_{0,3}(\pr,d))=d(r+1)+r,$$
            where we used equation (\ref{dimmbar}). But this is impossible unless $d\leq1$. For $d=0$, the above condition reads $i+j+k=r$. For $d=1$, it becomes $i+j+k=2r+1$.

            The following proposition is an opportunity to see the Reconstruction Algorithm \ref{reconpr} in action.
            \begin{proposition}\label{3ptp3p4}
                The following assertions hold:
                \begin{itemize}
                    \item In $\mathbb{P}^3$, we have $I_1(h^3\cdot h^2\cdot h^2)=1$.
                    \item In $\mathbb{P}^4$, we have $I_1(h^4\cdot h^3\cdot h^2)=1$ and $I_1(h^3\cdot h^3\cdot h^3)=1$.
                \end{itemize}
            \end{proposition}

            \begin{proof}
                We calculate $I_1(h^3\cdot h^2\cdot h^2)$ in $\mathbb{P}^3$ proceeding as in the proof of the Reconstruction Theorem. The three classes are already ordered by decreasing codimension. We rewrite $\gamma_n=h^2$ as $\lambda_1\cup \lambda_2=h^1\cup h^1$.
                In this setting, equation (\ref{reconpreq}) looks as follows:
                $$I_0(h^1\cdot h^1\cdot h^1)\cdot I_1(h^3\cdot h^2\cdot h^2)=I_1(h^3 \cdot h^1\cdot h^3)\cdot I_0(h^2\cdot h^1\cdot h^0),$$
                where $I_0(h^1\cdot h^1\cdot h^1)=I_0(h^2\cdot h^1\cdot h^0)=1$ by Lemma \ref{d=0inv} and $I_1(h^3 \cdot h^1\cdot h^3)=I_1(h^3 \cdot h^3)\cdot 1=1$ by Lemma \ref{divisoreq}. Thus $I_1(h^3\cdot h^2\cdot h^2)=1$ as claimed.

                We proceed analogously for $I_1(h^4\cdot h^3\cdot h^2)=1$ in $\mathbb{P}^4$. Again we write $h^2=h_1\cup h_1$ and consider equation (\ref{reconpreq}) in this situation:
                $$I_0(h^1\cdot h^1\cdot h^2)\cdot I_1(h^4\cdot h^3\cdot h^2)=I_1(h^1 \cdot h^4\cdot h^4)\cdot I_0(h^1\cdot h^3\cdot h^0).$$
                Both degree $0$ invariants are equal to $1$, and $I_1(h^1 \cdot h^4\cdot h^4)$ recovers $I_1(h^4\cdot h^4)=1$. We conclude that $I_1(h^4\cdot h^3\cdot h^2)=1$.

                For $I_1(h^3\cdot h^3\cdot h^3)$, the equation becomes: 
                $$I_0(h^1\cdot h^2\cdot h^1)\cdot I_1(h^3\cdot h^3\cdot h^3)=I_1(h^2 \cdot h^3\cdot h^4)\cdot I_0(h^1\cdot h^3\cdot h^0).$$
                The degree $0$ invariants are equal to $1$, and we just showed that $I_1(h^4\cdot h^3\cdot h^2)=1$. Therefore $I_1(h^3\cdot h^3\cdot h^3)=1$.
                
            \end{proof}

            We now generalize these statements.
            \begin{proposition}\label{3ptpr}
                For $r\geq2$, all non-zero three-point invariants $I_d(\gamma_1,\gamma_2,\gamma_3)$ of $\pr$ are equal to 1.
            \end{proposition}

            \begin{proof}%adapted from pf of recursion algo
                Since $n=3$ is fixed, we need only retrace the first two recursions from the proof of Theorem \ref{reconpr}, namely those over the degree $d$ and over the smallest codimension of a class appearing in the invariant (i.e. $\text{codim}(\lambda_i)$ in the language of that proof). We will reformulate these recursions into inductions. 
                %Assume w.l.o.g. that $\gamma_1,\gamma_2,\gamma_3$ are ordered by decreasing codimension.

                We first do induction over the degree $d$. From Lemma \ref{d=0inv}, we know that the statement holds in the case $d=0$. Pick a $d>0$ and assume that all non-zero three-point invariants of smaller degree are equal to $1$. Set $b_i:=\text{codim}(\gamma_i)$ for $i\in\{1,2,3\}$.
                Assume w.l.o.g. that $b_1\geq b_2 \geq b_3$. We rewrite $\gamma_3$ as $\lambda_1\cup \lambda_2$, with $c_i:=\text{codim}(\lambda_i)$ for $i\in\{1,2\}$ and $\text{codim}(\gamma_3)>c_2\geq c_1$ (w.l.o.g.).

                We now begin our second, 'internal' induction over $\text{codim}(\gamma_3)$. If $\text{codim}(\gamma_3)=0$, then we know by Lemma \ref{fundclassinv} that all invariants are zero, so the statement holds trivially. By Lemma \ref{divisoreq}, we also know that the claim holds when $\text{codim}(\gamma_3)=1$. Now let $\text{codim}(\gamma_3)>1$ and assume that all three-point invariants with a class of smaller codimension are equal to $1$. 

                We again look at equation (\ref{reconpreq}) in this setting. 
                By the induction hypothesis on $d$, all summands with $d_A<d$ and $d_B<d$ are equal to $1$.
                We examine the case $d_A=d$ or $d_B=d$. Arguing as in the proof of Proposition \ref{3ptp3p4} and writing $\gamma_3=h^{c_1+c_2}$, we find summands $I_0(\lambda_1\cdot\lambda_2\cdot h^{r-c_1-c_2})\cdot I_d(\gamma_1\cdot \gamma_2 \cdot h^{c_1+c_2})$ and $I_d(\lambda_1\cdot\lambda_2\cdot h^{b_1+b_2})\cdot I_0(\gamma_1\cdot \gamma_2 \cdot h^{r-b_1-b_2})$ on the left hand side; and $I_0(\lambda_1\cdot\gamma_1\cdot h^{r-c_1-b_1})\cdot I_d(\lambda_2\cdot \gamma_2 \cdot h^{c_1+b_1})$ and $I_d(\lambda_1\cdot\gamma_1\cdot h^{c_2+b_2})\cdot I_0(\lambda_2\cdot \gamma_2 \cdot h^{r-c_2-b_2})$ on the right hand side. Applying Lemma \ref{d=0inv} and reordering the classes turns these summands into $I_d(\gamma_1\cdot \gamma_2 \cdot h^{c_1+c_2})$, $I_d(h^{b_1+b_2} \cdot\lambda_2\cdot \lambda_1)$, $I_d(h^{c_1+b_1}\cdot\gamma_2\cdot\lambda_2)$, and $I_d(h^{c_2+b_2}\cdot\gamma_1\cdot\lambda_1)$, respectively. The first summand is the invariant we would like to prove is equal to $1$, and the remaining three are equal to $1$ by the induction hypothesis on $\text{codim}(\gamma_3)$.

                To conclude, note that both sides of equation (\ref{reconpreq}) have the same number of summands. By the above arguments, all summands on the right hand side are equal to $1$. On the left hand side, one summand is equal to $I_d(\gamma_1,\gamma_2,\gamma_3)$ and all the others are equal to $1$ as well. Thus $I_d(\gamma_1,\gamma_2,\gamma_3)=1$, and this concludes both of our inductions and thereby the proof.
                 
            \end{proof}

            % \textcolor{red}{add corresponding result abt 3pt invariants in ponex needed for calc of small q ring here or keep it down in ring comp? ie result about certaain invs being equal to 1}

        \subsection{Invariants of Rules in $\ponex$}\label{ruleinv}
            %TODO (maybe): choose more props: isolate computations from later calculations
            We can easily calculate invariants of bidegrees corresponding to (single or multiple) horizontal or vertical rules in cases where the divisor equation \ref{divisoreqponex} is applicable.
            \begin{proposition}\label{ruleinvs} %from quantum notes p.14
                Let $n\geq0$. Then $I_{(0,1)}(T_1^{\bullet n}\cdot T_3)=I_{(1,0)}(T_2^{\bullet n}\cdot T_3)=1$, and for any $d>1$ and $e>1$, it holds that $I_{(0,e)}(T_1^{\bullet n}\cdot T_3^{2e-1})= I_{(d,0)}(T_2^{\bullet n}\cdot T_3^{2d-1})=0$.
                % \begin{itemize}
                %     \item $I_{(0,1)}(T_1^{\bullet n}\cdot T_3)=1$,
                %     \item $I_{(1,0)}(T_2^{\bullet n}\cdot T_3)=1$,
                % \end{itemize}
                % and for any $d>1$, $e>1$:
                % \begin{itemize}
                %     \item $I_{(0,e)}(T_1^{\bullet n}\cdot T_3^{2e-1})=0$,
                %     \item $I_{(d,0)}(T_2^{\bullet n}\cdot T_3^{2d-1})=0$.
                % \end{itemize}
            \end{proposition}

            \begin{proof}
                In each case, we apply first Lemma \ref{divisoreqponex}, then Proposition \ref{gwndeponex} and finally Fact \ref{n10} (for the first two calculations), respectively Proposition \ref{nd0} (for the last two).
                \begin{itemize}
                    \item $I_{(0,1)}(T_1^{\bullet n}\cdot T_3)=I_{(0,1)}(T_3)=N_{(0,1)}=1$,
                    \item $I_{(1,0)}(T_2^{\bullet n}\cdot T_3)=I_{(1,0)}(T_3)=N_{(1,0)}=1$,
                    \item $I_{(0,e)}(T_1^{\bullet n}\cdot T_3^{2e-1}) = e^n \cdot N_{(0,e)}=0$,
                    \item $I_{(d,0)}(T_2^{\bullet n}\cdot T_3^{2d-1}) = d^n \cdot N_{(d,0)}=0$.
                \end{itemize}
                
            \end{proof}

    \section{Potentials}\label{comppot}
        %\subsection{Potentials in $\pr$ for $r\in \{1,2,3\}$}

            %\begin{proposition}[Classical Potentials]
            We compute the classical potentials for $\pone$, $\ptwo$, $\mathbb{P}^3$ and $\ponex$.

            For $\pone$, we have: 
            $$\Phi_{\pone}^{\text{cl}}(x_0,x_1)=\sum_{i,j,k}\frac{x_i x_j x_k}{3!}I_0(h^i\cdot h^j\cdot h^k)=3\cdot\frac{x_0^2 x_1}{3!}=\frac{1}{2}x_0^2 x_1,$$
            where we used Lemma \ref{d=0inv} (or Proposition \ref{p1invs}) and the factor $3$ appears as the number of selections of $(i,j,k)$ such that $i+j+k=1$, i.e. the number of permutations of $(1,0,0)$.

            For $\ptwo$, we have:
            $$ \Phi_{\ptwo}^{\text{cl}}(x_0,x_1,x_2)=\sum_{i,j,k}\frac{x_i x_j x_k}{3!}I_0(h^i\cdot h^j\cdot h^k)=3\cdot\frac{x_0 x_1^2}{3!} + 3\cdot\frac{x_0^2 x_2}{3!} =\frac{1}{2}(x_0^2 x_1 + x_0^2 x_2),$$
            which again follows from Lemma \ref{d=0inv}. The first factor $3$ is the number of permutations of $(1,1,0)$ and the second is the number of permutations of $(2,0,0)$.

            For $\mathbb{P}^3$, we proceed analogously:
            $$\Phi_{\mathbb{P}^3}^{\text{cl}}(x_0,x_1,x_2,x_3)=\frac{x_1^3}{3!}+3!\cdot\frac{x_0 x_1 x_2}{3!} =\frac{1}{6}x_1^3 + x_0 x_1 x_2,$$
            since there is only one permutation of $(1,1,1)$ whereas there are $3!=6$ of $(2,1,0)$.

            Finally, for $\ponex$, we just use equation (\ref{phiijkponex}) and the previously mentioned constraints on the inputs (i.e. that their codimensions must add up to two and that there can be no two identical divisor classes). We find the following expression:
            $$\Phi_{\ponex}^{\text{cl}}(x_0,x_1,x_2,x_3)=\frac{1}{2}x_0^2x_3+x_0x_1x_2,$$
            since there are $3$ was to permute $(0,0,3)$ and $3!$ ways to permute $(0,1,2)$.

            From Proposition \ref{p1invs}, we know all the GW-invariants of $\pone$. We can therefore explicitly determine its GW-potential.
            
            \begin{proposition}[GW-potential for $\pone$] % add $r\in \{1,2,3\}$?
                The GW-potential for $\pone$ is given by the following expression:
                $$\Phi_{\pone}(x_0,x_1)=\frac{1}{2}x_0^2 x_1 + exp(x_1).$$
            \end{proposition}

            \begin{proof}
                We decompose the GW-potential as $\Phi_{\pone}(x_0,x_1)=\Phi_{\pone}^{\text{cl}}(x_0,x_1) + \Gamma_{\pone}(x_0,x_1)$, where the quantum potential $\Gamma_{\pone}$ is given by $\sum_{\mathbf{a}}\frac{\mathbf{x}^{\mathbf{a}}}{\mathbf{a}!} I_+(\mathbf{h}^{\mathbf{a}})$. Using Proposition \ref{p1invs}, we can simplify this sum to  $\sum_{n}\frac{x_1^{n}}{n!} I_1(h^{\bullet n})=\sum_{n}\frac{x_1^{n}}{n!}=\text{exp}(x_1)$.

                This combined with $\Phi_{\pone}^{\text{cl}}(x_0,x_1)=\frac{x_0^2 x_1}{2}$ from above yields the result.
                
            \end{proof}

            %find link between p1 and p1xp1. can we compute all gwis as in p1 case? is there some kind of functor?            

    \section{Rings}\label{comprings}
        %\subsection{$\pr$}
            %add?
            % \begin{proposition}[Cohomology Ring of $\ptwo$]   
            % \end{proposition}
            In this section, we calculate the small quantum rings of $\pr$ and $\ponex$. In the case of $\pr$, our knowledge of its three-point invariants will play an important role.

            \begin{proposition}[Cohomology Ring of $\pr$]   
                It holds that: 
                $$H^*(\pr)\cong \frac{\Q[h]}{(h^{r+1})}.$$
            \end{proposition}
    
            \begin{proof}
                By definition, we have $H^*(\pr)=\langle h^0, h^1,\dots,h^r\rangle$. Writing $1=h^0$ and $h=h^1$, we have again by definition that $h\cup1=h$, $h^i\cup h^j=h^{i+j}$, and $h^i=0$ for $i>r$. This implies the statement.
                
            \end{proof}
    
            \begin{proposition}[Small Quantum Ring of $\pr$]\label{qringpr}  
                It holds that:
                $$Q^*(\pr)\cong \frac{\Q[h,q]}{(h^{r+1}-q)}.$$
            \end{proposition}
    
            \begin{proof}
                We must compute the structure constants $\Phi_{ijk}(x):=\Phi_{ijk}(\mathbf{x})|_{x_0=x_2=\dots=x_n=0}$ of the (small) quantum product.
                Using equation (\ref{phiijk}) in combination with Lemma \ref{divisoreq} in this setting yields that $\Phi_{ijk}(x)=\sum_{n=0}^\infty\frac{x^n}{n!}\sum_d^\infty d^n I_d(h^i\cdot h^j\cdot h^k)$. By the dimension constraint (\ref{dimconspr}), we have that $I_d(h^i\cdot h^j\cdot h^k)=0$ unless: 
                $$(d+1)r+d=\text{dim}(\overline{M}_{0,3}(\pr,d))=\text{codim}(h^i)+\text{codim}(h^j)+\text{codim}(h^k)=i+j+k.$$
                Since $i+j+k\leq3r$, this implies that $d\in\{0,1\}$. The above sum can then be simplified to:
                $$I_0(h^i\cdot h^j\cdot h^k) + \sum_{n=0}^\infty\frac{x^n}{n!} I_1(h^i\cdot h^j\cdot h^k).$$
                The dimension constraint reads $i+j+k=r$ for the degree $0$ term and $i+j+k=2r+1$ for the degree $1$ term. Set $q:= \text{exp}(x)=\sum_{n=0}^\infty\frac{x^n}{n!}$.

                We are now ready to compute the quantum product of two generators $h^i$ and $h^j$.
                This is by definition equal to $h^i*h^j=\sum_{e+f=r}\Phi_{ije}h^f=\sum_{f=0}\Phi_{ij(r-f)}h^f.$
                From our computations and dimension constraints above (applied to the case $k=(r-f)$), we deduce that only the degree $0$ term with $f=i+j$ and the degree $1$ term with $f=i+j-r-1$ survive. Since $0\leq f\leq r$, this implies that:
                \begin{equation*}
                    h^i*h^j=
                    \begin{cases}
                      I_0(h^i\cdot h^j\cdot h^{r-(i+j)})h^{i+j}, & \text{if}\ i+j\leq r \\
                      q\cdot I_1(h^i\cdot h^j\cdot h^{2r+1-(i+j)})h^{i+j-r-1}, & \text{if}\ r<i+j\leq 2r+1.
                    \end{cases}
                \end{equation*}
                By Proposition \ref{3ptpr}, this means that:
                \begin{equation*}
                    h^i*h^j=
                    \begin{cases}
                      h^{i+j}, & \text{if}\ i+j\leq r \\
                      q\cdot h^{i+j-r-1}, & \text{if}\ r+1\leq i+j\leq 2r+1.
                    \end{cases}
                \end{equation*}
                In particular, by the associativity of the (small) quantum product (Proposition \ref{asspr}), the $(r+1)$-fold power $h*\dots *h$ is equal to $h*h^r=q\cdot h^{1+r-r-1}=q$. This implies the statement of the proposition.
                
            \end{proof}

            \begin{remark}
                The quantum product is not a deformation of the classical product, in the sense that setting $\mathbf{x}=\mathbf{0}$ in the structure constants of the quantum product does not retrieve the classical product. Indeed, the product with $\mathbf{x}=\mathbf{0}$ corresponds to the ring $\frac{\Q[h]}{(h^{r+1}-1)}$, i.e. the above case with $q=1$.
            \end{remark}

        %\subsection{$\ponex$}

            \begin{proposition}[Cohomology Ring of $\ponex$]   
                It holds that: 
                $$H^*(\ponex)\cong \frac{\Q[h,v]}{(h^2,v^2)}.$$
            \end{proposition}
    
            \begin{proof}
                By definition, $H^*(\ponex)=\langle T_0, T_1, T_2, T_3\rangle$. Recall that $\text{codim}(T_0)=0$, $\text{codim}(T_1)=\text{codim}(T_2)=1$, and $\text{codim}(T_3)=2$ (in particular, the subscripts do not correspond to the codimensions, contrary to the superscripts for $\pr$). More precisely, $T_0$ is the fundamental class, $T_1$ and $T_2$ are classes of a vertical and a horizontal rule, respectively, and $T_3$ is the point class.

                By definition of the cup product as the Poincar\'e dual of the intersection product, we have $T_1\cup T_2= T_3$, $T_0\cup T=T$ for any class $T$, and $T_i\cup T_i=0$ unless $i=0$. This concludes the proof.
            \end{proof}

            As mentioned in Lemma \ref{00inv}, this implies that any bidegree $(0,0)$ invariant containing more than one of either $T_1$ or $T_2$ is zero. 

            \begin{proposition}[Small Quantum Ring of $\ponex$]\label{qringponex}   
                It holds that:
                $$Q^*(\ponex)\cong \frac{\Q[h,v,q_h,q_v]}{(h^2-q_h, v^2-q_v)}.$$
            \end{proposition}
    
            \begin{proof}
                We compute the structure constants $\Phi_{ijk}(x_1,x_2)=\Phi_{ijk}(x_0,x_1,x_2,x_3)|_{x_0=x_3=0}$ of the quantum product. By definition, they are given by: 
                $$\sum_{n_1=0}^\infty \sum_{n_2=0}^\infty \frac{x_1^{n_1}x_2^{n_2}}{n_1!n_2!}\sum_{(d,e)} I_{(d,e)}(T_1^{\bullet n_1}\cdot T_2^{\bullet n_2} \cdot T_i \cdot T_j\cdot T_k).$$
                Lemma \ref{00inv} allows us to split up this sum into $I_{(0,0)}(T_i\cdot T_j\cdot T_k) + \Gamma_{ijk}(x_1,x_2)$, where
                $$\Gamma_{ijk}(x_1,x_2)=\sum_{n_1=0}^\infty \sum_{n_2=0}^\infty \frac{x_1^{n_1}x_2^{n_2}}{n_1!n_2!}\sum_{\substack{(d,e) \\ d+e>0}} I_{(d,e)}(T_1^{\bullet n_1}\cdot T_2^{\bullet n_2} \cdot T_i \cdot T_j\cdot T_k).$$

                We compute the quantum product of two arbitrary generators $T_i$ and $T_j$, i.e. $T_i*T_j=T_i\cup T_j + \sum_{k+f=3}\Gamma_{ijk}(x_1,x_2)T_f$. 

                Dimension constraint (\ref{dimconsponex}) requires that:
                % The following dimension constraint applies to $\Gamma_{ijk}$:
                $$\text{dim}\bigl(\overline{M}_{0,n_1+n_2+3}(\ponex,(d,e))\bigr)=n_1+n_2+\text{codim}(T_i)+\text{codim}(T_j)+\text{codim}(T_k),$$
                and the left hand side is equal to $n_1+n_2+2(d+e+1)$ by equation (\ref{dimmbarone}).

                The resulting constraint, namely $2(d+e+1)=\text{codim}(T_i)+\text{codim}(T_j)+\text{codim}(T_k)\leq6$, implies that the sum of the three codimensions is even and that $d+e\leq2$.
                
                Let us now examine $\Gamma_{ijk}$ in more detail. The constraint reveals that the sum over $(d,e)$ with $d+e>0$ is actually just a sum over $d+e \in\{1,2\}$, namely:
                $$\sum_{d=1}^2I_{(d,0)}(T_1^{\bullet n_1}\cdot T_2^{\bullet n_2} \cdot T_i \cdot T_j\cdot T_k)\ + \ \sum_{e=1}^2I_{(0,e)}(T_1^{\bullet n_1}\cdot T_2^{\bullet n_2} \cdot T_i \cdot T_j\cdot T_k)\ +\ I_{(1,1)}(T_1^{\bullet n_1}\cdot T_2^{\bullet n_2} \cdot T_i \cdot T_j\cdot T_k).$$

                Applying Lemma \ref{divisoreqponex} (the divisor equation) to each of the three summands leads to:
                \begin{equation}\label{summsumm}
                     \sum_{d=1}^2d^{n_2}I_{(d,0)}(T_1^{\bullet n_1} \cdot T_i \cdot T_j\cdot T_k)\ + \ \sum_{e=1}^2e^{n_1}I_{(0,e)}( T_2^{\bullet n_2} \cdot T_i \cdot T_j\cdot T_k)\ +\ I_{(1,1)}(T_i \cdot T_j\cdot T_k).
                \end{equation}
                Now since $(d,e)\neq(0,0)$ in $\Gamma_{ijk}$, we know by Lemma \ref{fundclassinvponex} that only $T_i$, $T_j$ and $T_k$ of positive codimension can appear. The only potentially non-zero $\Gamma_{ijk}$ are thus $\Gamma_{113}$, $\Gamma_{123}$, $\Gamma_{223}$ and $\Gamma_{333}$.

                We examine each of these four possibilities.

                Ad $\Gamma_{113}$:
                Let's see what becomes of our sum (\ref{summsumm}). 
                Consider the invariant $I_{(d,0)}(T_1^{\bullet n_1} \cdot T_1 \cdot T_1\cdot T_3)$. The codimensions of the inputs add up to $n_1+4$, and by the dimension constraint (\ref{dimconsponex}), this must be equal to $n_1+2d+2$. This implies that the invariant is zero unless $d=1$. The same goes for the invariant $I_{(0,e)}( T_2^{\bullet n_2} \cdot T_1 \cdot T_1\cdot T_3)$, where we must have $e=1$. The sum therefore becomes:
                
                \begin{equation*}
                    I_{(1,0)}(T_1^{\bullet n_1+2}\cdot T_3)\ + \ I_{(0,1)}( T_2^{\bullet n_2} \cdot T_3)\ +\ I_{(1,1)}(T_3),
                \end{equation*}
                    
                where we appplied the divisor equation to the second and third summands. The third summand is zero since $I_{(1,1)}(T_3)$ fails to obey the dimension constraint: otherwise we would have $4=\text{dim}(\overline{M}_{0,1}(\ponex,(1,1)))=\text{codim}(T_3)=2$, see equation (\ref{dimmbarone}).
                By Lemma \ref{interslemma}, the first summand is zero and the second can only survive if $n_2=0$. Indeed, the image curve of a bidegree $(d,0)$ map is determined by $d$ points, and imposing another intersection condition thus yields zero maps.

                The sum therefore evaluates to $I_{(0,1)}(T_3)=N_{(0,1)}=1$.
                We conclude that $\Gamma_{113}(x_1,x_2)=\sum_{n_1=0}^\infty \frac{x_1^{n_1}}{n_1!}=\text{exp}(x_1)=:q_v$.
                \\

                Ad $\Gamma_{223}$:
                
                Proceeding in the same manner, we find that $\Gamma_{223}(x_1,x_2)=\sum_{n_2=0}^\infty \frac{x_2^{n_2}}{n_2!}=\text{exp}(x_2)=:q_h$.
                \\

                Ad $\Gamma_{123}$:

                This time, the sum evaluates to $I_{(1,0)}(T_1^{\bullet n_1+1}\cdot T_3)+I_{(0,1)}( T_2^{\bullet n_2+1} \cdot T_3)$, which is zero as we've seen above. Thus $\Gamma_{123}(x_1,x_2)=0$.
                \\

                Ad $\Gamma_{333}$: 

                Using the dimension constraint (\ref{dimconsponex}),
                we can deduce that only the invariants with $d=2$ and $e=2$ survive in sum (\ref{summsumm}):
                $$ 2^{n_2}I_{(2,0)}(T_1^{\bullet n_1} \cdot T_3^{\bullet3})\ + \ 2^{n_1}I_{(0,2)}( T_2^{\bullet n_2} \cdot T_3^{\bullet3})\ +\ I_{(1,1)}(T_3^{\bullet3}).$$
                Since a bidegree $(0,2)$ or $(2,0)$ curve is determined by two points but the corresponding invariants contain three point classes, the first two summands vanish. By Proposition \ref{gwndeponex}, the final summand is equal to $N_{(1,1)}=1$.
                Thus $\Gamma_{333}(x_1,x_2)=\sum_{n_1=0}^\infty \sum_{n_2=0}^\infty\frac{x_1^{n_1}}{n_1!}\frac{x_2^{n_2}}{n_2!}=q_v\cdot q_h$.
                \\
                
                Combining these findings, we have:
                $$T_i*T_j=T_i\cup T_j \ + \sum_{k+f=3}\Gamma_{ijk}(x_1,x_2)T_f = T_i\cup T_j + \Gamma_{ij3}(x_1,x_2)T_0.$$ 

                By associativity of $*$ (Proposition \ref{assponex}), this corresponds to the following multiplication table:
                \begin{center}
                    \renewcommand\arraystretch{1.3}
                    \setlength\doublerulesep{0pt}
                    \begin{tabular}{|c||*{4}{c|}}
                        \hline$*$ & $T_0$ & $T_1$ & $T_2$ & $T_3$ \\
                        \hline\hline
                        $T_0$ & $T_0$ & $T_1$ & $T_2$ & $T_3$ \\ 
                        \hline
                        $T_1$ & $T_1$ & $q_v$ & $T_3$ & $q_v T_2$ \\ 
                        \hline
                        $T_2$ & $T_2$ & $T_3$ & $q_h$ & $q_h T_1$ \\ 
                        \hline
                        $T_3$ & $T_3$ & $q_v T_2$ & $q_h T_1$ & $q_v q_h$ \\ 
                        \hline
                    \end{tabular}
                \end{center}

                Rewriting $T_0$ as $1$, $T_1$ as $v$, $T_2$ as $h$ and  $T_3$ as $vh$ concludes the proof.

            \end{proof}

\chapter{Kontsevich's Formula via Associativity of the Quantum Product} %%as application of theory 2
In this chapter, we shall witness the power of the quantum formalism by using it for a more concise proof of Kontsevich's formula. We will do so by deducing a differential equation of the quantum potentials from the associativity of the quantum product and then translating this differential equation into a recursion for the GW-invariants. 
        
    \section{The $\ptwo$ Case}\label{konassp2}
        %deduce from assoc. of *

        % \begin{proposition} see theory 2
        %     %associativity of * (for ptwo )
        % \end{proposition}

        % \begin{proof}%similat to classical nd proof
        % \end{proof}

        \begin{proof}[Alternate Proof of Theorem \ref{nd}]
            Recall that by definition of the (big) quantum product, we have the following equalities:
            \begin{alignat*}{3}
                h^1*h^1 & = h^2 \; + \; & \Gamma_{111}h^1 + \Gamma_{112}h^0&\\
                h^1*h^2 & = & \Gamma_{121}h^1+ \Gamma_{122}h^0&\\
                h^2*h^2 & = & \Gamma_{221}h^1+ \Gamma_{222}h^0&
            \end{alignat*}
            Let us use these in order to compute the products $(h^1*h^1)*h^2$ and $h^1*(h^1*h^2)$ in terms of the generators $h^i$:
            \begin{align*}
                (h^1*h^1)*h^2 = \; & (h^2 + \Gamma_{111}h^1 + \Gamma_{112}h^0)*h^2\\ 
                = \; & h^2*h^2 + \Gamma_{111}h^1*h^2 + \Gamma_{112}h^2\\
                = \; & \Gamma_{221}h^1+ \Gamma_{222}h^0 + \Gamma_{111}(\Gamma_{121}h^1+ \Gamma_{122}h^0) + \Gamma_{112}h^2 \\
                = \; &  (\Gamma_{222}+\Gamma_{111}\Gamma_{122})h^0 + (\Gamma_{122}+\Gamma_{111}\Gamma_{112})h^1+\Gamma_{112} h^2,\\
                 \\
                h^1*(h^1*h^2) = \; & h^1*(\Gamma_{121}h^1+ \Gamma_{122}h^0) \\
                = \; & (\Gamma_{112}\Gamma_{112})h^0 + (\Gamma_{122}+\Gamma_{111}\Gamma_{112})h^1+\Gamma_{112} h^2.
                \end{align*}
                By the associativity of the quantum product (Proposition \ref{assponex}), both products are equal. In particular, the coefficients of the $h^0$-terms coincide, i.e.:
                \begin{equation}\label{h0eq}
                    \Gamma_{222}+\Gamma_{111}\Gamma_{122} = \Gamma_{112}\Gamma_{112}.
                \end{equation}
            Translating this differential equation back into a recursion for the involved GW-invariants will yield Kontsevich's formula for $\ptwo$. Fortunately, we can make several simplifications without sacrificing generality. 
            Recall that  
            $\Gamma_{ijk}(\mathbf{x})=\sum_{\mathbf{a}}\frac{\textbf{x}^{\textbf{a}}}{\textbf{a}!} I_+(\textbf{h}^{\textbf{a}}\cdot h^i \cdot h^j \cdot h^k).$
            By Lemma \ref{fundclassinv}, we have $I_+=0$ if $\textbf{h}^{\textbf{a}}$ contains a class $h^0$. Thus only the case $a_0=0$ contributes, i.e. $\Gamma_{ijk}$ is independent of $x_0$. 
            Furthermore, we know from Lemma \ref{divisoreq} that the $I_+$ containing classes $h^1$ are determined entirely by those without. We may therefore assume w.l.o.g. that $x_0=x_1=1$. After rewriting $x:=x_2$, we get:
            $$\Gamma_{ijk}(x)=\sum_{n=0}^\infty\frac{x^n}{n!} I_+((h^2)^{\bullet n}\cdot h^i\cdot h^j \cdot h^k),$$
            the generating function of the invariants
            $I_+((h^2)^{\bullet n}\cdot h^i\cdot h^j \cdot h^k).$
    
            With this in mind, let us unpack equation (\ref{h0eq}) using the product rule for exponential power series (Definition \ref{prodrule}). More precisely, we compare the coefficients on both sides at some index $n$: 
            \begin{align*}
                I_+((h^2)^{\bullet n}\cdot h^2\cdot h^2 \cdot h^2) & + \sum_{n_A+n_B=n}\binom{n}{n_A}I_+((h^2)^{\bullet n_A}\cdot h^1\cdot h^1 \cdot h^1)\cdot I_+((h^2)^{\bullet n_B}\cdot h^1\cdot h^2 \cdot h^2) \\
                & = \sum_{n_A+n_B=n}\binom{n}{n_A}I_+((h^2)^{\bullet n_A}\cdot h^1\cdot h^1 \cdot h^2)\cdot I_+((h^2)^{\bullet n_B}\cdot h^1\cdot h^1 \cdot h^2).
            \end{align*}
    
            Recall that by definition, $I_+((h^2)^{\bullet n}\cdot h^i\cdot h^j \cdot h^k)=\sum_{d>0}I_{d}((h^2)^{\bullet n}\cdot h^i\cdot h^j \cdot h^k)$.
            We now invoke dimension constraint (\ref{dimconspr}) to examine which invariants are non-trivial. 
            Here, the constraint reads:
            $$2n_*+i+j+k=\text{dim}(\overline{M}_{0,n_*+3}(\ptwo,d_*))=3d_*+2+n_*,$$
            where the subscript $*$ can be omitted or replaced by $A$ or $B$.
            Thus for fixed $d_*$, only the following case survives:
            \begin{equation}\label{dimcons}
                n_*=3d_*+2-i-j-k.
            \end{equation}
            
            Writing out the $I+$-terms leads to:
            \begin{align*}
                & \sum_{d>0}I_{d}((h^2)^{\bullet n}\cdot h^2\cdot h^2 \cdot h^2) \\
                & + \sum_{n_A+n_B=n}\binom{n}{n_A}\sum_{d_A>0}I_{d_A}((h^2)^{\bullet n_A}\cdot h^1\cdot h^1 \cdot h^1)\cdot \sum_{d_B>0}I_{d_B}((h^2)^{\bullet n_B}\cdot h^1\cdot h^2 \cdot h^2) \\
                & = \sum_{n_A+n_B=n}\binom{n}{n_A}\sum_{d_A>0}I_{d_A}((h^2)^{\bullet n_A}\cdot h^1\cdot h^1 \cdot h^2)\cdot \sum_{d_B>0}I_{d_B}((h^2)^{\bullet n_B}\cdot h^1\cdot h^1 \cdot h^2).
            \end{align*}
            We then apply Lemma \ref{divisoreq} in order to extract the $h^1$-classes:
            \begin{align*}
                \sum_{d>0}I_{d}((h^2)^{\bullet n+3}) 
                &\; + \sum_{n_A+n_B=n}\binom{n}{n_A}\sum_{d_A>0}d_A^3I_{d_A}((h^2)^{\bullet n_A})\cdot \sum_{d_B>0}d_BI_{d_B}((h^2)^{\bullet n_B+2}) \\
                &\; = \sum_{n_A+n_B=n}\binom{n}{n_A}\sum_{d_A>0}d_A^2I_{d_A}((h^2)^{\bullet n_A+1})\cdot \sum_{d_B>0}d_B^2I_{d_B}((h^2)^{\bullet n_B+1}).
            \end{align*}
            Now we apply the dimension constraint (\ref{dimcons}) from above in order to express $n_*$ in terms of $d_*$. It follows from the constraint that $n_A+n_B=n$ if and only if $d_A+d_B=d$. Fixing the only $d$ contributing to $I_+$ (recall Proposition-Definition \ref{collinvpr}) allows us to rewrite the equation again:
            \begin{align*}
                I_{d}((h^2)^{\bullet 3d-1}) 
                \;& + \sum_{d_A+d_B=d}\binom{3d-4}{3d_A-1}d_A^3I_{d_A}((h^2)^{\bullet 3d_A-1})\cdot d_BI_{d_B}((h^2)^{\bullet 3d_B-1}) \\
                & = \sum_{n_A+n_B=n}\binom{2d-4}{3d_A-2}d_A^2I_{d_A}((h^2)^{\bullet 3d_A-1})\cdot d_B^2I_{d_B}((h^2)^{\bullet 3d_B-1}).
            \end{align*}
            But by Corollary \ref{gwndpr}, this is precisely Kontesvich's formula for $\ptwo$, i.e. equation (\ref{kontsformp2}).

        \end{proof}

    \section{The $\ponex$ Case}\label{konassponex} 
        %deduce from assoc. of *

        % \begin{proposition} see theory 2
        %     %associativity of * (for ponex)
        % \end{proposition}

        % \begin{proof}
        % \end{proof}
        
        \begin{proof}[Alternate Proof of Theorem \ref{nde}] %exc 5.8 WRITE THIS ONE FIRST
            Recall that by definition of the (big) quantum product, we have the following equalities:
            \begin{alignat*}{3}
                T_1*T_2 & = T_3 \; + \; & \Gamma_{121}T_2 + \Gamma_{122}T_1+ \Gamma_{123}T_0&\\
                T_1*T_3 & = & \Gamma_{131}T_2 + \Gamma_{132}T_1+ \Gamma_{133}T_0&\\
                T_2*T_3 & = & \Gamma_{231}T_2 + \Gamma_{232}T_1+ \Gamma_{233}T_0&\\
                T_3*T_3 & = & \Gamma_{331}T_2 + \Gamma_{332}T_1+ \Gamma_{333}T_0&.
            \end{alignat*}
            Recall that in the original proof of Kontsevich's formula, our data was comprised of a vertical rule, a horizontal rule and several points. We choose the corresponding cohomology classes, namely $T_1$, $T_2$ and $T3$.
            In terms of the generators $T_i$, the products $(T_1*T_2)*T_3$ and $T_1*(T_2*T_3)$ look as follows:
            \begin{align*}
                (T_1*T_2)*T_3 = \; & (T_3 + \Gamma_{121}T_2 + \Gamma_{122}T_1+ \Gamma_{123}T_0)*T_3 \\
                 = \; & T_3*T_3 + \Gamma_{121}T_2*T_3 + \Gamma_{122}T_1*T_3 + \Gamma_{123}T_0*T_3 \\
                 = \; & (\Gamma_{331}T_2 + \Gamma_{332}T_1+ \Gamma_{333}T_0) + \Gamma_{121}(\Gamma_{231}T_2 + \Gamma_{232}T_1+ \Gamma_{233}T_0) \\ 
                & + \Gamma_{122}(\Gamma_{131}T_2 + \Gamma_{132}T_1+ \Gamma_{133}T_0) + \Gamma_{123}T_3 \\
                 = \; & (\Gamma_{333}+\Gamma_{112}\Gamma_{233}+\Gamma_{122}\Gamma_{133})T_0 + (\Gamma_{233}+\Gamma_{112}\Gamma_{223}+\Gamma_{122}\Gamma_{123})T_1 \\
                 & + (\Gamma_{133}+\Gamma_{112}\Gamma_{123}+\Gamma_{122}\Gamma_{113})T_2 + \Gamma_{123}T_3
                ,
            \end{align*}
            \begin{align*}
               T_1*(T_2*T_3) = \; & T_1 * (\Gamma_{231}T_2 + \Gamma_{232}T_1+ \Gamma_{233}T_0) \\
                 = \; & \Gamma_{123}(T_3 + \Gamma_{121}T_2+ \Gamma_{122}T_1 + \Gamma_{123}T_0) + \Gamma_{223}(\Gamma_{111}T_2 + \Gamma_{112}T_1+ \Gamma_{113}T_0) +\Gamma_{233}T_1 \\ 
                 = \; & (\Gamma_{123}\Gamma_{123}+\Gamma_{223}\Gamma_{113})T_0 + (\Gamma_{123}\Gamma_{122}+\Gamma_{223}\Gamma_{112} + \Gamma_{233})T_1 \\
                 & + (\Gamma_{123}\Gamma_{112}+\Gamma_{223}\Gamma_{111})T_2 + \Gamma_{123}T_3
                 .
            \end{align*}
            By the associativity of the quantum product (Proposition \ref{assponex}), both products are equal. In particular, the coefficients of the $T_0$-terms coincide:%- and $T_2$, i.e. \textcolor{red}{remove unnecessary one}:
            \begin{equation}\label{t0eq}
                \Gamma_{333}+\Gamma_{112}\Gamma_{233}+\Gamma_{122}\Gamma_{133} = \Gamma_{123}\Gamma_{123}+\Gamma_{223}\Gamma_{113}.
            \end{equation}
            % \begin{equation}\label{t2eq}
            %     \Gamma_{133}+\Gamma_{122}\Gamma_{123}+\Gamma_{122}\Gamma_{113} = \Gamma_{123}\Gamma_{112}+\Gamma_{223}\Gamma_{111}.
            % \end{equation}
            We translate the differential equation (\ref{t0eq}) into a recursion for the involved GW-invariants. Again, we can make several simplifications without sacrificing generality. 
            Recall that $\Gamma(\mathbf{x})=\sum_{\mathbf{a}}\frac{\textbf{x}^{\textbf{a}}}{\textbf{a}!} I_+(\textbf{T}^{\textbf{a}})$ and $\Gamma_{ijk}(\mathbf{x})=\sum_{\mathbf{a}}\frac{\textbf{x}^{\textbf{a}}}{\textbf{a}!} I_+(\textbf{T}^{\textbf{a}}\cdot T_i \cdot T_j \cdot T_k)$, where $\textbf{T}^{\textbf{a}}=T_0^{\bullet a_0}\cdot T_1^{\bullet a_1}\cdot T_2^{\bullet a_2}\cdot T_3^{\bullet a_3}$.
            By Lemma \ref{fundclassinvponex}, we have $I_+=0$ if $\textbf{T}^{\textbf{a}}$ contains a class $T_0$. Thus only the case $a_0=0$ contributes, i.e. $\Gamma_{ijk}$ is independent of $x_0$. 
            Note that this time, we cannot omit the variables corresponding to the divisor classes: since both degrees $d$ and $e$ could hypothetically be pulled out of an invariant using the divisor equation, there no longer is a one to one correspondence between those invariants with divisor classes and those without. The potentials of interest therefore look as follows:
            $$\Gamma(x_1,x_2,x_3)=\sum_{n_1,n_2,n_3\ge0}\frac{x_1^{n_1}x_2^{n_2}x_3^{n_3}}{n_1!n_2!n_3!} I_+(T_1^{\bullet n_1}\cdot T_2^{\bullet n_2}\cdot T_3^{\bullet n_3}).$$
            % the generating function of the invariants
            % $I_+(T_3^{\bullet n}\cdot T_i\cdot T_j\cdot T_k).$

            Using the definition of $I_+$, the divisor equation and arguments similar to those used in the discussion of $\Gamma_{333}$ in the proof of Proposition \ref{qringponex}, and Proposition \ref{gwndeponex}, we rewrite this as follows:
            \begin{align*}
                \Gamma(x_1,x_2,x_3)&=\sum_{d+e>0}\;\sum_{n_1\ge0}e^{n_1}\frac{x_1^{n_1}}{n_1!}\sum_{n_2\ge0}e^{n_2}\frac{x_2^{n_2}}{n_2!}\sum_{n_3\ge0}I_{(d,e)}(T_3^{\bullet n_3})\frac{x_3^{n_3}}{n_3!} \\
                &=\sum_{d+e>0}N_{(d,e)}\frac{x_3^{2(d+e)-1}}{(2(d+e)-1)!}\exp(ex_1+dx_2).
            \end{align*}
            Note that differentiating $\Gamma$ in $x_1$ simply produces an extra factor $e$ in the sum over $d+e\ge0$ (and eliminates the $e=0$ term). In the same way, the derivative in $x_2$ contains an extra factor $d$. 
            Plugging this information into equation (\ref{t0eq}) and comparing coefficients at some pair $(d,e)$ yields Kontsevich's formula. 

            To visualize this, we treat the terms $\Gamma_{333}$ and $\Gamma_{112}\Gamma_{233}$ in detail.
            By differentiating our new characterization of $\Gamma$ three times in $x_3$, we obtain the following:
            $$\Gamma_{333}=\sum_{d+e>0}N_{(d,e)}\frac{x_3^{2(d+e)-4}}{(2(d+e)-4)!}\exp(ex_1+dx_2).$$
            Working out the product of 
            $$\Gamma_{112}=\sum_{d+e>0}d^2eN_{(d,e)}\frac{x_3^{2(d+e)-1}}{(2(d+e)-1)!}\exp(ex_1+dx_2)$$
            and
            $$\Gamma_{233}=\sum_{d+e>0}eN_{(d,e)}\frac{x_3^{2(d+e)-3}}{(2(d+e)-3)!}\exp(ex_1+dx_2) $$
            yields the following rather large expression:
            $$\sum_{d+e>0}\bigl(\sum_{\substack{d_A+d_B=d\\e_A+e_B=e}}\frac{1}{(2(d_A+e_A)-1)!}\frac{1}{(2(d_B+e_B)-3)!}d_Ae_A^2N_{(d_A,e_A)}d_BN_{(d_B,e_B)}\bigr)x_3^{2(d+e)-4}\exp(ex_1+dx_2).$$
            Multiplying both results by $(2(d+e)-4)!$ and comparing the coefficients at $$x_3^{2(d+e)-4}\exp(ex_1+dx_2)$$ turns them into the familiar terms from Kontsevich's formula. Indeed, the coefficient from $\Gamma_{333}$ is simply $N_{(d,e)}$, and the coefficient from $\Gamma_{112}\Gamma_{233}$ looks as follows:
            $$\sum_{\substack{d_A+d_B=d\\e_A+e_B=e}}\binom{2(d+e)-4}{2(d_A+e_A)-1}d_Ae_A^2N_{(d_A,e_A)}d_BN_{(d_B,e_B)}.$$
            Repeating these arguments for the remaining terms of our differential equation (\ref{t0eq}) yields Kontsevich's formula for $\ponex$, i.e. (\ref{kontsponex}).
            
        \end{proof}

        \begin{remark}[Equivalence]
            Note that in the 2-dimensional spaces $\ptwo$ and $\ponex$, the only non-divisor classes are the point class and the fundamental class. By Lemmas \ref{fundclassinv} and \ref{fundclassinvponex}, the GW-invariants containing the fundamental class are determined entirely by those without. The same applies to codimension $1$ classes by Lemmas \ref{divisoreq} and \ref{divisoreqponex}, so one is left with invariants containing only point classes. But these are determined entirely by the dimension constraints \ref{dimcons} and \ref{dimconsponex} and by Kontsevich's formula due to Lemmas \ref{gwndpr} and \ref{gwndeponex}. We conclude that there is no loss of generality when passing from the big to the small quantum cohomology for these spaces, since both are determined entirely by Kontsevich's formula.
            
            But then the same holds for the quantum product, which is defined using the GW-ivariants, and thereby also for its associativity.

            Combining this with the above proofs, one could say that Kontsevich's formula is in fact equivalent to the associativity of the quantum product.
        \end{remark} %as application of theory 2
%what else can be done with theory 2

\chapter{Perspectives}\label{aus}

%stuff from generalization chapters in book

%what else can be done with quantum cohomology. 
%maybe frobienius manifolds from 5.5.4 as connection to DG

%wiki for qh: Because it expresses a structure or pattern for Gromov–Witten invariants, quantum cohomology has important implications for enumerative geometry. It also connects to many ideas in mathematical physics and mirror symmetry. In particular, it is ring-isomorphic to symplectic Floer homology.

%higher genus, classes instead of degrees
%other targets, e.g. more general del Pezzo surfaces
%See \cite{Kontsevich_1994} for more general formula (with $\beta$ etc.).

%products of other proj spaces, can knowledge of pr be used? or simply the more general forms from rahul etc
Since the possibilities for further study and generalization are too numerous to be listed here, we will briefly mention only a few. 

Given what we already have, studying the convergence behaviour of our invariants $N_d$ and $N_{(d,e)}$ might lead to interesting results.

As a first generalization, one could of course consider curves of positive genus $g>0$. For $g=1$, there is in fact a recursive formula for the numbers $E_d$ of degree $d$ plane curves of genus $1$ passing through $3d+g-1=3d$ general points, see \cite{panda97}.

One could also consider new target spaces, such as more general del Pezzo surfaces or a general smooth projective variety $X$. It would then become necessary to generalize our notion of a degree to classes $\beta\in A_1(X)$, see for example \cite{Kontsevich_1994}. But even in this case, there exists a projective coarse moduli space parametrizing maps $\mu:C\rightarrow X$ with $\mu_*[C]=\beta$. 

Looking to different horizons, one could study how quantum cohomology relates to mathematical physics, in particular to string theory and mirror symmetry. As previously mentioned, a seminal reference is \cite{Witten:1990hr}. Another interesting context is that of Riemannian geometry: here one can use the structure constants $\Phi_{ijk}$ as Christoffel symbols in order to define a formal connection, see \cite{dub}. The flatness of this connection turns out to be equivalent to the associativity of the quantum product. In this way, one constructs a class of examples of so called Frobenius manifolds.

There is also the possibility of transferring the entire theory into the relatively new language of stacks, whereby automorphisms can be avoided by bypassing the definition of the moduli functor.
%convergence rates of invariants %what else can be done with new theory
\chapter{Acknowledgements}

I owe a great debt of gratitude to my two advisors, Rahul Pandharipande and Junliang Shen. It is due to Rahul's willingness to lend an ear to an unknown (and nervous) student at his office door and to Junliang's immense openness and hospitality that I was able to embark upon the exhilarating adventure of my semester at Yale. 
I would like to particularly thank Junliang for providing me with guidance and direction throughout my stay, and for sharing with me his valuable insights into mathematical research. 

Another essential figure I am indebted to is Yvette Barnard, the Yale mathematics registrar: I might have gone astray in bureaucracy had it not been for her ever swift and effective assistance.

Finally, I am grateful to the mathematics departments of ETH and Yale at large: the former for its indispensable organizational and financial assistance, and the latter for providing fertile soil for fruitful discussion, each exchange with one of its members teaching me something new. 
It was an honor to temporarily count myself among them and to thereby see my love of mathematics reinvigorated.

\appendix
\chapter{Supplementary Materials}\label{appa}

    %\section{Supplements}
        %maybe N11, N12 and N1d from ponex classical

         \begin{lemma}\label{qijcodim}%lemma 3.5.2
            Let $n\ge2$ and consider the locus of maps whose marks $p_i\neq p_j$ have the same image:
            $$Q_{ij}=\{\mu\in{M}_{0,n}(\pr,d)|\mu(p_i)=\mu(p_j)\}.$$
            Then, writing $M:={M}_{0,n}(\pr,d)$, we have that $\text{codim}_M(Q_{ij})=r$.
        \end{lemma}    
        
        \begin{proof}
             %reduction to n geq 3 from lemma 2.8.3
            We can assume w.l.o.g. that $n\ge3$. Indeed, consider the commutative diagram (\ref{diag}).
            Denote by $\Delta$ the diagonal in $\pr\times\pr$. Then the result for $n+1$ implies the result for $n$ by flatness of $\epsilon$:
            \begin{align*}
                \hat Q_{ij} &= (\hat\nu_i\times\hat\nu_j)^{-1}(\Delta)\\
                &= \epsilon^{-1}(\nu_i\times\nu_j)^{-1}(\Delta)\\
                &= \epsilon^{-1}(Q_{ij}).
            \end{align*}

            We can now write $\m=M_{0,n}\times W(r,d)$. Recall that $W(r,d)$ is the space of $(r+1)$-tuples of degree $d$ binary forms. Let $\sum_{i=0}^da_{ki}x^{d-i}y^i$ be the $k$-th form for $0\le k\le r$. Assume w.l.o.g. that $p_i=[0:1]$ and $p_j=[1:0]$. Then $\mu(p_i)=\mu(p_j)$ if and only if $(a_{0d},a_{1d},\dots,a_{rd})=\lambda((a_{00},a_{10},\dots,a_{r0}))$ for some $\lambda\in\C^*$. Assuming w.l.o.g. that $a_{00}$ and $a_{0d}$ are non-zero, these are $r$ independent conditions in the $a_{ij}$. 
            Thus the fiber of the projection $\m\rightarrow W(r,d)$ has dimension $rd+d$ (recall that $\dim(W(r,d)=rd+d+r$). Adding to this the dimension of $M_{0,n}$, we get that $\dim(Q_{ij})=rd+r+n-3$, which is equivalent to the desired codimension because of equation (\ref{dimmbar}). %see 2.1.2. needs also normal crossings, right?
            
        \end{proof}

        \begin{lemma}\label{355}%3.5.5 book
            Let $\Gamma_1,\dots,\Gamma_n\subset\pr$ be general subvarieties with $\codim \Gamma_i)\geq2$ for all $i$ and such that $\sum\codim\Gamma_i)=\dim\mbar)$. 
            Then for any $\mu\in\underline\nu^{-1}(\underline\Gamma)$, the intersection of $\Gamma_i$ with the image curve $C=\mu(\pone)$ consists only of the point $\mu(p_i)$.
        \end{lemma}
        For the proof, see Lemma 3.5.5 in \cite{invitation}.
        % \begin{proof}
        %     \textcolor{red}{add?}
        % \end{proof}

        \begin{proposition}\label{n3}%move back into chapter?
            There are exactly $N_3=12$ rational cubics passing through 8 points in general position in $\ptwo$.
        \end{proposition}

        %this is preparation for Nd proof. 
        \begin{proof}[Proof] %do this first. base prev pf on this
        
            We now work in $\overline{M}_{0,9}(\ptwo,3)$ and denote the 9 marks by $m_1$, $m_2$, $p_1$, \dots, $p_{7}$. As our data, we pick two lines $L_1$ and $L_2$ and $7$ points $Q_1, \dots, Q_{7}$ in general position in $\ptwo$.
            We denote by $Y$ the subvariety of $\overline{M}_{0,9}(\ptwo,3)$ consisting of stable maps $(C;m_1, m_2, p_1, \dots, p_{7};\mu)$ sending $m_1$ to $L_1$, $m_2$ to $L_2$ and $p_i$ to $Q_i$ for $1\leq i\leq 7$:
            $$Y=\nu_{m_1}^{-1}(L_1) \cap \nu_{m_2}^{-1}(L_2) \cap \nu_{p_1}^{-1}(Q_1) \cap \dots \cap \nu_{p_{7}}^{-1}(Q_{7}).$$
            By Lemmas \ref{lemma1} and \ref{lemma2} and Remark \ref{onY}, $Y$ is a curve and its intersection with each boundary divisor is transversal %in the sense of satisfying Kleiman's eq 
            and takes place in $\overline{M}_{0,9}^*(\ptwo,3)$.

            We again use the fundamental boundary relation (\ref{fundrel}) to conclude that
            \begin{equation}\label{equivn3}
                \abs{Y\cap D(m_1,m_2|p_1,p_2)} = \abs{Y\cap D(m_1,p_1|m_2,p_2)}.
            \end{equation}
            where $D(m_1,m_2|p_1,p_2)$ and $D(m_1,p_1|m_2,p_2)$ represent sums of the form $$\sum_{\substack{m_1,m_2\in A \\ p_1,p_2\in B \\d_A+d_B=3}}D(A,B;d_A,d_B) \; \text{ and} \sum_{\substack{m_1,p_1\in A \\ m_2,p_2\in B \\d_A+d_B=3}}D(A,B;d_A,d_B),$$ respectively.
            As in Proposition \ref{n1}, counting up then equating the contributions of the left and the right hand side will yield the desired result.

            LHS:\\
            Again we begin by illustrating $D_{AB}:=D(A,B;d_A,d_B)$ before the distribution of the five spare marks $p_3$, \dots, $p_7$. 
            
            \definecolor{rvwvcq}{rgb}{0.08235294117647059,0.396078431372549,0.7529411764705882}
            \definecolor{wrwrwr}{rgb}{0.3803921568627451,0.3803921568627451,0.3803921568627451}
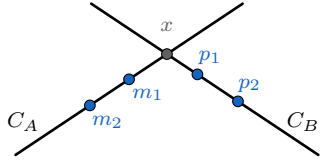
            \begin{tikzpicture}[line cap=round,line join=round,>=triangle 45,x=1cm,y=1cm]
                \clip(-5,0.5) rectangle (4.5,4); %size of rectangle containing image (given by 2pts)
                \draw [line width=1pt] (1,3)-- (4,1);
                \draw [line width=1pt] (3,3)-- (0,1);
                \begin{scriptsize}
                    \draw[color=black] (3.8,1.44) node {$C_{B}$};
                    
                    \draw [fill=rvwvcq] (2.4025391414130275,2.0649739057246483) circle (2pt);
                    \draw[color=rvwvcq] (2.57,2.2890429805862818) node {$p_{1}$};
                    
                    \draw [fill=rvwvcq] (2.938269319967983,1.7078204533546781) circle (2pt);
                    \draw[color=rvwvcq] (3.092541216255241,1.9233276303855245) node {$p_{2}$};
                    
                    \draw[color=black] (0.1,1.44) node {$C_{A}$};
                    
                    \draw [fill=rvwvcq] (1.4993503144908957,1.9995668763272638) circle (2pt);
                    \draw[color=rvwvcq] (1.75,1.8) node {$m_{1}$};
                    
                    \draw [fill=rvwvcq] (0.9819670598590557,1.654644706572704) circle (2pt);
                    \draw[color=rvwvcq] (1.2162624630513585,1.44) node {$m_{2}$};
                    
                    \draw [fill=wrwrwr] (2,2.3333333333333335) circle (2pt);
                    \draw[color=wrwrwr] (2,2.7) node {$x$};
                \end{scriptsize}
            \end{tikzpicture}

            As before, $x$ denotes the intersection node of $C_A$ and $C_B$. The picture is the same as in the proof of Proposition \ref{n2}, however there are more degrees to consider and more spare marks.
            
            Since we must have $d_A+d_B=3$, the only four possible pairs of degrees are $(d_A,d_B)\in\{(0,3),(1,2),(2,1),(3,0)\}$. The five marks $p_3$, \dots, $p_7$ remain to be distributed between $C_A$ and $C_B$. There are $2\cdot(\binom{5}{0}+\binom{5}{1}+\binom{5}{2})=32$ such distributions. The number of irreducible components of the divisor $Y\cap D(m_1,m_2|p_1,p_2)$ is therefore equal to $4\cdot32=128$, i.e. there are 128 summands of the form $D_{AB}$ on the left hand side. 
            
            We now study the intersection of $Y$ with each of those divisors $D_{AB}$ by considering each pair of degrees $(d_A,d_B)$.\\

            Ad $(3,0)$:
            This pair does not contribute to the intersection. Indeed, $Q_1=\mu(p_1)=\mu(p_2)=Q_2$ contradicts the genericity of our data. \\

            Ad $(0,3)$: %could put a tilde on gamma1, cf lemmas1,2 in notes. but we've probably explained well enoug in lemmma 1 why we can assume Gamma without tilde to have the desired properties
            In this case we get $\{\mu(m_1)\}=\{\mu(m_2)\}=L_1\cap L_2$. Genericity can only be guaranteed if all spare marks lie on $C_B$. We can apply Lemmas \ref{lemma1} and \ref{lemma2} with ${\Gamma}_1=\mu(m_1)$, $\Gamma_2=Q_1$, \dots, $\Gamma_8=Q_7$ to deduce a contribution of $N_3$, the number we are looking for. \\

            Ad $(1,2)$: 
            Let us examine first the possible distributions of the spare $p_i$'s, then the positions of the $pi$'s and $m_i$'s, and finally the position of the intersection node $x$.
            \begin{itemize}
                %number and distribution of marks
                \item We must add 2 marks to $C_A$ and 3 marks to $C_B$. Indeed, any other distribution would contradict genericity: putting more than 2 marks on $C_A$ (or equivalently, less than 3 on $C_B$) would force more than 2 of the $Q_i$'s onto one line. Similarly, putting more than 3 marks on $C_B$ (or equivalently, less than 2 on $C_B$) would force more than 5 of the $Q_i$'s onto one conic.
                There are thus $\binom{5}{2}=10$ possible distributions of the marks onto $C_A$ and $C_B$.
                %positions of marks
                \item Once distributed, the positions of the points $p_1$, \dots, $p_7$ on their respective twigs are uniquely determined. Indeed, for each $i\in\{1, \dots,7\}$, the point $p_i$ is the only element of $\mu^{-1}(Q_i)$ by Lemma \ref{lemma2}, and the order of the $p_i$'s is fixed because the maps $C_A \rightarrow \ptwo$ and $C_B \rightarrow \ptwo$ are each birational onto their respective image. 
                In particular, setting $\Gamma_i=p_i$ for the two marks $p_i$ of $C_A$ and applying Lemmas \ref{lemma1} and \ref{lemma2} yields a contribution of $N_1$. The analogous procedure on the five marks of $C_B$ yields a contribution of $N_2$. 

                The position of the mark $m_1$ is determined by the intersection $\mu(C_A)\cap L_1$. By B\'ezout's Theorem, there is thus precisely $d_A\cdot1=1$ possible position for $m_1$. The same argument applies to $m_2$, where $\mu(C_A)\cap L_2$ yields another contribution of $d_A$, so the total contribution from the positions of $m_1$ and $m_2$ is equal to $d_A^2=1$.
                %gluing
                \item Finally, there are two ways of 'gluing' $\mu(C_A)$ and $\mu(C_B)$ together: again by B\'ezout's Theorem, we have $d_A\cdot d_B=1\cdot2=2$ possible choices for $\mu(x)$.
            \end{itemize} 

            The total contribution from this pair of degrees is therefore equal to $\binom{5}{2}\cdot d_A^2 N_1\cdot d_A d_B \cdot N_2=20$.
            
            Ad $(2,1)$: We proceed as for $(1,2)$.
            \begin{itemize}
                %number and distribution of marks
                \item By the same reasoning as above, we must add 5 marks to $C_A$ and no marks to $C_B$.
                There is thus $\binom{5}{5}=1$ possible distribution of the marks onto $C_A$ and $C_B$.
                %positions of marks
                \item Once distributed, the positions and orders of the points $p_1$, \dots, $p_7$ on their respective twigs are uniquely determined by Lemma \ref{lemma2} and the birationality of the maps $C_A \rightarrow \ptwo$ and $C_B \rightarrow \ptwo$. 
                In particular, we get a contribution of $N_2$ from $C_A$ and a contribution of $N_1$ from $C_B$. 

                The positions of $m_1$ and $m_2$ are determined by the intersections $\mu(C_A)\cap L_1$ and $\mu(C_A)\cap L_2$, respectively. As before, the number of possible positions can be computed with B\'ezout's Theorem. In particular, each intersection contributes a factor $d_A\cdot1=d_A$ for a total contribution of $d_A^2=4$.

                \item As before, there are $d_A\cdot d_B=2\cdot1=2$ ways of 'gluing' $C_A$ and $C_B$ together by B\'ezout's Theorem. 
            \end{itemize} 
            This pair therefore contributes a factor $\binom{5}{5}\cdot d_A^2 N_1\cdot d_A d_B \cdot N_2=8$.
            \\
            \\
            As discussed above, these are all possible pairs of degrees. The left hand side of equality (\ref{equivn3}) thus yields the number $N_3+28$.\\

            RHS:\\  
            Before the distribution of the spare marks, the $D_{AB}$ for the right hand side look the same as in the proof of Proposition \ref{n2}.
            
            \definecolor{rvwvcq}{rgb}{0.08235294117647059,0.396078431372549,0.7529411764705882}
            \definecolor{wrwrwr}{rgb}{0.3803921568627451,0.3803921568627451,0.3803921568627451}
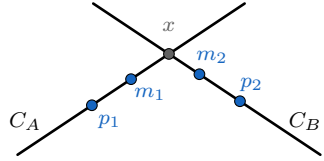
            \begin{tikzpicture}[line cap=round,line join=round,>=triangle 45,x=1cm,y=1cm]
                \clip(-5,0.5) rectangle (4.5,4); %size of rectangle containing image (given by 2pts)
                \draw [line width=1pt] (1,3)-- (4,1);
                \draw [line width=1pt] (3,3)-- (0,1);
                \begin{scriptsize}
                    \draw[color=black] (3.8,1.44) node {$C_{B}$};
                    
                    \draw [fill=rvwvcq] (2.4025391414130275,2.0649739057246483) circle (2pt);
                    \draw[color=rvwvcq] (2.57,2.2890429805862818) node {$m_2$};
                    
                    \draw [fill=rvwvcq] (2.938269319967983,1.7078204533546781) circle (2pt);
                    \draw[color=rvwvcq] (3.092541216255241,1.9233276303855245) node {$p_{2}$};
                    
                    \draw[color=black] (0.1,1.44) node {$C_{A}$};
                    
                    \draw [fill=rvwvcq] (1.4993503144908957,1.9995668763272638) circle (2pt);
                    \draw[color=rvwvcq] (1.75,1.8) node {$m_{1}$};
                    
                    \draw [fill=rvwvcq] (0.9819670598590557,1.654644706572704) circle (2pt);
                    \draw[color=rvwvcq] (1.2162624630513585,1.44) node {$p_1$};
                    
                    \draw [fill=wrwrwr] (2,2.3333333333333335) circle (2pt);
                    \draw[color=wrwrwr] (2,2.7) node {$x$};
                \end{scriptsize}
            \end{tikzpicture}
            
        The possible pairs of degrees $(d_A,d_B)$ are the same as before, and we must again distribute the five marks $p_3, \dots,p_7$. Again we proceed to an examination of the divisors $D_{AB}$ for each pair of degrees.

        Ad $(0,3)$, $(3,0)$: If $(d_A,d_B)=(0,3)$, then $Q_1=\mu(p_1)=\mu(m_1)\in L_1$. Similarly, if $(d_A,d_B)=(3,0)$, then $Q_2\in L_2$. Both conclusions contradict the genericity of the data, so neither case contributes to the right hand side.\\

        Ad $(1,2)$: We proceed as for the left hand side. 
        \begin{itemize}
            %number and distribution of marks
            \item By the same reasoning as before, we need one extra $p_i$ on $C_A$ and four more $p_i$-marks on $C_B$.
            There are thus $\binom{5}{1}=5$ possible distributions of the marks onto $C_A$ and $C_B$.
            %positions of marks
            \item Again, the positions and orders of the points $p_1$, \dots, $p_7$ are uniquely determined. 
            We get a contribution of $N_1$ from $C_A$ and a contribution of $N_2$ from $C_B$. 

            The positions of $m_1$ and $m_2$ are determined by the intersections $\mu(C_A)\cap L_1$ and $\mu(C_B)\cap L_2$, respectively. By B\'ezout's Theorem, we get respective contributions of $d_A\cdot1=d_A$ and $d_B\cdot1=d_B$. In total, there are thus $d_A\cdot d_B=2$ possible positions for $m_1$ and $m_2$. 

            \item As in the previous examinations, there are $d_A\cdot d_B=2\cdot1=2$ ways of 'gluing' $C_A$ and $C_B$ together. 
        \end{itemize} 
        This pair therefore contributes a factor $\binom{5}{1}\cdot d_A\cdot N_1\cdot d_A d_B \cdot d_B N_2=20$.
        \\

        Ad $(2,1)$: Due to the symmetry of the $D_{AB}$ on the right hand side, this case is identical to the previous one up to renaming the points and twigs. We get a contribution of $\binom{5}{1}\cdot d_A N_2\cdot d_A d_B \cdot d_B N_1=20$.
        \\
        \\
        In conclusion, the right hand side thus yields the number $40$.
        Combining the results of both counts, equality (\ref{equivn3}) implies $N_3+28=40$, or $N_3=12$, as claimed.
        
        \end{proof}
\chapter{Supplementary Theory}\label{appb}

    %see if it might be better to just add reference in chapter.

    %genus formula for nodal plane curves for p2 and for ponex
    %https://encyclopediaofmath.org/wiki/Genus_of_a_curve#:~:text=Curves%20of%20genus%20g%3D1,non%2Dhyper%2Delliptic%20curves.
    \begin{theorem}\label{genus} %see pg 9 of n3 notes. shafarevic 3.3
        The genus of a nodal curve $C$ of degree $d$ in $\ptwo$ is given by $g_C=\frac{(d-1)(d-2)}{2}-\delta$, where $\delta$ is the number of nodes of $C$.
    \end{theorem}
    This is the genus-degree formula, a classical result.

    \begin{theorem}\label{genusponex} %see pg 9 of n3 notes. shafarevic 2.3
        The genus of a smooth curve $C\subset\ponex$ of bidegree $(d,e)$ is given by
        $g_C=(d-1)(e-1)$.
    \end{theorem}
    This is example 4.11 in section 2.2.3 of \cite{shafarevich}, where it is shown using the adjunction formula.

    The following theorem is taken from \cite{Fulton} (Appendix B. 9.2) and \cite{invitation} (section 3.4). A proof can be found in \cite{Kleiman}.
    
    Let $G$ be a connected algebraic group, $A$ an irreducible variety with a transitive $G$-action and $B$ and $C$ two more irreducible varieties. Let $f:B\rightarrow A$ and $g:C\rightarrow A$ be morphisms of irreducible varieties:
    \begin{center}
        \begin{tikzcd}
        && C \arrow[dd,"g"]\\\\
        B\arrow[rr,"f"] && A
        \end{tikzcd}
    \end{center}
    
    Given $\sigma\in G$, we denote by $B^{\sigma}$ the variety $B$ considered as a variety over $A$ via $\sigma \circ f$.
    
    \begin{theorem}[Kleiman]\label{kleiman}
        There exists a dense open subset $U\subset G$ such that for every $\sigma\in U$, the fiber product $B^{\sigma}\times_A C$ is either empty or satisfies $$\text{dim}(B^{\sigma}\times_A C)=\text{dim} (B) + \text{dim} (C) - \text{dim} (A).$$
        If $B$ and $C$ are smooth, then $U$ can be chosen such that $B^{\sigma}\times_A C$ is smooth for every $\sigma\in U$.
    \end{theorem}

\bibliographystyle{alpha}
\bibliography{ref}

\begin{thebibliography}{Dub94}

\bibitem[Dub94]{dub}
Boris Dubrovin.
\newblock Geometry of 2d topological field theories, 1994.

\bibitem[EH16]{eh2}
David Eisenbud and Joe Harris.
\newblock {\em 3264 and All That: A Second Course in Algebraic Geometry}.
\newblock Cambridge University Press, 2016.

\bibitem[FM11]{greg}
Greg Friedman and James McClure.
\newblock Cup and cap products in intersection (co)homology.
\newblock 2011.

\bibitem[FP96]{FPnotes}
W.~Fulton and R.~Pandharipande.
\newblock Notes on stable maps and quantum cohomology, 1996.

\bibitem[Ful75]{fulton2}
William Fulton.
\newblock Rational equivalence on singular varieties.
\newblock {\em Publications Math\'ematiques de l'IH\'ES}, 45:147--167, 1975.

\bibitem[Ful98]{Fulton}
W.~Fulton.
\newblock {\em Intersection Theory}.
\newblock Ergebnisse der Mathematik und ihrer Grenzgebiete. Springer New York,
  1998.

\bibitem[Itz94]{Itzy}
C.~Itzykson.
\newblock {Counting rational curves on rational surfaces}.
\newblock {\em Int. J. Mod. Phys. B}, 8:3703--3724, 1994.

\bibitem[Kle74]{Kleiman}
Steven~L. Kleiman.
\newblock The transversality of a general translate.
\newblock {\em Compositio Mathematica}, 28(3):287--297, 1974.

\bibitem[KM94]{Kontsevich_1994}
M.~Kontsevich and Yu. Manin.
\newblock Gromov-witten classes, quantum cohomology, and enumerative geometry.
\newblock {\em Communications in Mathematical Physics}, 164(3):525--562, aug
  1994.

\bibitem[KV06]{invitation}
J.~Kock and I.~Vainsencher.
\newblock {\em An Invitation to Quantum Cohomology: Kontsevich's Formula for
  Rational Plane Curves}.
\newblock Progress in Mathematics. Birkh{\"a}user Boston, 2006.

\bibitem[Pan97]{panda97}
Rahul Pandharipande.
\newblock A geometric construction of getzler's relation, 1997.

\bibitem[SR13]{shafarevich}
I.R. Shafarevich and M.~Reid.
\newblock {\em Basic Algebraic Geometry 1: Varieties in Projective Space}.
\newblock SpringerLink : B{\"u}cher. Springer Berlin Heidelberg, 2013.

\bibitem[Wit91]{Witten:1990hr}
Edward Witten.
\newblock {Two-dimensional gravity and intersection theory on moduli space}.
\newblock {\em Surveys Diff. Geom.}, 1:243--310, 1991.

\end{thebibliography}

\backmatter

\end{document}